\documentclass[11pt,reqno,a4paper]{amsart}

\usepackage[T1]{fontenc}
\usepackage{lmodern, comment}
\usepackage[a4paper,left=0.85in,right=0.85in,top=0.9in,bottom=0.9in]{geometry}
\usepackage{amsmath,amsthm,amssymb,mathtools,mathrsfs,enumitem}
\usepackage{xcolor}
\usepackage[expansion=false]{microtype}
\usepackage[hidelinks]{hyperref}
\renewcommand{\epsilon}{\varepsilon}
\newcommand{\Z}{\mathbb Z}

\newcommand{\ad}{\operatorname{ad}}
\newcommand{\B}{\mathbb B}

\newcommand{\D}{\mathbb D}

\newcommand{\R}{\mathbb R}
\newcommand{\C}{\mathbb C}

\newcommand{\Aut}{\operatorname{Aut}}
\newcommand{\Ric}{\operatorname{Ric}}
\newcommand{\Lie}{\mathcal L}
\newcommand{\tr}{\operatorname{tr}}
\newcommand{\ddbar}{\sqrt{-1}\,\partial\bar\partial}
\newcommand{\ii}{\sqrt{-1}}
\newcommand{\id}{\operatorname{id}}
\newcommand{\PSH}{\operatorname{PSH}}

\newcommand{\dist}{\operatorname{dist}}

\newcommand{\Isom}{\operatorname{Isom}}

\newtheorem{theorem}{Theorem}[section]
\newtheorem{proposition}[theorem]{Proposition}
\newtheorem{lemma}[theorem]{Lemma}
\newtheorem{corollary}[theorem]{Corollary}
\newtheorem{conjecture}[theorem]{Conjecture}

\theoremstyle{definition}
\newtheorem{remark}[theorem]{Remark}
\newtheorem{remarks}[theorem]{Remarks}

\title[Nontrivial expanders on Stein domains]{Nontrivial Expanding K\"ahler--Ricci Solitons and Automorphism Groups of Stein Domains}
\hypersetup{pdftitle={Nontrivial Expanding Kahler-Ricci Solitons and Automorphism Groups of Stein Domains},
  pdfauthor={Soumya Ganguly},
  pdfsubject={Complete expanding Kahler-Ricci solitons on Stein domains},
  pdfkeywords={Kahler-Ricci soliton, positive closed current, Stein domain, automorphism group, Kobayashi visibility}}

\author{Soumya Ganguly}
\address{Department of Mathematics, Rutgers University--New Brunswick,
110 Frelinghuysen Rd, Piscataway, NJ 08854, USA}
\email{soumya.ganguly@rutgers.edu}

\date{}

\begin{document}

\begin{abstract}
We study complete nontrivial expanding K\"ahler--Ricci
solitons with complete real holomorphic soliton fields.
Let $\Omega$ be a relatively compact domain in a complex
surface such that a neighborhood of $\overline\Omega$
admits a holomorphic immersion into Euclidean space.
We prove that such a soliton exists if and only if
$\Omega$ is Stein and the identity component of
$\Aut(\Omega)$ is noncompact. On such a Stein domain,
a complete real holomorphic vector field is the
soliton field of a complete nontrivial expander
if and only if its flow has noncompact closure
in $\Aut(\Omega)$. In higher dimensions, we establish existence criteria
using the adjoint representation of automorphism groups,
boundary geometry, and homogeneous positive closed
$(1,1)$-currents. We give sufficient conditions for the holomorphic
$\C^{n-1}$-action on the universal globalization of a
proper holomorphic $\R^{n-1}$-action on a Stein
$n$-manifold to be free and proper. Applications include bounded convex and
homogeneous domains, domains satisfying weak Kobayashi
visibility, and certain finite-type domains.
\end{abstract}

\maketitle

\section{Introduction and Main Results}\label{sec:introduction}

The K\"ahler--Ricci flow
\begin{equation}\label{eq:KRF}
\frac{\partial\omega(t)}{\partial t}=-\Ric_{\omega(t)}
\end{equation}
evolves a K\"ahler metric on a fixed complex manifold. Its
self-similar solutions, called K\"ahler--Ricci solitons, evolve
by rescaling and pullback under biholomorphisms. We study
complete expanding solitons whose real holomorphic soliton
vector fields are also complete.

The expanding case arises naturally on bounded pseudoconvex
domains. Mok--Yau~\cite[Main Theorem, p.~42]{MokYau} proved that
a bounded domain in \(\C^n\) is pseudoconvex if and only if it
admits a complete K\"ahler--Einstein metric with negative Einstein
constant. For any \(K>0\), a constant rescaling of this metric
gives a complete metric \(\omega_E\) satisfying
\begin{equation}\label{eq:KE-normalization}
\Ric_{\omega_E}=-K\omega_E,\qquad K>0.
\end{equation}
Then \(\omega(t)=(1+Kt)\omega_E\) solves \eqref{eq:KRF}. More
generally, an expanding K\"ahler--Ricci soliton satisfies
\begin{equation}\label{eq:expanding-soliton}
\Ric_\omega+K\omega=\Lie_Y\omega,\qquad K>0,
\end{equation}
where \(Y\) is a real holomorphic vector field and \(K\) is the
\emph{expansion constant}. Throughout this paper, the metric
\(\omega\) is required to be complete and the flow of \(Y\) to
be defined for all real time. The soliton is \emph{trivial} if
\(\Lie_Y\omega=0\), equivalently if
\(\Ric_\omega=-K\omega\), and \emph{nontrivial} otherwise.
Thus \((\omega_E,0)\) is a trivial expander.

Conlon--Deruelle--Sun
\cite[Theorem~A]{ConlonDeruelleSun2024}
study complete gradient K\"ahler--Ricci solitons
under quadratic curvature-decay assumptions,
including bounds on covariant derivatives of
curvature in the expanding case. Our focus is on
complete expanding solitons with complete real
holomorphic vector fields, without requiring the
gradient condition or curvature-decay assumptions.

In \cite[Theorem~1.1]{GangulyKRsol2026Paper21}, we proved that
for a relatively compact domain \(\Omega\) in a Stein manifold
\(Z\), the existence of a complete expanding soliton with complete
soliton field is equivalent to \(\Omega\) being Stein. The same
equivalence holds when \(Z\) is an arbitrary complex manifold,
provided a neighborhood of \(\overline\Omega\) admits a holomorphic
immersion into some \(\C^N\)
\cite[Remarks~3.6(2)]{GangulyKRsol2026Paper21}.

The existence of a complete expander does not by itself imply
the existence of a nontrivial one. Indeed, a complete real
holomorphic Killing field $Y$ of $\omega_E$ gives a trivial
expander $(\omega_E,Y)$, even when $Y\ne0$.

More generally, let $\{\Phi_s\}_{s\in\R}\subset\Aut_0(\Omega)$
be the flow of a complete real holomorphic vector field $Y$,
where $\Aut_0(\Omega)$ denotes the identity component of
$\Aut(\Omega)$. If this flow has compact closure in
$\Aut(\Omega)$, then every complete expanding soliton with
soliton field $Y$ is trivial
\cite[Lemma~3.1]{GangulyKRsol2026Paper21}.

Noncompactness of $\Aut_0(\Omega)$ is therefore necessary for
a complete nontrivial expander. We ask when it is sufficient
and, for a given complete real holomorphic vector field $Y$,
when an expander exists with soliton field $Y$.

\begin{conjecture}\label{conj:aut0}
Let \(Z\) be a complex manifold and let \(\Omega\Subset Z\) be a
domain of complex dimension \(n\ge2\). Suppose that an open
neighborhood \(V\subset Z\) of \(\overline\Omega\) admits a
holomorphic immersion \(F:V\to\C^N\). Then \(\Omega\) admits
a complete nontrivial expanding K\"ahler--Ricci soliton with
complete field if and only if \(\Omega\) is Stein and
\(\Aut_0(\Omega)\) is noncompact.
\end{conjecture}

We prove Conjecture~\ref{conj:aut0} in complex dimension
two and identify complete real holomorphic fields that occur
as soliton fields in higher dimensions. The constructions rely on
positive closed $(1,1)$-currents that transform
exponentially under the flow of a complete real
holomorphic vector field.

Let \((\omega,Y)\) satisfy \eqref{eq:expanding-soliton}, and
write \(\Phi_s\) for the flow of \(Y\). The associated
self-similar K\"ahler--Ricci flow, normalized at \(t=1/K\),
is
\begin{equation}\label{eq:selfsimilar-flow}
\omega(t)=Kt\,\Phi_{-K^{-1}\log(Kt)}^*\omega,\qquad t>0
\end{equation}
\cite[equation~(4)]{GangulyKRsol2026Paper21}. By the group
law \eqref{eq:real-action-group-law},
\begin{equation}\label{eq:selfsimilar-homogeneity}
\Phi_s^*\omega(t)
 =Kt\,\Phi_{s-K^{-1}\log(Kt)}^*\omega
 =e^{Ks}\omega(e^{-Ks}t).
\end{equation}
If \(\omega(t)\rightharpoonup T\) as \(t\downarrow0\), then
\(T\) is a positive closed \((1,1)\)-current and
\eqref{eq:selfsimilar-homogeneity} gives
\begin{equation}\label{eq:current-homogeneity}
\Phi_s^*T=e^{Ks}T,\qquad s\in\R.
\end{equation}
Let $\mathcal T_{Y,K}(\Omega)$ denote the cone of positive closed
real $(1,1)$-currents satisfying \eqref{eq:current-homogeneity}.
If $\Omega\Subset Z$ is Stein and a neighborhood of
$\overline\Omega$ admits a holomorphic immersion into
$\C^N$, then a nonzero element of this cone yields
a complete nontrivial expander with field $Y$ and
expansion constant $K$
\cite[Theorem~1.2]{GangulyKRsol2026Paper21}.

We use the adjoint representation of $\Aut_0(\Omega)$ to construct
homogeneous currents. With the usual bracket of vector fields,
$\ad_Y$ denotes the complex-linear extension of $V\mapsto[Y,V]$.

\begin{theorem}\label{stein:thm:ambient-unified}\label{thm:unified-current}
Let \(\Omega\) be a connected complex manifold of dimension \(n\ge2\)
whose bounded holomorphic functions separate tangent directions.
Suppose that \(G=\Aut_0(\Omega)\) is a noncompact Lie group acting
properly, and let \(\mathfrak g\) be its Lie algebra of infinitesimal
automorphisms.
\begin{enumerate}[label=\textup{(\roman*)},leftmargin=0pt,
                  labelindent=0pt,labelsep=0.45em,itemindent=*,
                  widest=ii,align=left]
\item Let \(Y\in\mathfrak g\) generate a proper real flow.
If \(\ad_Y\) has an eigenvalue with positive (respectively negative)
real part, then \(\mathcal T_{Y,K}(\Omega)\)
(respectively \(\mathcal T_{-Y,K}(\Omega)\)) contains a nonzero
current for every \(K>0\). If \(\ad_Y\) has a nontrivial Jordan
block corresponding to a purely imaginary eigenvalue, both
\(\mathcal T_{Y,K}(\Omega)\) and \(\mathcal T_{-Y,K}(\Omega)\)
contain nonzero smooth forms for every \(K>0\).
\item Suppose that for every proper real flow generated by
\(Y\in\mathfrak g\), the operator \(\ad_Y\) is diagonalizable over
\(\C\) with spectrum in \(\ii\R\). Then
\begin{equation}\label{eq:central-decomposition}
G\simeq\R^p\times C,\qquad 1\le p\le n,
\end{equation}
with \(C\) a compact connected Lie group. The real factor is central
and acts freely and properly, and the holomorphic parts of its
generators are pointwise complex-linearly independent.
If \(\Omega\) is Stein, its Bergman kernel form is nowhere zero,
its Bergman metric is positive definite, and \(p\ge n-1\), then for
every \(Y\in\mathfrak g\) with nonzero \(\R^p\)-component
and every \(K>0\), both \(\mathcal T_{Y,K}(\Omega)\) and
\(\mathcal T_{-Y,K}(\Omega)\) contain nonzero smooth forms.
\end{enumerate}

Suppose in addition that $\Omega\Subset Z$ is Stein and that a
neighborhood of $\overline\Omega$ admits a holomorphic immersion
into Euclidean space. For every $Y\in\mathfrak g$ and $K>0$ for
which \textup{(i)} or \textup{(ii)} establishes
$\mathcal T_{Y,K}(\Omega)\ne\{0\}$, there is a complete
nontrivial expander with soliton field $Y$ and expansion constant
$K$ \cite[Theorem~1.2]{GangulyKRsol2026Paper21}.
\end{theorem}

Theorem~\ref{thm:unified-current} yields the following existence criterion.

\begin{corollary}\label{cor:unified-existence}
Let $\Omega\Subset Z$ be a Stein domain of complex
dimension $n\ge2$ such that a neighborhood of
$\overline\Omega$ admits a holomorphic immersion into
$\C^N$. Suppose that $G=\Aut_0(\Omega)$ is noncompact
and is not isomorphic to
\[
\R^p\times C,\qquad 1\le p\le n-2,
\]
for any compact connected Lie group $C$.
Then, for every $K>0$, $\Omega$ admits a complete
nontrivial expanding K\"ahler--Ricci soliton with
complete soliton vector field and expansion constant $K$.
\end{corollary}

For a Stein domain in a complex surface satisfying
the immersion assumption, Lemma~\ref{stein:lem:realflow}
shows that every complete real holomorphic field
with noncompact flow closure generates a proper
$\R$-action. Applying
Theorem~\ref{thm:central-corank-one}
gives the following characterization.

\begin{corollary}\label{thm:stein-surface-intro}\label{thm:surface-prescribed}
\label{stein:cor:dimtwo}\label{thm:two-dimensional-characterization}
Conjecture~\ref{conj:aut0} holds in complex dimension two.
Moreover, let \(\Omega\) be a two-dimensional Stein domain satisfying
the hypotheses of the conjecture, and let \(Y\) be a complete real
holomorphic vector field on \(\Omega\). Then \(Y\) is the
soliton field of a complete nontrivial expander if and only
if its flow has noncompact closure in \(\Aut(\Omega)\).
For every such \(Y\) and every \(K>0\), there exist
smooth nonzero positive closed real \((1,1)\)-forms
\(T_+\) and \(T_-\), each of rank one at every point,
such that
\[
\Phi_s^*T_+=e^{Ks}T_+,
\qquad
\Phi_s^*T_-=e^{-Ks}T_-.
\]
In particular, both \(Y\) and \(-Y\) occur as soliton
fields of complete nontrivial expanders with expansion
constant \(K\).
\end{corollary}

Thus, by Corollary~\ref{thm:stein-surface-intro}
and \cite[Remarks~3.6(1)]{GangulyKRsol2026Paper21},
a bounded domain in $\C^2$ admits a complete
nontrivial expander with complete soliton field
if and only if it is pseudoconvex and
$\Aut_0(\Omega)$ is noncompact.

Yau conjectured that the Bergman metric of a bounded
pseudoconvex domain is complete and K\"ahler--Einstein if and only
if the domain is homogeneous \cite{Yau1982Problem}.
Our characterization by the existence of complete
nontrivial expanding K\"ahler--Ricci solitons gives
a related classification in complex dimension two.
In \cite[Theorem~1.4]{GangulyKRsol2026Paper21},
we proved that a relatively compact domain of complex
dimension at least two, whose bounded holomorphic
functions separate points and whose boundary orbit
accumulation points are all $C^2$ strongly pseudoconvex,
is biholomorphic to the unit ball if and only if it
admits a complete nontrivial expanding K\"ahler--Ricci
soliton with complete soliton field. For smoothly
bounded strongly pseudoconvex domains in $\C^n$,
$n\ge2$, Huang and Xiao characterized the ball
by the K\"ahler--Einstein property of the Bergman
metric \cite[Theorem~1.1]{HuangXiao2016}.

We also construct homogeneous currents from bounded
nonnegative plurisubharmonic functions that decay locally
uniformly along the backward flow.

\begin{theorem}\label{thm:backward-psh-prescribed-intro}
Let \(\Omega\Subset Z\) be Stein, and suppose that a neighborhood
of $\overline\Omega$ admits a holomorphic immersion into $\C^N$.
Let \(Y\) generate a complete real holomorphic flow. Suppose a bounded nonnegative
\(u\in\PSH(\Omega)\), \(u\not\equiv0\), satisfies
\begin{equation}\label{eq:backward-psh-decay}
 u\circ\Phi_s\longrightarrow0
 \quad\text{locally uniformly as }s\to-\infty .
\end{equation}
Then for every \(K>0\) there is a complete nontrivial expander with
field \(Y\) and expansion constant \(K\).
\end{theorem}

The next theorem uses the time-one displacement of a holomorphic
immersion $\iota=(f_1,\ldots,f_N)$ defined near $\overline\Omega$.
Under either of the following boundary hypotheses, every proper real holomorphic
flow $\Phi_t$ satisfies
\[
(f_j\circ\Phi_1-f_j)\circ\Phi_t\longrightarrow0
\quad\text{locally uniformly as }t\to\pm\infty
\]
for each $j$; see \eqref{eq:unit-displacement-decay}.
Properness and the immersion ensure that at least one
displacement is not identically zero, so
Lemma~\ref{lem:two-ended-holomorphic} applies.

\begin{theorem}\label{thm:intro-manifold-geometric}
\label{thm:stein-no-boundary-discs}\label{prop:disc-free-current}
\label{thm:disc-free-prescribed}\label{thm:stein-weak-visibility}
\label{prop:visibility-current}\label{thm:visibility-prescribed}
\label{thm:disc-free-main}\label{cor:stein-no-discs-aut}
Let \(\Omega\Subset Z\) be a domain such that a neighborhood
of $\overline\Omega$ admits a holomorphic immersion into $\C^N$.
Assume either
\begin{enumerate}[label=\textup{(\roman*)}]
\item \(\partial_Z\Omega\) contains no nonconstant analytic disc, or
\item $\Omega$ has $(1,\kappa_0)$-visibility relative to
$\partial_Z\Omega$ for some $\kappa_0>0$.
\end{enumerate}
For every proper holomorphic real flow on \(\Omega\), with generator
\(Y\), and every \(K>0\), both \(\mathcal T_{Y,K}(\Omega)\) and
\(\mathcal T_{-Y,K}(\Omega)\) contain nonzero smooth forms.
If \(\Omega\) is Stein, both \(Y\) and \(-Y\) are soliton fields
of complete nontrivial expanders with expansion constant \(K\).
A complete nontrivial expander with complete soliton field exists
if and only if \(\Omega\) is Stein and
\(\Aut_0(\Omega)\) is noncompact.
\end{theorem}

Suppose that $\Omega$ is Kobayashi hyperbolic and satisfies
relative $(1,\kappa_0)$-visibility for some $\kappa_0>0$. If
$F\in\Aut(\Omega)$ and each pair $E,L\Subset\Omega$ of compact sets
satisfies $F^m(E)\cap L=\varnothing$ for all sufficiently large
$m$, then $F^m$ converges locally uniformly to one point of
$\partial_Z\Omega$ (Lemma~\ref{lem:weak-visibility-iterates}).
For a proper real flow this gives a limit at each end of real
time (Remark~\ref{rem:single-kappa-manifold}). For the iteration
theorem on bounded taut weakly visible submanifolds, see
\cite[Theorem~1.14]{ChandelMaitraSarkar2024}.
Bounded Goldilocks domains satisfy the visibility hypothesis
by \cite[Theorem~1.4]{BharaliZimmer2017}.

A nonconstant analytic disc in the boundary has infinite order of contact
with the boundary along its image. In particular, finite
D'Angelo $1$-type at every boundary point rules out nonconstant
boundary discs.

\begin{corollary}\label{cor:intro-global-finite-type}
\label{ft:thm:higher-global}\label{cor:finite-type-higher}
Let \(\Omega\Subset Z\) be a Stein domain such that a neighborhood
of $\overline\Omega$ admits a holomorphic immersion into $\C^N$.
Suppose \(\partial_Z\Omega\)
is \(C^\infty\) and of finite D'Angelo \(1\)-type at every point.
A complete real holomorphic field has noncompact flow closure if and
only if it is the field of a complete nontrivial expander for
every \(K>0\).
\end{corollary}

Suppose that $\Psi_t$ is a proper holomorphic $\R$-action
on $M$ and that $M/\langle\Psi_1\rangle$ is Stein.
Proposition~\ref{prop:stein-cyclic-character} gives a
holomorphic function $F$ on $M$ satisfying
$F\circ\Psi_t=F+t$. Hence the forms
$\ddbar e^{\pm K\Re F}$ transform by the characters
$e^{\pm Kt}$.

\begin{theorem}\label{thm:equivariant-stein-target}
Let $\{\Phi_t\}_{t\in\R}$ be a holomorphic $\R$-action on a
complex manifold $\Omega$, with generator $Y$. Let
$\{\Psi_t\}_{t\in\R}$ be a proper holomorphic $\R$-action
on a connected complex manifold $M$. Suppose there is a holomorphic
map $f:\Omega\to M$ satisfying
\begin{equation}\label{eq:flow-map-relation}
f\circ\Phi_t=\Psi_t\circ f,\qquad t\in\R.
\end{equation}
If $M/\langle\Psi_1\rangle$ is Stein, then for every $K>0$
there exist smooth positive closed real $(1,1)$-forms
$T_+$ and $T_-$ on $\Omega$, each of rank one at
every point, satisfying
\[
\Phi_t^*T_+=e^{Kt}T_+,
\qquad \Phi_t^*T_-=e^{-Kt}T_-.
\]
If $\Omega$ is a relatively compact Stein domain such that
a neighborhood of its closure admits a holomorphic immersion
into $\C^N$, then both $Y$
and $-Y$ are soliton fields of complete nontrivial expanders
with expansion constant $K$.
\end{theorem}

For a proper flow on a bounded homogeneous domain,
the time-one quotient is Stein by
\cite[Theorem~6.5]{MiebachHomogeneous2010}.
Theorem~\ref{thm:equivariant-stein-target} therefore
applies with the identity map. For bounded convex
domains, we instead use the estimate for the diagonal
Bergman kernel in Lemma~\ref{lem:convex-kernel}
to control the logarithmic Jacobians of automorphisms,
and then apply Proposition~\ref{prop:frame-halfline}.

\begin{theorem}\label{thm:geometric-prescribed}
Let $M$ be biholomorphic to a bounded convex or bounded homogeneous
domain in $\C^n$. For every proper real holomorphic flow $\Phi_t$ on
$M$ with generator $Y$, and every $K>0$, there are nonzero smooth
positive closed real $(1,1)$-forms $T_\pm$ with
\[
\Phi_t^*T_\pm=e^{\pm Kt}T_\pm .
\]
Both $Y$ and $-Y$ are soliton fields of complete nontrivial expanders
with expansion constant $K$. If $n\ge2$, $M$ admits a complete
nontrivial expander with complete soliton field if and only if
$\Aut_0(M)$ is noncompact.
\end{theorem}

In the homogeneous case, the nontrivial soliton metrics
in Theorem~\ref{thm:geometric-prescribed} cannot themselves
be homogeneous. Indeed, by
\cite[Theorem~1.1(i) and the proof of
Proposition~4.1]{LoiMossa2023}, a
K\"ahler--Ricci soliton admitting a holomorphic isometric
immersion into a bounded homogeneous domain with a
homogeneous K\"ahler metric is K\"ahler--Einstein,
and its soliton field is Killing.

Under the finite-type and Levi-rank hypotheses below,
Bedford--Pinchuk \cite[main theorem, p.~165]{BedfordPinchukEgg1991}
classify bounded pseudoconvex domains with noncompact
automorphism group. The real translations on a Siegel
realization of each egg give the converse construction.

\begin{corollary}\label{eggdomain_char}
\label{cor:finite-type-bp-all-dim}
\label{pcvx:prop:bedford-pinchuk}
Let $\Omega\Subset\C^{n+1}$, $n\ge1$, have
$C^\infty$ boundary of finite D'Angelo type and
Levi rank at least $n-1$ everywhere. Then $\Omega$
admits a complete nontrivial expanding K\"ahler--Ricci
soliton with complete soliton field if and only if it
is biholomorphic to
\begin{equation}\label{eq:egg-model-family}
E_m=\left\{(w,z_1,\ldots,z_n):
|w|^2+|z_1|^{2m}+\sum_{j=2}^n|z_j|^2<1
\right\}
\end{equation}
for some integer $m\ge1$.

In this case, if $Y$ generates a complete real holomorphic flow
with noncompact closure, then for every $K>0$ there are nonzero
smooth positive closed $(1,1)$-forms in
$\mathcal T_{Y,K}(\Omega)$ and $\mathcal T_{-Y,K}(\Omega)$.
Both fields occur for complete nontrivial expanders of expansion
constant $K$.
\end{corollary}

For a domain in a complex surface with a smooth pseudoconvex
finite-type boundary orbit accumulation point, the polynomial
model of Seshadri--Verma
\cite[Theorem~1.1]{SeshadriVerma2008}
is pseudoconvex, and hence the domain is Stein. If, moreover, $\Omega\Subset Z$ and a neighborhood of
$\overline\Omega$ admits a holomorphic immersion into $\C^N$, then Theorem~\ref{ft:thm:two-d-prescribed}
gives nonzero smooth positive closed $(1,1)$-forms
of rank one in $\mathcal T_{Y,K}(\Omega)$ and
$\mathcal T_{-Y,K}(\Omega)$ for every complete
real holomorphic field $Y$ with noncompact
flow closure and every $K>0$. 

On a weighted polynomial model, we give explicit currents for
every nonzero real linear combination of dilation and translation
fields (Proposition~\ref{ft:prop:weighted}). This applies to
bounded pseudoconvex domains in $\C^n$ with smooth finite-type
boundary and at least two points in the boundary orbit
accumulation set, by the uniformization theorem of Zimmer
\cite[Theorem~1.2]{Zimmer2020}. Section~\ref{sec:finite-type-extensions}
also treats nontangential orbit accumulation at an $h$-extendible
point.

Section~\ref{sec:preliminaries} fixes notation and conventions.
Section~\ref{sec:positive-current-methods} constructs homogeneous
positive closed currents and establishes the estimates needed
for Fourier integrals. Section~\ref{sec:central-reduction}
proves properness of the automorphism-group action under the
immersion hypothesis and studies the adjoint representation.
Section~\ref{sec:corank-complexification} proves Hausdorffness
of the universal globalization, using the moment map associated
with a complete K\"ahler--Einstein metric of negative Einstein
constant or the positive definite Bergman metric.
Fourier transforms of tensor squares of square-integrable
holomorphic $n$-forms then give equivariant holomorphic sections
whose ratios establish properness of the complexified action
in corank one (Theorem~\ref{stein:thm:corank}).
The full-rank case gives a translation action on $\C^n$
(Corollary~\ref{stein:cor:unconditional-corank}).
The associated character line bundles yield the rank-one
homogeneous forms of Theorem~\ref{thm:central-corank-one}.
Section~\ref{sec:boundary-dynamics} studies normal limits of
automorphisms and boundary convergence of iterates under
relative Kobayashi visibility. Section~\ref{sec:geometric-applications}
applies these results to construct homogeneous currents and
complete expanders under geometric hypotheses.
Finally, Section~\ref{sec:finite-type-extensions} applies
uniformization theorems and constructs explicit homogeneous
currents on weighted polynomial models.

\section{Conventions and Preliminaries}\label{sec:preliminaries}

Throughout the paper, a domain is a connected open set.
For a subset $A$ of a complex manifold $Z$ we write $A\Subset Z$
if $\overline A$ is compact in $Z$. For a domain $\Omega\subset Z$ put
\begin{equation}\label{eq:relative-boundary-notation}
\partial_Z\Omega=\overline\Omega\setminus\Omega.
\end{equation}
We put $n=\dim_\C\Omega=\dim_\C Z$. We write
$\Omega\Subset\C^n$ for a bounded Euclidean domain. A point $\xi\in\partial_Z\Omega$ is a \emph{boundary orbit
accumulation point} if there exist $z_0\in\Omega$ and
$F_j\in\Aut(\Omega)$ such that $F_j(z_0)\to\xi$ in $Z$.

For a K\"ahler manifold $(M,J,g)$ we write
$\omega_g=g(J\cdot,\cdot)$ and, when no confusion is likely, use
$\omega$ for both the K\"ahler metric and its associated form. Our
local conventions are
\begin{equation}\label{eq:metric-ricci-conventions}
\omega=\ii g_{j\bar k}\,dz^j\wedge d\bar z^k,
\qquad
\Ric_\omega=-\ddbar\log\det(g_{j\bar k}).
\end{equation}
The Riemannian Ricci tensor satisfies
$\Ric_\omega(V,W)=\Ric_g(JV,W)$ for real tangent vectors $V,W$.
On $\C^n$ set
\begin{equation}\label{eq:euclidean-volume-convention}
\omega_{\mathrm{Euc}}=\ii\sum_{j=1}^n dz^j\wedge d\bar z^j,
\qquad dV_{\mathrm E}=\omega_{\mathrm{Euc}}^n.
\end{equation}
In the coordinates of \eqref{eq:metric-ricci-conventions},
$\omega^n/dV_{\mathrm E}=\det(g_{j\bar k})$. For a real $(1,1)$-form
$\alpha$ we set
\begin{equation}\label{eq:trace-laplacian-conventions}
\tr_\omega\alpha
=n\frac{\alpha\wedge\omega^{n-1}}{\omega^n},
\qquad
\Delta_\omega f=\tr_\omega(\ddbar f)
=g^{j\bar k}f_{j\bar k},
\qquad
R_\omega=\tr_\omega\Ric_\omega.
\end{equation}
If $\Delta_g=\tr_g\nabla^2$ and $S_g$ are the Riemannian Laplacian
and scalar curvature, respectively, then
$\Delta_g=2\Delta_\omega$ and $S_g=2R_\omega$. For real $f$,
$|\nabla^g f|_g^2=2|\partial f|_\omega^2$; the Riemannian volume
form is $dV_g=\omega^n/n!$. Finally,
\begin{equation}\label{eq:ddc-convention}
d^c=\ii(\bar\partial-\partial),\qquad dd^c=2\ddbar.
\end{equation}

For a real vector field $Y$ with local flow $\{\Phi_t\}$, our Lie
derivative convention is
\begin{equation}\label{eq:lie-derivative-convention}
\Lie_Y A
=\left.\frac{d}{dt}\right|_{t=0}\Phi_t^*A.
\end{equation}
The pullback of a vector field $V$ means
$\Phi_t^*V=(\Phi_{-t})_*V$; in particular,
$\Lie_YV=[Y,V]$. The field $Y$ is \emph{real holomorphic} if
$\Lie_YJ=0$, equivalently if its local flow consists of
biholomorphisms. It is \emph{complete} if this flow is defined for
all real time. Completeness of a K\"ahler metric always means
geodesic completeness.

A holomorphic $\R$-action on $M$ is a smooth map
$\Phi:\R\times M\to M$ whose time-$t$ maps $\Phi_t$ are
holomorphic and satisfy
\begin{equation}\label{eq:real-action-group-law}
\Phi_0=\id_M,\qquad
\Phi_{s+t}=\Phi_s\circ\Phi_t\quad(s,t\in\R).
\end{equation}
Thus $\Phi_t^{-1}=\Phi_{-t}$, and its infinitesimal generator is a
complete real holomorphic vector field. Conversely, each such field
generates a holomorphic $\R$-action. The action is \emph{free} if
$\Phi_t(z)=z$ implies $t=0$; it is \emph{proper} if
\begin{equation}\label{eq:proper-action-convention}
\R\times M\longrightarrow M\times M,
\qquad (t,z)\longmapsto (z,\Phi_t(z))
\end{equation}
is proper. The same definition applies to a Lie group $H$ acting on
$M$, with $H$ in place of $\R$ and $(h,z)\mapsto(z,h\cdot z)$ in
\eqref{eq:proper-action-convention}.

We write $\Aut(M)$ for the group of biholomorphisms of $M$ equipped
with the compact-open topology, and $\Aut_0(M)$ for its identity
component. When $G=\Aut_0(\Omega)$ is a Lie group, $\mathfrak g$
denotes its real Lie algebra of complete real holomorphic vector
fields, with the usual vector-field bracket, as in
\eqref{eq:lie-derivative-convention}. We write
$\mathfrak g_\C=\mathfrak g\otimes_\R\C$, and for
$Y\in\mathfrak g$ set
\begin{equation}\label{eq:adjoint-and-holomorphic-part}
\ad_Y(V)=[Y,V],\qquad
Y^{1,0}=\tfrac12(Y-\ii JY).
\end{equation}
We say that bounded holomorphic functions on a complex
manifold $M$ \emph{separate tangent directions} if,
for every $z\in M$ and $0\ne v\in T_z^{1,0}M$,
there exists a bounded $f\in\mathcal O(M)$ such that
$df_z(v)\ne0$.
The flow has \emph{noncompact closure} if
$\overline{\{\Phi_t:t\in\R\}}$ is noncompact in $\Aut(M)$.
For a holomorphic $\R^p$-action the parameter
$t=(t_1,\ldots,t_p)$ ranges over $\R^p$ in the same smooth action
definition. The generator of $s\mapsto\Phi_{sv}$ for $v\in\R^p$
is $X_v$. The universal complexification of $A=\R^p$ is $\C^p$.

For a complex manifold $M$, $\PSH(M)$ denotes its
plurisubharmonic functions, with the value $-\infty$ allowed.
A domain $\Omega\subset\C^n$ is pseudoconvex if it has a
continuous plurisubharmonic exhaustion $\rho$; thus
\begin{equation}\label{eq:psh-exhaustion-convention}
\{\rho<c\}\Subset\Omega\qquad(c\in\R).
\end{equation}

In Euclidean space this is equivalent to Steinness
\cite{Hormander1990}.
A smooth real function $h$ is pluriharmonic if $\ddbar h=0$.

For a domain $\Omega\Subset Z$ with $C^\infty$ boundary, let
$p\in\partial_Z\Omega$ and choose local holomorphic coordinates
centered at $p$ and a smooth real defining function $\rho$ near $p$.
The D'Angelo $1$-type of the boundary at $p$ is
\[
\Delta_1(\partial_Z\Omega,p)
=\sup_\gamma\frac{\operatorname{ord}_0(\rho\circ\gamma)}
{\operatorname{ord}_0(\gamma-\gamma(0))},
\]
where $\gamma$ ranges over nonconstant holomorphic curve germs
with $\gamma(0)=p$. For a smooth function germ, the order is
the least degree of a nonzero Taylor term, or $+\infty$ if
all Taylor coefficients vanish.
A change of defining function multiplies $\rho$ by a
nowhere-vanishing smooth function, while a change of
holomorphic coordinates preserves the vanishing order
of the curve germ. Thus $\Delta_1(\partial_Z\Omega,p)$
is independent of these choices
\cite[Definitions~2.16 and~2.18]{DAngelorealhypordcontc};
see also \cite[Section~6]{DAngeloKohn1999}.

A real $(1,1)$-current $T$ on an $n$-dimensional complex manifold is
\emph{positive} if it takes nonnegative values on smooth,
compactly supported strongly positive $(n-1,n-1)$-forms. It is
\emph{closed} if $dT=0$ distributionally. A positive closed
$(1,1)$-current has local plurisubharmonic potentials
$T=\ddbar u$. Weak convergence $T_j\rightharpoonup T$ is
convergence on smooth, compactly supported $(n-1,n-1)$-forms.
For a smooth real
$(1,1)$-form, positivity means semipositivity. The rank
of a smooth semipositive $(1,1)$-form at a point is the
rank over $\C$ of its Hermitian coefficient matrix in
local holomorphic coordinates.

For a complete real holomorphic field $Y$ on $M$ with flow $\Phi_s$
and $K>0$, define the cone
\begin{equation}\label{eq:homogeneous-current-cone}
\mathcal T_{Y,K}(M)=\Bigl\{T:\
\begin{aligned}&T\text{ is a positive closed real }(1,1)\text{-current on }M,\\
&\Phi_s^*T=e^{Ks}T\quad(s\in\R)\end{aligned}\Bigr\}.
\end{equation}
More generally, let $G$ act holomorphically on $M$. A tensor
field or current $A$ on $M$ is said to transform by a
\emph{character} $\chi:G\to\C^*$ if $\chi$ is a continuous
group homomorphism and
\[
\Phi_g^*A=\chi(g)A\qquad(g\in G).
\]
In particular, every current in $\mathcal T_{Y,K}(M)$
is homogeneous with respect to $\Phi_s$, with positive
character $\chi(s)=e^{Ks}$. The cone
\eqref{eq:homogeneous-current-cone} may consist only of
the zero current. Reversing the flow gives
$\mathcal T_{-Y,K}(M)$, whose elements satisfy
$\Phi_s^*T=e^{-Ks}T$.

We denote by $k_M$ and $\kappa_M$ the Kobayashi pseudodistance
and the infinitesimal Kobayashi--Royden pseudometric, respectively.
Thus $k_M$ is defined by chains of holomorphic discs and
$M$ is Kobayashi hyperbolic if $k_M$ is a distance.
For a Hermitian metric $h$ on an ambient manifold $Z$,
$d_h$ is its induced distance.
Our normalization on the disc is
\[
\kappa_{\D}(\zeta;v)=\frac{|v|}{1-|\zeta|^2},
\qquad
k_{\D}(0,r)=\operatorname{arctanh} r \quad (0\le r<1).
\]
For $\lambda\ge1$ and $\kappa\ge0$, a
$(\lambda,\kappa)$-almost-geodesic in $\Omega$ is an absolutely
continuous path $\gamma:[a,b]\to\Omega$ with
$\kappa_\Omega(\gamma(t);\gamma'(t))\le\lambda$ almost
everywhere and
\begin{equation}\label{eq:almost-geodesic-convention}
\lambda^{-1}|s-t|-\kappa
\le k_\Omega(\gamma(s),\gamma(t))\le\lambda|s-t|+\kappa
\qquad(a\le s,t\le b)
\end{equation}
\cite[Definition~2.3]{Masanta2025Visibility}.

Following \cite[Definition~1.1(c)]{Masanta2025Visibility}, a
Kobayashi hyperbolic domain $\Omega\Subset Z$ has
\emph{relative weak visibility} if, for every pair of distinct
points $\xi,\eta\in\partial_Z\Omega$, there exist neighborhoods
$U_\xi,U_\eta$ in $Z$ with disjoint closures such that, for every
$\kappa\ge0$, there is a compact set $Q\Subset\Omega$ meeting
every $(1,\kappa)$-almost-geodesic $\gamma:[a,b]\to\Omega$ with
\[
\gamma(a)\in U_\xi\cap\Omega,\qquad
\gamma(b)\in U_\eta\cap\Omega.
\]

For a fixed $\kappa_0>0$, relative $(1,\kappa_0)$-visibility
means that the same condition holds with $\kappa=\kappa_0$ only.
For a bounded domain $\Omega\Subset\C^n$, this agrees with
the $(1,\kappa_0)$-visibility condition of
\cite[Definition~1.1 and Theorem~2.8]{ChandelMaitraSarkar2024}.

A sequence of maps $F_j:\Omega\to\Omega$ is
\emph{compactly divergent} if, for every pair of compact sets
$E,L\Subset\Omega$,
\[
F_j(E)\cap L=\varnothing
\]
for all sufficiently large $j$. An \emph{ambient normal limit}
of such a sequence is a map $F:\Omega\to Z$ obtained as a
locally uniform limit in $Z$ of a subsequence.

For $n=\dim_\C M$, let $K_M=\Lambda^n(T^{1,0}M)^*$
be the canonical bundle and $K_M^{\otimes2}$ its tensor
square. A holomorphic frame of $K_M$ over $U\subset M$
is a nowhere-vanishing holomorphic $n$-form on $U$,
such as $dz^1\wedge\cdots\wedge dz^n$ in local coordinates.
A global frame is one defined on all of $M$.
We write $\mathcal O(M)$ for holomorphic functions on $M$
and $H^0(M,L)$ for holomorphic sections of a bundle $L$.
For $\sigma\in H^0(M,K_M)$ set
\begin{equation}\label{eq:canonical-volume}
dV_\sigma=\ii^{n^2}\sigma\wedge\overline\sigma,\qquad
\|\sigma\|_{2,M}^2=\int_M dV_\sigma .
\end{equation}
If $\sigma=\rho\,dz^1\wedge\cdots\wedge dz^n$, then
$dV_\sigma=|\rho|^2dV_{\mathrm E}/n!$ by
\eqref{eq:euclidean-volume-convention}.
The Bergman space $A^2(M,K_M)$ consists of
holomorphic $n$-forms with finite norm. Its Bergman kernel form on the diagonal
is $\mathcal B_M=\sum_j dV_{\sigma_j}$ for an orthonormal
basis $\{\sigma_j\}$ of $A^2(M,K_M)$.
In a local holomorphic frame $e$ write
$\mathcal B_M=B_e\,dV_e$. Equivalently $B_e(z)$ is the
squared norm of the evaluation functional in that frame.
The Bergman kernel form $\mathcal B_M$ is nonzero at $z$ exactly
when some $\sigma\in A^2(M,K_M)$ satisfies $\sigma(z)\ne0$.
On this set the Bergman $(1,1)$-form is
\begin{equation}\label{eq:bergman-form}
\omega_B=\ddbar\log B_e .
\end{equation}
Since $B_{ge}=|g|^{-2}B_e$ for a holomorphic frame change
$e\mapsto ge$, \eqref{eq:bergman-form} is frame independent.
For a Euclidean domain, $K_\Omega(z,w)$ is the scalar
Bergman kernel for functions square-integrable with respect
to $dV_{\mathrm E}$, so
$\mathcal B_\Omega=K_\Omega(z,z)\,dV_{\mathrm E}$.

We write $\omega_E$ for a complete K\"ahler--Einstein form,
normalized by \(\Ric_{\omega_E}=-\omega_E\) unless a different
Einstein constant is specified.
\section{Construction of Positive Closed Currents}
\label{sec:positive-current-methods}

\subsection{Plurisubharmonic constructions}

We construct a plurisubharmonic potential with character
$e^{Ks/2}$. Its square has character $e^{Ks}$.

\begin{lemma}\label{lem:positive-character}
Let $M$ be a connected complex manifold, and let
$\{\Phi_s\}_{s\in\R}$ be a holomorphic $\R$-action
on $M$. Suppose that
$U\in\PSH(M)\cap L^\infty_{\mathrm{loc}}(M)$
is nonnegative and not identically zero, and that
for some $a>0$,
\begin{equation}\label{eq:positive-character-law}
U\circ\Phi_s=e^{as}U
\qquad(s\in\R).
\end{equation}
Then
\[
T=\ddbar(U^2)
\]
is a nonzero positive closed real $(1,1)$-current
on $M$ satisfying
\[
\Phi_s^*T=e^{2as}T
\qquad(s\in\R).
\]
If, in addition, $\log U\in\PSH(M)$, with
$\log 0=-\infty$, then
$S=\ddbar U$ is a nonzero positive closed real
$(1,1)$-current satisfying
\[
\Phi_s^*S=e^{as}S
\qquad(s\in\R).
\]
If \eqref{eq:positive-character-law} is assumed only
for $s=1$, the conclusions hold with the character
identities restricted to $s\in\Z$.
\end{lemma}

\begin{proof}
Since $U$ is nonnegative and plurisubharmonic,
and the function $t\mapsto t^2$ is convex and
increasing on $[0,\infty)$, the function $U^2$
is plurisubharmonic. Consequently,
$T=\ddbar(U^2)$ is a positive closed real
$(1,1)$-current.

We show that $T$ is nonzero. Suppose that
$T=0$, and set $H=U^2$. Then
$\ddbar H=0$ in the distributional sense.
Taking the Euclidean trace in a coordinate ball gives
$\Delta H=0$ in the distributional sense.
By Weyl's lemma \cite[Theorem~IX.25]{ReedSimonII1975}, $H$ has a smooth harmonic representative.
Since $H$ is plurisubharmonic, it agrees everywhere
with this representative. Hence $H$ is smooth
and pluriharmonic.

We claim that $H>0$ on $M$. Indeed, if $H(p)=0$
at some point $p\in M$, the strong minimum principle
for harmonic functions
\cite[Section~2.2.3]{PDEbookEvans_1998}
implies that $H$ vanishes in a neighborhood of $p$.
Thus $\{H=0\}$ is both open and closed in $M$.
By connectedness, $H$ would vanish identically,
contrary to the hypothesis on $U$.
Therefore $H>0$ and $U=\sqrt H$ is smooth. Since $\ddbar H=0$, differentiation gives
\[
\ddbar U
=
-\frac{\ii}{4H^{3/2}}\,
\partial H\wedge\bar\partial H.
\]
The left-hand side is nonnegative because
$U$ is plurisubharmonic, whereas the
right-hand side is nonpositive. Therefore
\[
\partial H=0.
\]
Since $M$ is connected, $H$ is a positive constant.
Evaluating \eqref{eq:positive-character-law} at $s=1$ gives
$H\circ\Phi_1=e^{2a}H$, a contradiction. Thus $T\ne0$.

Finally, since $\Phi_s$ is holomorphic,
\eqref{eq:positive-character-law} yields
\[
\begin{aligned}
\Phi_s^*T
&=\ddbar\bigl((U\circ\Phi_s)^2\bigr)\\
&=\ddbar(e^{2as}U^2)\\
&=e^{2as}T.
\end{aligned}
\]
Suppose now that $\log U$ is plurisubharmonic.
Then $S=\ddbar U$ is positive and closed.
If $S=0$, then $\ddbar U=0$ in the distributional sense. By Weyl's lemma,
$U$ is smooth and pluriharmonic.
Since $U\ge0$ and $U\not\equiv0$, the
strong minimum principle and connectedness
of $M$ imply that $U>0$ on $M$.
Consequently,
\[
0\leq\ddbar\log U
=-\frac{\ii}{U^2}\,
\partial U\wedge\bar\partial U.
\]
It follows that $\partial U=0$, so $U$ is constant on $M$,
contrary to $U\circ\Phi_1=e^a U$. Thus $S\ne0$. Furthermore,
\[
\Phi_s^*S
=\ddbar(U\circ\Phi_s)
=\ddbar(e^{as}U)
=e^{as}S.
\]

Suppose instead that \eqref{eq:positive-character-law} is assumed
only for $s=1$. The preceding nonvanishing arguments still apply.
Since $\Phi_m=\Phi_1^m$ for $m\in\Z$, we have
\[
U\circ\Phi_m=e^{am}U,\qquad
\Phi_m^*T=e^{2am}T
\qquad(m\in\Z).
\]
If $\log U$ is plurisubharmonic, we also have
$\Phi_m^*S=e^{am}S$ for every $m\in\Z$.
\end{proof}

We first observe that a positive closed current
satisfying a character identity at integer times
can be averaged to obtain the corresponding
identity for the full real action.

Suppose a nonzero positive closed current $T_0$ satisfies
$\Phi_1^*T_0=e^{k}T_0$ for some $k\in\R$. Set
\begin{equation}\label{eq:time-one-average}
 T=\int_0^1e^{-ks}\Phi_s^*T_0\,ds .
\end{equation}
For a test form $\beta$, the pairing
$(\Phi_s^*T_0)(\beta)=T_0(\Phi_{-s}^*\beta)$ depends continuously on $s$.
The supports of the transported forms for $0\le s\le1$ lie in a
fixed compact set. The average is positive and closed, and it is
nonzero because a positive test form pairing positively with $T_0$
also does so for $s$ close to zero. The group law
\eqref{eq:real-action-group-law} gives
\begin{equation}\label{eq:time-one-average-character}
\Phi_t^*T=e^{kt}\int_t^{t+1}e^{-kr}\Phi_r^*T_0\,dr=e^{kt}T.
\end{equation}
The second equality follows from the period-one identity
$e^{-k(r+1)}\Phi_{r+1}^*T_0=e^{-kr}\Phi_r^*T_0$.

Let $u$ be a positive pluriharmonic function on $M$ satisfying
\[
u\circ\Phi_s=e^{\lambda s}u
\qquad(s\in\R),
\]
for some $\lambda\in\R\setminus\{0\}$.
For any $k\in\R$, put $q=k/\lambda$ and define
\begin{equation}\label{eq:pluriharmonic-signed-character}
T_k=u^{q-2}\ii\,\partial u\wedge\bar\partial u.
\end{equation}
Since $u>0$, this is a smooth semipositive real
$(1,1)$-form for every $q\in\R$. Moreover,
$\partial\bar\partial u=0$ implies
$d(\partial u)=d(\bar\partial u)=0$. Consequently,
\[
\begin{aligned}
dT_k
&=\ii(q-2)u^{q-3}\,
du\wedge\partial u\wedge\bar\partial u =0.
\end{aligned}
\]
Thus $T_k$ is closed.

Differentiating the identity
$u\circ\Phi_s=e^{\lambda s}u$ at $s=0$ gives
$Yu=\lambda u$. Since $\lambda\ne0$ and $u>0$,
the differential $du$ is nowhere zero.
Hence $\partial u$ is nowhere zero, and $T_k$
has rank one at every point.

Since $\Phi_s$ is holomorphic, we have
\[
\begin{aligned}
\Phi_s^*T_k
&=\ii\,(u\circ\Phi_s)^{q-2}
  \partial(u\circ\Phi_s)
  \wedge\bar\partial(u\circ\Phi_s)\\
  &=\ii\,(e^{\lambda s}u)^{q-2}
  \partial(e^{\lambda s}u)
  \wedge\bar\partial(e^{\lambda s}u) =e^{q\lambda s}T_k
=e^{ks}T_k.
\end{aligned}
\]

Suppose instead that $u\circ\Phi_1=e^\lambda u$.
By the group law \eqref{eq:real-action-group-law},
\[
u\circ\Phi_m=e^{\lambda m}u
\qquad(m\in\Z),
\]
and the same calculation gives
\[
\Phi_m^*T_k=e^{km}T_k
\qquad(m\in\Z).
\]
The form $T_k$ is nonzero: otherwise
$\partial u=0$, so $u$ is constant on the
connected manifold $M$, contradicting
$u\circ\Phi_1=e^\lambda u$.

Applying \eqref{eq:time-one-average} to $T_k$,
we obtain a nonzero smooth positive closed
real $(1,1)$-form
\[
\widetilde T_k
=\int_0^1 e^{-ks}\Phi_s^*T_k\,ds
\]
satisfying
\[
\Phi_t^*\widetilde T_k=e^{kt}\widetilde T_k
\qquad(t\in\R),
\]
by \eqref{eq:time-one-average-character}.

We next construct plurisubharmonic functions
satisfying the character identity
\eqref{eq:positive-character-law}
by integrating nonnegative plurisubharmonic
functions along the flow.

\begin{proposition}\label{prop:laplace-general}
Let $M$ be a connected complex manifold, and let
$\{\Phi_s\}_{s\in\R}$ be a holomorphic $\R$-action
on $M$. Fix $K>0$, and suppose that
$b\in\PSH(M)\cap L^\infty_{\mathrm{loc}}(M)$
is nonnegative and not identically zero. Assume that
\begin{equation}\label{eq:laplace-psh-potential}
U(z)=\int_{\R}e^{-Ks/2}b(\Phi_s(z))\,ds
\end{equation}
converges locally uniformly on $M$. Then
$U$ is a nonnegative, locally bounded,
plurisubharmonic function, not identically
zero, satisfying
\begin{equation}\label{eq:laplace-character}
U\circ\Phi_t=e^{Kt/2}U
\qquad(t\in\R).
\end{equation}
Consequently,
$T=\ddbar(U^2)$ is a nonzero positive closed
real $(1,1)$-current on $M$ satisfying
\eqref{eq:current-homogeneity}.
\end{proposition}

\begin{proof}
For $R>0$, set
\[
U_R(z)=\int_{-R}^{R}
e^{-Ks/2}b(\Phi_s(z))\,ds.
\]
We first show that $U_R$ is plurisubharmonic.
For every compact set $C\Subset M$, the set
\[
C_R=\{\Phi_s(z):z\in C,\ |s|\le R\}
\]
is compact. Since $b$ is locally bounded,
there exists $B_{C,R}>0$ such that
\[
0\le b(\Phi_s(z))\le B_{C,R},
\qquad z\in C,\quad |s|\le R.
\]

Let $z_j\to z$ and choose a compact neighborhood $C$ containing
$z$ and all sufficiently late $z_j$. Apply Fatou's lemma to
$e^{-Ks/2}\{B_{C,R}-b(\Phi_s(z_j))\}$ on $[-R,R]$ to obtain
\[
\begin{aligned}
\limsup_{j\to\infty}U_R(z_j)
&\le
\int_{-R}^{R}e^{-Ks/2}
\limsup_{j\to\infty}
b(\Phi_s(z_j))\,ds \le
\int_{-R}^{R}e^{-Ks/2}
b(\Phi_s(z))\,ds =U_R(z),
\end{aligned}
\]
where the second inequality follows from
the upper semicontinuity of $b$ and the
continuity of $\Phi_s$. Thus $U_R$ is
upper semicontinuous. Let $f$ be a holomorphic map from a
neighborhood of $\overline{\Delta}$ into $M$.
Applying the mean-value inequality to
$b\circ\Phi_s\circ f$ and Tonelli's theorem
give
\[
\begin{aligned}
U_R(f(0))
&\le \int_{-R}^{R}\frac{e^{-Ks/2}}{2\pi}
\int_0^{2\pi}
b(\Phi_s(f(e^{\ii\theta})))\,d\theta\,ds =\frac{1}{2\pi}\int_0^{2\pi}
U_R(f(e^{\ii\theta}))\,d\theta.
\end{aligned}
\]
Together with upper semicontinuity,
this shows that $U_R$ is plurisubharmonic. By our assumption,
the functions $U_R$ converge locally
uniformly to $U$ as $R\to\infty$.
Consequently, $U$ is nonnegative,
plurisubharmonic, and locally bounded.

We next show that $U\not\equiv0$.
Choose $p\in M$ with $b(p)>0$.
In a holomorphic coordinate neighborhood
of $p$, the function $b$ is subharmonic.
The mean-value inequality on sufficiently
small Euclidean balls centered at $p$
therefore shows that $b$ is positive
on a set of positive measure.
Let $dV$ be any smooth positive volume
form on $M$. Since each $\Phi_s$ is a
diffeomorphism,
\[
\int_M b(\Phi_s(z))\,dV(z)>0
\qquad(s\in\R),
\]
where the integrals may be infinite.
By Tonelli's theorem,
\[
\begin{aligned}
\int_M U(z)\,dV(z)
&=\int_M\int_{\R}
e^{-Ks/2}b(\Phi_s(z))\,ds\,dV(z) =\int_{\R}e^{-Ks/2}
\left(\int_M b(\Phi_s(z))\,dV(z)\right)ds
>0.
\end{aligned}
\]
Hence $U$ is not identically zero. For $t\in\R$, the group law
\eqref{eq:real-action-group-law}
and the change of variables $r=s+t$ give
\[
\begin{aligned}
U(\Phi_t(z))
&=\int_{\R}e^{-Ks/2}
b(\Phi_{s+t}(z))\,ds =e^{Kt/2}\int_{\R}e^{-Kr/2}
b(\Phi_r(z))\,dr =e^{Kt/2}U(z).
\end{aligned}
\]
This proves \eqref{eq:laplace-character}.
Applying Lemma~\ref{lem:positive-character}
with $a=K/2$, we conclude that
$T=\ddbar(U^2)$ is a nonzero positive closed
real $(1,1)$-current satisfying
\eqref{eq:current-homogeneity}.
\end{proof}

\begin{remark}\label{rem:laplace-convergence}
A sufficient condition for the local uniform
convergence of \eqref{eq:laplace-psh-potential}
is the following: for every compact set
$C\Subset M$, there exist constants
$A_C,B_C>0$ and $\gamma_C<K/2$ such that
\[
b(\Phi_s(z))=0
\qquad(z\in C,\ s<-A_C)
\]
and
\[
b(\Phi_s(z))\le B_Ce^{\gamma_Cs}
\qquad(z\in C,\ s\ge0).
\]
For $s<-A_C$, the integrand in
\eqref{eq:laplace-psh-potential} vanishes on $C$.
For $s\ge0$, we have
\[
e^{-Ks/2}b(\Phi_s(z))
\le B_Ce^{-(K/2-\gamma_C)s},
\qquad z\in C.
\]
Since $\gamma_C<K/2$, it follows that
\[
\begin{aligned}
\sup_{z\in C}\int_R^\infty
e^{-Ks/2}b(\Phi_s(z))\,ds
&\le B_C\int_R^\infty
e^{-(K/2-\gamma_C)s}\,ds =\frac{B_C}{K/2-\gamma_C}
e^{-(K/2-\gamma_C)R}
\longrightarrow0
\end{aligned}
\]
as $R\to\infty$. Moreover, the set
\[
\{\Phi_s(z):z\in C,\ -A_C\le s\le0\}
\]
is compact. Since $b$ is locally bounded,
the integral over $[-A_C,0]$ is uniformly
bounded on $C$. Combining this with the vanishing condition
$b(\Phi_s(z))=0$ for $z\in C$ and $s<-A_C$,
we conclude that the integral in
\eqref{eq:laplace-psh-potential}
converges uniformly on $C$.
\end{remark}

When $b$ is bounded, the forward-time estimate in
Remark~\ref{rem:laplace-convergence} is automatic.
Thus only the vanishing of $b\circ\Phi_s$ for
sufficiently negative $s$ needs to be verified.
This gives the following construction.

\begin{lemma}\label{lem:backward-cutoff-current}
Let $\Omega$ be a connected complex manifold, let
$\{\Phi_s\}_{s\in\R}\subset\Aut(\Omega)$ be a holomorphic $\R$-action, and
let $K>0$.  Suppose that
\[
 b\in\PSH(\Omega)\cap L^\infty(\Omega),\qquad b\ge0,
 \qquad b\not\equiv0,
\]
and that for every compact set $C\Subset\Omega$ there is $A_C>0$ such that
\begin{equation}\label{eq:backward-cutoff-vanishing}
 b(\Phi_s z)=0\qquad(z\in C,\ s<-A_C).
\end{equation}
Then there is a nonzero positive closed $(1,1)$-current $T$ on $\Omega$
satisfying
\begin{equation}\label{eq:backward-current-character}
 \Phi_t^*T=e^{Kt}T\qquad(t\in\R).
\end{equation}
\end{lemma}

\begin{proof}
Since $b$ is bounded, for every compact
$C\Subset\Omega$ and every $R\ge0$,
\[
\begin{aligned}
\sup_{z\in C}\int_R^\infty
e^{-Ks/2}b(\Phi_s(z))\,ds
&\le \|b\|_\infty
\int_R^\infty e^{-Ks/2}\,ds =\frac{2\|b\|_\infty}{K}e^{-KR/2}.
\end{aligned}
\]
Together with
\eqref{eq:backward-cutoff-vanishing},
Remark~\ref{rem:laplace-convergence}
shows that the integral
\begin{equation}\label{eq:backward-laplace}
\Psi(z)=\int_{\R}e^{-Ks/2}b(\Phi_s(z))\,ds
\end{equation}
converges locally uniformly on $\Omega$.
By Proposition~\ref{prop:laplace-general},
$T=\ddbar(\Psi^2)$ is a nonzero positive
closed real $(1,1)$-current satisfying
\eqref{eq:backward-current-character}.
\end{proof}

\begin{proof}[Proof of Theorem~\ref{thm:backward-psh-prescribed-intro}]
Choose $p\in\Omega$ with $u(p)>0$ and fix $0<\delta<u(p)$. Set
\[
b=(u-\delta)_+=\max\{u-\delta,0\}.
\]
The function $b$ is plurisubharmonic, bounded and nonnegative,
and $b(p)=u(p)-\delta>0$. For each compact set $C\Subset\Omega$,
\eqref{eq:backward-psh-decay} gives $A_C>0$ such that
\[
u(\Phi_s(z))<\delta\qquad(z\in C,\ s<-A_C).
\]
Thus $b$ satisfies \eqref{eq:backward-cutoff-vanishing}.
Lemma~\ref{lem:backward-cutoff-current} gives a nonzero
$T\in\mathcal T_{Y,K}(\Omega)$, and
\cite[Theorem~1.2]{GangulyKRsol2026Paper21} gives the soliton.
\end{proof}

\subsection{Stein cyclic quotients and holomorphic maps}

The next construction replaces the backward-decay hypothesis
by the Steinness of the quotient under the time-one map of a
proper holomorphic $\R$-action. Cartan's Theorem B then
gives a holomorphic time function.

\begin{proposition}\label{prop:stein-cyclic-character}
Let $M$ be a connected complex manifold,
and let $\{\Psi_t\}_{t\in\R}$ be a proper
holomorphic $\R$-action on $M$.
Suppose that the quotient
\[
S=M/\langle\Psi_1\rangle
\]
is Stein. Then there exists a holomorphic
function $F\in\mathcal O(M)$ satisfying
\begin{equation}\label{eq:stein-cyclic-time}
F\circ\Psi_t=F+t
\qquad(t\in\R).
\end{equation}

Consequently, for every $K>0$, the form
\[
T_K=\ddbar e^{K\operatorname{Re}F}
\]
is a smooth positive closed real
$(1,1)$-form of rank one at every point
of $M$, satisfying
\[
\Psi_t^*T_K=e^{Kt}T_K
\qquad(t\in\R).
\]
\end{proposition}

\begin{proof}
The stabilizers of the proper real action are compact subgroups
of $\R$, hence trivial. Restricting the proper action map to
the closed subgroup $\Z$ gives a free, properly discontinuous
action by the maps $\Psi_m$, $m\in\Z$. Thus
\[
 \pi:M\longrightarrow S
\]
is a holomorphic covering. Choose an open cover $\{V_i\}$ of $S$
over which $\pi$ admits holomorphic
sections
\[
\sigma_i:V_i\longrightarrow M.
\]
Since the fibers of $\pi$ are the orbits
of $\langle\Psi_1\rangle$, there is a unique
integer $n_{ij}(x)$ such that
\[
\sigma_i(x)=\Psi_{n_{ij}(x)}(\sigma_j(x)),
\qquad x\in V_i\cap V_j.
\]
As $\pi$ is a covering map, the functions
$n_{ij}:V_i\cap V_j\to\Z$ are locally
constant. By \eqref{eq:real-action-group-law}
and the freeness of the action,
\[
n_{ij}+n_{jk}=n_{ik}
\qquad\text{on }V_i\cap V_j\cap V_k.
\]

Since each $n_{ij}$ is locally constant,
the functions $n_{ij}$ form a
\v{C}ech $1$-cocycle with values in
$\mathcal O_S$.
As $S$ is Stein, Cartan's Theorem B
\cite[Chapter~IV]{GrauertRemmert2004}
gives
\[
H^1(S,\mathcal O_S)=0.
\]
Consequently, after refining the cover
if necessary, there exist holomorphic
functions $a_i\in\mathcal O(V_i)$ such that
\[
n_{ij}=a_i-a_j
\qquad\text{on }V_i\cap V_j.
\]

For $x\in V_i$ and $m\in\Z$, define
\[
h(\Psi_m(\sigma_i(x)))=m+a_i(x).
\]
These definitions agree on overlaps.
Indeed, for $x\in V_i\cap V_j$ and
$m\in\Z$, we have
\[
\Psi_m(\sigma_i(x))
=\Psi_{m+n_{ij}(x)}(\sigma_j(x))
\]
and
\[
m+a_i(x)=m+n_{ij}(x)+a_j(x).
\]
Thus $h$ is a globally defined
holomorphic function on $M$ satisfying
\begin{equation}\label{eq:cyclic-discrete-time}
h\circ\Psi_1=h+1.
\end{equation}

Now define
\[
F(z)=\int_0^1h(\Psi_s(z))\,ds.
\]
The function $F$ is holomorphic, since
$h\circ\Psi_s$ is holomorphic in $z$
and depends continuously on $s$. By \eqref{eq:real-action-group-law},
\[
F(\Psi_t(z))
=\int_t^{t+1}h(\Psi_s(z))\,ds.
\]
Differentiating with respect to $t$
and using \eqref{eq:cyclic-discrete-time},
we obtain
\[
\begin{aligned}
\frac{d}{dt}F(\Psi_t(z))
&=h(\Psi_{t+1}(z))-h(\Psi_t(z)) =1.
\end{aligned}
\]
Since $F(\Psi_0(z))=F(z)$, integration
proves \eqref{eq:stein-cyclic-time}. Let $Y$ be the infinitesimal generator
of $\{\Psi_t\}$. Differentiating
\eqref{eq:stein-cyclic-time} at $t=0$
gives
\[
YF=1.
\]
In particular, $\partial F$ is nowhere
vanishing. Since $F$ is holomorphic,
\[
\begin{aligned}
T_K
&=\ddbar e^{K\operatorname{Re}F}=\frac{K^2}{4}e^{K\operatorname{Re}F}
  \ii\,\partial F\wedge
  \bar\partial\overline F.
\end{aligned}
\]
Thus $T_K$ is a smooth positive closed real
$(1,1)$-form of rank one at every point.

Finally, since $\Psi_t$ is holomorphic,
\eqref{eq:stein-cyclic-time} gives
\[
\begin{aligned}
\Psi_t^*T_K
&=\ddbar e^{K\operatorname{Re}(F\circ\Psi_t)}=\ddbar e^{K(\operatorname{Re}F+t)}=e^{Kt}T_K.
\end{aligned}
\]
\end{proof}

\begin{proof}[Proof of Theorem~\ref{thm:equivariant-stein-target}]
Let $F\in\mathcal O(M)$ be given by
Proposition~\ref{prop:stein-cyclic-character},
and set $H=F\circ f$.
Equations~\eqref{eq:flow-map-relation} and
\eqref{eq:stein-cyclic-time} give
\[
H\circ\Phi_t=H+t
\qquad(t\in\R).
\]

Differentiating at $t=0$ gives $YH=1$. For $K>0$,
\begin{equation}\label{eq:equivariant-target-forms}
\begin{aligned}
T_\pm
&=\ddbar e^{\pm K\Re H}=\frac{K^2}{4}e^{\pm K\Re H}
  \ii\,\partial H\wedge\bar\partial\overline H.
\end{aligned}
\end{equation}
Since $YH=1$, the holomorphic one-form $\partial H$
is nowhere vanishing. Thus both forms have rank
one at every point. Moreover,
\[
\begin{aligned}
\Phi_t^*T_\pm
&=\ddbar e^{\pm K\operatorname{Re}(H\circ\Phi_t)} =e^{\pm Kt}T_\pm.
\end{aligned}
\]
The soliton assertion follows from
\cite[Theorem~1.2]{GangulyKRsol2026Paper21}.
\end{proof}

The forms constructed in the preceding proof are
pullbacks of the corresponding forms on $M$.
This argument applies to any homogeneous positive
closed form on the target, independently of the
cyclic-quotient construction.

\begin{remark}\label{rem:app-pullback}
Let $f:\Omega\to M$ be holomorphic, and let $\Phi_t$
and $\Psi_t$ be holomorphic $\R$-actions on $\Omega$
and $M$, with generators $Y$ and $X$, respectively.
Suppose that \eqref{eq:flow-map-relation} holds.
If $S\in\mathcal T_{X,K}(M)$ is a smooth form
for some $K>0$, then $f^*S$ is positive and closed,
and
\[
\begin{aligned}
\Phi_t^*(f^*S)
&=(f\circ\Phi_t)^*S =(\Psi_t\circ f)^*S=f^*(\Psi_t^*S) =e^{Kt}f^*S.
\end{aligned}
\]
Thus $f^*S\in\mathcal T_{Y,K}(\Omega)$.
If $f$ is a surjective submersion and $S\ne0$,
then $f^*S\ne0$, since a tangent vector on which
$S$ is positive can be lifted through $df$.

More generally, whenever $f^*S\ne0$ and
$\Omega\Subset Z$ is a Stein domain of complex
dimension $n\ge2$ such that a neighborhood of
$\overline\Omega$ admits a holomorphic immersion
into $\C^N$, the form $f^*S$ yields a complete
nontrivial expander with field $Y$ and expansion
constant $K$, by
\cite[Theorem~1.2]{GangulyKRsol2026Paper21}.
\end{remark}

\begin{corollary}
\label{cor:app-circular-base}
Let $\Omega\Subset\C^m_z\times\C^r_w$ be a bounded
pseudoconvex domain invariant under
$(z,w)\mapsto(z,e^{i\theta}w)$ and satisfying
$(z,w)\in\Omega\implies(z,0)\in\Omega$. Let $\Phi_t$ be a holomorphic $\R$-action on $\Omega$ commuting
with the circle action, and let $Y$ be its generator. Put
\[
D=\{z:(z,0)\in\Omega\},
\]
and denote by $\Psi_t$ the induced flow on $D$. If, for some
$K>0$, there exists a nonzero smooth positive closed real
$(1,1)$-form $S$ on $D$ satisfying
\[
\Psi_t^*S=e^{Kt}S,
\]
then $\Omega$ admits a complete nontrivial expanding
K\"ahler--Ricci soliton with soliton field $Y$ and
expansion constant $K$.
\end{corollary}

\begin{proof}
The fixed set $\{(z,0):z\in D\}$ is preserved by $\Phi_t$.
Write $\Phi_t=(A_t,B_t)$. Since $\Phi_t$ commutes with
the circle action,
\[
A_t(z,e^{\ii\theta}w)=A_t(z,w).
\]
For each $z\in D$, expand $A_t$ in the $w$-variables
near $(z,0)$:
\[
A_t(z,w)=\sum_{\alpha\in\mathbb Z_{\ge0}^r}
a_\alpha(z)w^\alpha.
\]
The invariance above implies
\[
(e^{\ii|\alpha|\theta}-1)a_\alpha(z)=0
\qquad(\theta\in\R).
\]
Thus $a_\alpha(z)=0$ whenever $|\alpha|>0$.
Consequently, $A_t(z,w)=A_t(z,0)$ in a
neighborhood of the zero section.

By hypothesis, $(z,w)\in\Omega$ implies $(z,0)\in\Omega$,
so $A_t(z,w)-A_t(z,0)$ is holomorphic on $\Omega$.
Since $\Omega$ is connected, the identity theorem gives
\[
A_t(z,w)=A_t(z,0)
\qquad((z,w)\in\Omega).
\]
Thus the holomorphic map $f(z,w)=z$ satisfies
\eqref{eq:flow-map-relation} for the flows $\Phi_t$ and $\Psi_t$.
It is a surjective submersion. Since $\Omega$ is a bounded
pseudoconvex Euclidean domain, the metric conclusion of
Remark~\ref{rem:app-pullback} applies.
\end{proof}

\subsection{Holomorphic Jacobians and one-sided estimates}

The Laplace construction also applies to holomorphic functions
that decay along both ends of a flow. Such functions will arise
below from the pullbacks of square-integrable
holomorphic $n$-forms under proper actions.

\begin{lemma}\label{lem:two-ended-holomorphic}
Let $M$ be a connected complex manifold, and let
$\{\Phi_t\}_{t\in\R}$ be a holomorphic $\R$-action
on $M$. Suppose that a nonzero
$f\in\mathcal O(M)$ satisfies
\begin{equation}\label{eq:two-ended-holomorphic}
f\circ\Phi_t\longrightarrow0
\quad\text{locally uniformly on }M
\quad\text{as }t\to\pm\infty.
\end{equation}
Then, for every $K>0$, there exist nonzero
smooth positive closed real $(1,1)$-forms
$T_{K,+}$ and $T_{K,-}$ on $M$ satisfying
\[
\Phi_t^*T_{K,+}=e^{Kt}T_{K,+},
\qquad
\Phi_{-t}^*T_{K,-}=e^{Kt}T_{K,-}
\qquad(t\in\R).
\]
\end{lemma}

\begin{proof}
Choose $p\in M$ with $f(p)\ne0$ and
$\delta\in(0,|f(p)|^2)$.
Let $\chi:[0,\infty)\to[0,\infty)$ be
a smooth convex nondecreasing function
such that
\[
\chi=0\quad\text{on }[0,\delta],
\qquad
\chi(|f(p)|^2)>0.
\]
Then
\[
b=\chi(|f|^2)
\]
is a smooth, nonnegative, nonzero
plurisubharmonic function on $M$.

By \eqref{eq:two-ended-holomorphic},
for every compact set $C\Subset M$,
there exists $A_C>0$ such that
\begin{equation}\label{eq:two-ended-cutoff-vanishing}
b(\Phi_t(z))=0
\qquad(z\in C,\ |t|\ge A_C).
\end{equation}
For $K>0$, define
\[
U_{K,\pm}(z)=\int_{\R}
e^{-Kt/2}b(\Phi_{\pm t}(z))\,dt.
\]
Let $V\Subset M$ be open and put $C=\overline V$.
By \eqref{eq:two-ended-cutoff-vanishing},
\[
U_{K,\pm}(z)=\int_{-A_C}^{A_C}
e^{-Kt/2}b(\Phi_{\pm t}(z))\,dt,
\qquad z\in V.
\]
Thus both integrals converge locally uniformly on $M$.
Since the integrands are smooth in $(t,z)$,
differentiation under the integral on the fixed
compact interval gives $U_{K,\pm}\in C^\infty(V)$.
As $V$ is arbitrary, both functions are smooth on $M$.
Proposition~\ref{prop:laplace-general} now gives
the nonzero smooth positive closed forms
\[
T_{K,\pm}=\ddbar(U_{K,\pm}^2)
\]
with the asserted transformation laws.
\end{proof}

\begin{lemma}\label{lem:jacobian-log}
Let $D$ be a connected complex manifold of dimension $n$,
and let $\Phi$ be a holomorphic $\R^p$-action on $D$.
Suppose that $\sigma$ is a nowhere-vanishing
holomorphic $n$-form on $D$. For each $t\in\R^p$,
define the nowhere-vanishing holomorphic function $L_t$
by $\Phi_t^*\sigma=L_t\sigma$.
There are holomorphic functions \(j_t\), continuous in \((t,z)\),
such that
\begin{equation}\label{eq:exact-jacobian-cocycle}
 e^{j_t}=L_t,\qquad j_0=0,\qquad
 j_{t+s}=j_t\circ\Phi_s+j_s.
\end{equation}
If \(Z_\alpha=X_{e_\alpha}^{1,0}\), write
\(\Lie_{Z_\alpha}\sigma=h_\alpha\sigma\).
Then
\[
j_t(z)=\int_0^1\sum_{\alpha=1}^p
t_\alpha h_\alpha(\Phi_{rt}z)\,dr .
\]
\end{lemma}

\begin{proof}
For fixed \(t\), differentiation of the pullback gives
\[
\frac{d}{dr}(\Phi_{rt}^*\sigma)
 =\Phi_{rt}^*(\Lie_{X_t}\sigma)
 =\left(\sum_\alpha t_\alpha h_\alpha\circ\Phi_{rt}\right)
       \Phi_{rt}^*\sigma ,
\]
since $\iota_{\overline{Z_\alpha}}\sigma=0$ and $d\sigma=0$.
The initial condition $L_0=1$ gives the radial integral for $j_t$.
For fixed $z$, the cocycle identity for $L_t$ implies that
$j_{t+s}(z)-j_t(\Phi_sz)-j_s(z)$ has values in $2\pi\ii\Z$.
It is continuous in $(t,s)$ and vanishes at $(0,0)$.
Since $\R^{2p}$ is connected, it vanishes identically.
This proves \eqref{eq:exact-jacobian-cocycle}.
\end{proof}

\begin{remarks}\label{rem:jacobian-frame-remarks}
\begin{enumerate}[
    label=\textup{(\roman*)},
    ref=\thetheorem\textup{(\roman*)},
    leftmargin=0pt,
    labelindent=0pt,
    labelsep=0.45em,
    itemindent=*
]
\item\label{rem:jacobian-coordinate}
Suppose that, in a holomorphic coordinate chart,
\[
\sigma=\rho\,dz^1\wedge\cdots\wedge dz^n,
\qquad
Z_\alpha=\sum_{k=1}^n
a_\alpha^k\frac{\partial}{\partial z^k}.
\]
When $z$ and $\Phi_t(z)$ belong to this chart,
the multiplier of $\sigma$ and the coefficients
$h_\alpha$ in Lemma~\ref{lem:jacobian-log}
are given by
\begin{equation}\label{eq:frame-jacobian}
\begin{aligned}
L_t
&=J_t\frac{\rho\circ\Phi_t}{\rho},
\qquad J_t=\det_\C D\Phi_t,\\
h_\alpha
&=\sum_{k=1}^n
\frac{\partial a_\alpha^k}{\partial z^k}
+\frac{Z_\alpha\rho}{\rho}.
\end{aligned}
\end{equation}
The first identity follows by pulling back
$\sigma$. To verify the Jacobian formula, take
a flow segment contained in the chart, and put
\[
v_t^k=\sum_{\alpha=1}^p t_\alpha a_\alpha^k,
\qquad A_r=D\Phi_{rt}.
\]
Differentiation gives
\[
A_r'=Dv_t(\Phi_{rt})A_r.
\]
By Jacobi's formula,
\[
\begin{aligned}
\frac{d}{dr}J_{rt}
&=J_{rt}\tr(A_r^{-1}A_r')=J_{rt}\sum_{\alpha,k}t_\alpha
\left(\frac{\partial a_\alpha^k}
{\partial z^k}\right)\circ\Phi_{rt}.
\end{aligned}
\]
The expression for $h_\alpha$ follows directly
by computing $\Lie_{Z_\alpha}\sigma$.
For the Euclidean coordinate frame $\rho=1$,
we have $e^{j_t}=J_t$.

\item\label{rem:jacobian-frame-change}
If $\sigma$ is replaced by $b\sigma$ with
$b\in\mathcal O^*(D)$, its multiplier becomes
\[
L_t'=L_t\frac{b\circ\Phi_t}{b}.
\]
Lemma~\ref{lem:jacobian-log} supplies a
holomorphic logarithm $j_t'$ normalized by
$j_0'=0$. When $b=e^g$ for some
$g\in\mathcal O(D)$, this logarithm is
\[
j_t'=j_t+g\circ\Phi_t-g.
\]
\end{enumerate}
\end{remarks}

Properness gives the following decay result for
square-integrable holomorphic $n$-forms and, when
$\Phi_t^*\sigma=L_t\sigma$, for the functions $L_t$.

\begin{lemma}\label{stein:lem:L2-decay}
Let $\Phi_t$ be a proper holomorphic $\R^p$-action on a complex
manifold $D$. If $\sigma\in H^0(D,K_D)$ is square-integrable,
then $\Phi_t^*\sigma\to0$ locally uniformly as $|t|\to\infty$.
If $\sigma$ is nowhere vanishing and
$\Phi_t^*\sigma=L_t\sigma$, then $L_t\to0$
locally uniformly as $|t|\to\infty$.
\end{lemma}

\begin{proof}
For any \(\epsilon>0\), choose a compact set \(L\subset D\) such
that \(\int_{D\setminus L}dV_\sigma<\epsilon\). Given open
\(V\Subset D\), properness makes the set of \(t\) with
\(\Phi_t(\overline V)\cap L\ne\varnothing\) compact. For all
sufficiently large \(|t|\), change of variables gives
\[
\int_V dV_{\Phi_t^*\sigma}=\int_{\Phi_t(V)}dV_\sigma
\le\int_{D\setminus L}dV_\sigma<\epsilon .
\]
In a coordinate chart \(A_t=\Phi_t^*\sigma=a_t\,dz^1\wedge\cdots
\wedge dz^n\). On concentric balls \(B_r\Subset B_{2r}\Subset V\)
in the chart, the submean inequality gives
\[
\sup_{B_r}|a_t|^2\le C_r\int_{B_{2r}}|a_t|^2\,dV_{\rm E}
\le n!C_r\int_V dV_{A_t}.
\]
A finite covering gives the asserted uniform local decay.
If $\sigma=\rho\,dz^1\wedge\cdots\wedge dz^n$
is nowhere vanishing in this coordinate chart, then
$\inf_{\overline{B_r}}|\rho|>0$.
Since $\Phi_t^*\sigma=a_t\,dz^1\wedge\cdots\wedge dz^n$,
we have $L_t=a_t/\rho\to0$ uniformly on $B_r$.
\end{proof}

Suppose $D\Subset Z$ and a holomorphic $n$-form
$\sigma$ is defined on a neighborhood of $\overline D$.
Its restriction to $D$ is square-integrable.
With any fixed smooth Hermitian metric on $K_D$,
Lemma~\ref{stein:lem:L2-decay} gives, for every $Q\Subset D$,
\begin{equation}\label{eq:intrinsic-canonical-decay}
\sup_Q\|\Phi_t^*\sigma\|\longrightarrow0
\qquad (|t|\longrightarrow\infty).
\end{equation}
For a bounded domain $D\subset\C^n$, take
$\sigma=dz^1\wedge\cdots\wedge dz^n$.
Then $\Phi_t^*\sigma=J_t\sigma$, and hence
\begin{equation}\label{eq:proper-jacobian-decay}
\sup_Q|J_t|\longrightarrow0
\qquad (|t|\longrightarrow\infty).
\end{equation}

We next combine this decay with an upper bound on
$\Re(j_t(z)-j_t(a))$ for $t\ge0$.

\begin{proposition}\label{prop:frame-halfline}
Let \(\Phi_t\) be a proper real holomorphic flow, with generator \(Y\),
on a connected complex manifold $D$ of dimension $n$
carrying a nowhere-vanishing holomorphic $n$-form
$\sigma$ with $\|\sigma\|_{2,D}<\infty$.
Suppose that $D$ admits
a smooth invariant volume form \(B_{\sigma}^{\,0}dV_\sigma\)
with \(B_{\sigma}^{\,0}>0\) and
\(\ddbar\log B_{\sigma}^{\,0}\ge0\). Let \(a\in D\) and let $j_t$ be the logarithm in Lemma~\ref{lem:jacobian-log}
for the frame $\sigma$. Consider the bound
\begin{equation}\label{eq:halfline-jacobian}
 \sup_{t\ge0,\ z\in L}\Re(j_t(z)-j_t(a))<\infty
 \qquad (L\Subset D).
\end{equation}
If this holds, then $\mathcal T_{Y,K}(D)$ contains a nonzero
smooth form for every $K>0$. If instead the bound holds with
$j_t$ replaced by $j_t+g\circ\Phi_t-g$ for some
$g\in\mathcal O(D)$, then $\mathcal T_{\varepsilon Y,K}(D)$
contains a nonzero smooth form for every $K>0$, with one
choice of $\varepsilon\in\{1,-1\}$ independent of $K$.
\end{proposition}

\begin{proof}
Write $\Phi_t^*\sigma=L_t\sigma$ with $L_t=e^{j_t}$, and put
\[
h_t=j_t+g\circ\Phi_t-g,\qquad q_t=h_t-h_t(a),
\]
taking $g=0$ in the first case. The bound is
$\sup_{t\ge0}\sup_L\Re q_t<\infty$ for every compact $L$.
On concentric coordinate balls $B(x,r)\Subset B(x,2r)$,
Borel--Carath\'eodory on complex lines through $x$
\cite[Chapter~XII, \S~3]{LangComplexAnalysis1999} gives
\[
 \sup_{B(x,r)}|q_t-q_t(x)|
 \le 2\bigl(\sup_{B(x,2r)}\Re q_t-\Re q_t(x)\bigr).
\]
Suppose that
\[
\sup_{t\ge0}|q_t(x)|\le A,
\qquad
M=\sup_{t\ge0}\sup_{B(x,2r)}\Re q_t<\infty.
\]
The preceding inequality gives
\[
\begin{aligned}
\sup_{t\ge0}\sup_{B(x,r)}|q_t|
&\le A+2(M+A) =2M+3A.
\end{aligned}
\]
Fix $z\in D$ and choose a path joining $a$ to $z$.
Cover its compact image by finitely many coordinate balls
$B(x_j,2r_j)\Subset D$, arranged so that $x_0=a$,
$x_m=z$, and $x_{j+1}\in B(x_j,r_j)$.
Starting with $q_t(a)=0$, the preceding estimate gives
uniform bounds for $q_t$ at each successive center
and then in a neighborhood of $z$.
Covering any compact set $L\Subset D$ by finitely
many such neighborhoods gives
\begin{equation}\label{eq:halfline-normal-bound}
\sup_{\substack{t\ge0\\z\in L}}|q_t(z)|<\infty.
\end{equation}

Set $G_N=N^{-1}\int_0^N q_t\,dt$. By Montel, along some
$N_j\to\infty$ we have $G_{N_j}\to G$ locally uniformly and
$G(a)=0$. Since $h_t$ satisfies the additive cocycle identity,
\[
q_{t+s}(z)=q_t(\Phi_s(z))-q_t(\Phi_s(a))+q_s(z).
\]
Integration gives
\[
G_N(\Phi_sz)-G_N(\Phi_sa)+q_s(z)-G_N(z)
=\frac1N\left(\int_N^{N+s}q_t(z)\,dt-\int_0^s q_t(z)\,dt\right).
\]
For each fixed $s\in\R$, the first integral is bounded on
compact sets by \eqref{eq:halfline-normal-bound} when $N>|s|$,
and the second by continuity on a compact interval, with both
integrals oriented when $s<0$. Passing to the limit shows
that $G\circ\Phi_s+h_s-G$ is the constant
$\ell(s)=G(\Phi_sa)+h_s(a)$. The cocycle makes $\ell$
continuous and additive, so $\ell(s)=\lambda s$. Thus
\begin{equation}\label{eq:halfline-canonical-multiplier}
G\circ\Phi_t+h_t-G=\lambda t .
\end{equation}
With $\widehat G=G+g$ this reads
\begin{equation}\label{eq:gauged-canonical-multiplier}
\widehat G\circ\Phi_t+j_t-\widehat G=\lambda t .
\end{equation}

If $g=0$, the cocycle gives
\[
G_N(\Phi_1a)+j_1(a)
=\frac1N\left(\int_N^{N+1}j_t(a)\,dt-\int_0^1j_t(a)\,dt\right).
\]
The quantity $\Re j_t(a)=\log|L_t(a)|$ is bounded above
for $t\ge0$ by Lemma~\ref{stein:lem:L2-decay} and
continuity. In the Euclidean coordinate frame,
Cauchy's estimate for the bounded maps $\Phi_t$
gives the same upper bound at $a$.
Passing to the limit along $N_j$ and using
\eqref{eq:halfline-canonical-multiplier}, we obtain
\[
\begin{aligned}
\Re\lambda
&=\lim_{j\to\infty}\frac1{N_j}
\left(
\int_{N_j}^{N_j+1}\Re j_t(a)\,dt
-\int_0^1\Re j_t(a)\,dt
\right)\le0.
\end{aligned}
\]

If $\Re\lambda=0$, set $f=e^{-\widehat G}$. Then
\begin{equation}\label{eq:halfline-decaying-function}
f\circ\Phi_t=e^{-\lambda t}L_t f .
\end{equation}
Here $|e^{-\lambda t}|=1$ and $L_t\to0$ locally uniformly as
$|t|\to\infty$ by Lemma~\ref{stein:lem:L2-decay} applied to
$\sigma$. Lemma~\ref{lem:two-ended-holomorphic} gives the
two nonzero smooth forms for $Y$ and $-Y$.

If \(b=\Re\lambda\ne0\), invariance of \(B_{\sigma}^{\,0}dV_\sigma\)
gives \(\log B_{\sigma}^{\,0}\circ\Phi_t+2\Re j_t=\log B_{\sigma}^{\,0}\).
Combining this with \eqref{eq:gauged-canonical-multiplier}, we
obtain, for $H=\log B_{\sigma}^{\,0}-2\Re\widehat G$,
\begin{equation}\label{eq:halfline-real-weight}
H\circ\Phi_t=H-2bt,\qquad \ddbar H=\ddbar\log B_{\sigma}^{\,0}.
\end{equation}
Put \(\varepsilon=-\operatorname{sgn}(b)\) and \(c=K/(2|b|)\).
When $g=0$, the last inequality gives $b<0$, so
$\varepsilon=1$. Thus
\[
T_K=\ddbar e^{cH}
=e^{cH}\left(c\,\ddbar H+c^2\ii\,\partial H\wedge\bar\partial H\right).
\]
By \eqref{eq:halfline-real-weight},
\[
\Phi_{\varepsilon t}^*T_K=\ddbar e^{c(H+2|b|t)}=e^{Kt}T_K.
\]
Also $YH=-2b$, so $\partial H$ is nowhere
zero and the second term is nonzero
everywhere on $D$.
If $\ddbar H>0$, the form is positive definite.
In particular, if the invariant volume form is
$\omega_E^n$ for a K\"ahler--Einstein metric with
$\Ric_{\omega_E}=-\omega_E$, then
$\ddbar\log B_\sigma^{\,0}=\omega_E>0$.
\end{proof}

\begin{remarks}\label{rem:halfline-observations}
\begin{enumerate}[
    label=\textup{(\roman*)},
    ref=\thetheorem\textup{(\roman*)},
    leftmargin=0pt,
    labelindent=0pt,
    labelsep=0.45em,
    itemindent=*
]
\item\label{rem:halfline-euclidean}
\label{prop:halfline-jacobian}
\label{cor:gauged-halfline}
For a bounded pseudoconvex domain
$D\Subset\C^n$, the coordinate frame
\[
\sigma=dz^1\wedge\cdots\wedge dz^n
\]
is square-integrable. The Bergman kernel form
$\mathcal B_D=B_\sigma^{\,0}dV_\sigma$
is invariant, with
\[
\ddbar\log B_\sigma^{\,0}>0.
\]
Thus Proposition~\ref{prop:frame-halfline}
applies to the holomorphic logarithms $j_t$
of the complex Jacobians of a proper
holomorphic real flow on $D$.

If \eqref{eq:halfline-jacobian} holds, then
for every $K>0$ there exists a complete
nontrivial expander with soliton field $Y$.
If the bound holds after replacing $j_t$ by
$j_t+g\circ\Phi_t-g$ for some
$g\in\mathcal O(D)$, the same conclusion
holds with field $\varepsilon Y$ for a fixed
$\varepsilon\in\{1,-1\}$ independent of $K$.
These conclusions follow from
\cite[Theorem~1.2]{GangulyKRsol2026Paper21}.

In the latter construction, $\Re\lambda=0$
gives nonzero smooth homogeneous forms for
both $Y$ and $-Y$, while $\Re\lambda\ne0$
gives a positive definite homogeneous form
for the selected field $\varepsilon Y$.

\item\label{rem:halfline-bergman}
\label{cor:halfline-bergman}
The bound \eqref{eq:halfline-jacobian} can also
be expressed in terms of the invariant volume
density. Put $R=\log B_\sigma^{\,0}$.
Invariance of $B_\sigma^{\,0}dV_\sigma$ gives
\[
R\circ\Phi_t+2\Re j_t=R.
\]
Subtracting this identity at $a$ from the
identity at $z$ gives
\[
\begin{aligned}
2\Re(j_t(z)-j_t(a))
&=R(z)-R(a) \quad+R(\Phi_t a)-R(\Phi_t z).
\end{aligned}
\]
Since $R(z)-R(a)$ is bounded on each
compact subset of $D$, condition
\eqref{eq:halfline-jacobian} is equivalent to
\[
\sup_{t\ge0}\sup_{z\in L}
\bigl(R(\Phi_t a)-R(\Phi_t z)\bigr)<\infty
\qquad(L\Subset D).
\]
For $D\Subset\C^n$ with
$\sigma=dz^1\wedge\cdots\wedge dz^n$,
we have $R=\log K_D(z,z)+\log(n!)$.
This additive constant cancels in the preceding difference.
\end{enumerate}
\end{remarks}

\subsection{Orbital \texorpdfstring{$L^2$}{L2} estimates for holomorphic
\texorpdfstring{$n$}{n}-forms}
\label{sec:finite-volume-frames}\label{stein:sec:intrinsic-L2}

In Section~\ref{sec:corank-complexification}, we take Fourier
transforms of $(\Phi_t^*\sigma)^{\otimes2}$ in $t\in\R^p$.
To ensure locally uniform convergence of these integrals,
we need estimates uniform on compact subsets of $M$ that
are stronger than the decay in Lemma~\ref{stein:lem:L2-decay}.

\begin{lemma}\label{stein:lem:intrinsic-mass}
\label{lem:shared-jacobian-L2}
Let $M$ be a connected complex manifold of
dimension $n$, and let
$\Phi:\R^p\times M\to M$ be a proper
holomorphic action. Suppose that
$\sigma\in H^0(M,K_M)$ satisfies
$\|\sigma\|_{2,M}<\infty$, and put
\[
A_t=\Phi_t^*\sigma.
\]
For every relatively compact open set
$V\Subset M$, define
\[
E_V=\{u\in\R^p:\Phi_u(\overline V)\cap\overline V\ne\varnothing\}.
\]
Then \(E_V\) is compact and, writing
\(|E_V|\) for its Lebesgue measure,
\begin{equation}\label{stein:eq:mass}
\int_{\R^p}\int_V \ii^{n^2}A_t\wedge\overline{A_t}\,dt
\le |E_V|\,\|\sigma\|_{2,M}^2<\infty.
\end{equation}
Moreover, fix a smooth Hermitian metric $h$
on $M$ and a smooth Hermitian metric and
connection on $K_M$. Let $\nabla^m A_t$
denote the iterated covariant derivative,
using the induced tensor-product connection
and the Levi-Civita connection of $h$.
For every $Q\Subset V$ and $m\ge0$,
\begin{equation}\label{stein:eq:derivative-mass}
\int_{\R^p}\sup_{z\in Q}|\nabla^m A_t(z)|^2\,dt
\le C_{Q,V,m}\,\|\sigma\|_{2,M}^2<\infty.
\end{equation}
The constant is independent of $t$ and $\sigma$.
In any holomorphic coordinate chart, write
$A_t=a_t(z)\,dz^1\wedge\cdots\wedge dz^n$.
Then $z\mapsto a_\cdot(z)$ is holomorphic with
values in $L^2(\R^p)$.
\end{lemma}

\begin{proof}
Properness of the action implies that $E_V$ is
compact. For $w\in M$, put
\[
B_w=\{t\in\R^p:w\in\Phi_t(V)\}.
\]
If $B_w\ne\varnothing$, choose $t_0\in B_w$.
For $t\in B_w$, write
\[
w=\Phi_t(z)=\Phi_{t_0}(z_0),
\qquad z,z_0\in V.
\]
Then $z_0=\Phi_{t-t_0}(z)$, so
$B_w\subset t_0+E_V$. Consequently,
\begin{equation}\label{stein:eq:multiplicity}
|B_w|\le |E_V|.
\end{equation}

Since $dV_{A_t}=\Phi_t^*dV_\sigma$,
change of variables and Tonelli's theorem give
\[
\begin{aligned}
\int_{\R^p}\int_V dV_{A_t}\,dt
&=\int_{\R^p}\int_{\Phi_t(V)}
dV_\sigma\,dt=\int_M |B_w|\,dV_\sigma(w)\le |E_V|\|\sigma\|_{2,M}^2.
\end{aligned}
\]
This proves \eqref{stein:eq:mass}. Fix $Q\Subset V$. Choose finitely many coordinate
balls covering $Q$, with slightly larger closures
contained in $V$. On each ball write
\[
A_t=a_t(z)\,dz^1\wedge\cdots\wedge dz^n.
\]
On each coordinate ball, the mean-value
inequality and Cauchy's estimates bound
the derivatives of $a_t$ of orders at most
$m$ by its $L^2$ norm on a larger ball.
The coefficients of the fixed connections
and their derivatives are bounded on
these balls. Consequently,
\[
\sup_Q|\nabla^m A_t|^2
\le C_{Q,V,m}\int_V dV_{A_t}.
\]
Here the constant also depends on the fixed
Hermitian metric and connection, but not on
$t$ or $\sigma$. Integrating in $t$ and using
\eqref{stein:eq:mass} proves
\eqref{stein:eq:derivative-mass}.

Finally, in a smaller coordinate ball, Taylor's
formula gives
\[
\begin{aligned}
&\left\|a_\cdot(z+h)-a_\cdot(z)
-\sum_jh^j\partial_j a_\cdot(z)
\right\|_{L^2(\R^p)} \le
C|h|^2
\left(
\int_{\R^p}\sup_Q|D^2a_t|^2\,dt
\right)^{1/2}
=O(|h|^2).
\end{aligned}
\]
The derivatives $\partial_j a_\cdot(z)$
belong to $L^2(\R^p)$ by
\eqref{stein:eq:derivative-mass}.
Thus the coefficient map is holomorphic
with values in $L^2(\R^p)$.
\end{proof}

If $\sigma$ is nowhere vanishing on $M$, write
$A_t=L_t\sigma$. Fix a holomorphic
coordinate chart $(U,z)$ and choose
$Q\Subset V\Subset U$, with $V$ open.
On $U$, write
\[
\sigma=\rho\,dz^1\wedge\cdots\wedge dz^n,
\qquad
A_t=a_t\,dz^1\wedge\cdots\wedge dz^n.
\]
Since $\rho$ is nowhere zero,
$L_t=a_t/\rho$. The derivatives of
$\rho^{-1}$ are bounded on $\overline V$.
The product rule and
\eqref{stein:eq:derivative-mass},
applied for orders $0,\ldots,m$, give
\begin{equation}\label{eq:shared-jacobian-L2}
\int_{\R^p}\sup_{z\in Q}
\sum_{|\alpha|=m}
|\partial_z^\alpha L_t(z)|^2\,dt
<\infty
\qquad(m\ge0).
\end{equation}
The same estimates show that
$z\mapsto L_\cdot(z)$ is holomorphic
with values in $L^2(\R^p)$.

On $D\subset\C^n$, write
$\sigma=\rho\,dz^1\wedge\cdots\wedge dz^n$.
Then $\sigma$ is square-integrable precisely when
$\int_D|\rho|^2\,dV_{\mathrm E}<\infty$, and its
multiplier is $L_t=J_t(\rho\circ\Phi_t)/\rho$.

The same estimates hold when $\R^p$ is replaced
by a lattice $\Lambda\subset\R^p$, with counting
measure in place of Lebesgue measure. Indeed,
if $B_w\cap\Lambda\ne\varnothing$ and
$t_0\in B_w\cap\Lambda$, then
\[
B_w\cap\Lambda
\subset t_0+(E_V\cap\Lambda).
\]
Hence
\[
\#(B_w\cap\Lambda)
\le \#(E_V\cap\Lambda)<\infty,
\]
uniformly in $w$. Summing over $\Lambda$
and applying the same Cauchy estimates
gives the discrete analogues of
\eqref{stein:eq:mass} and
\eqref{stein:eq:derivative-mass}.
\section{The Structure of the Automorphism Group}
\label{sec:central-reduction}

\subsection{Properness of the automorphism-group action}
\label{sec:bounded-realization}

We first use the Bergman metric to establish
properness of the automorphism-group action.
For Stein ambient manifolds, the immersion
hypothesis below follows from
\cite[Theorem~3]{NarasimhanEmbedding}.

\begin{lemma}\label{stein:lem:bounded}
Let $\Omega\Subset Z$ be a domain in a
complex manifold. Suppose that a neighborhood
$V$ of $\overline\Omega$ admits a holomorphic
immersion $F:V\to\C^N$.

Then bounded holomorphic functions on $\Omega$
separate tangent directions, and every
holomorphic map $\C\to\Omega$ is constant.
The Bergman kernel form $\mathcal B_\Omega$
is nowhere zero, and the associated Bergman
metric is positive definite.

Moreover, $\Aut(\Omega)$ is a finite-dimensional
real Lie group acting properly on $\Omega$.
If $X$ and $JX$ are both complete real
holomorphic vector fields on $\Omega$,
then $X=0$.
\end{lemma}
\begin{proof}
The components of $F|_\Omega$ are bounded since $F(\overline\Omega)$
is compact, and they separate tangent directions since $dF$ is
injective. If $h:\C\to\Omega$ is holomorphic, Liouville's theorem
applied to $F\circ h$ gives $dF_{h(\zeta)}\,dh_\zeta=0$.
Since $dF$ is injective, $dh=0$.

Fix \(x\in\Omega\) and \(0\ne v\in T^{1,0}_x\Omega\). Choose
indices with \(\sigma=dF_{j_1}\wedge\cdots\wedge dF_{j_n}\)
nonzero at \(x\), and a component \(f=F_j\) with \(df_x(v)\ne0\).
Set \(\tau=(f-f(x))\sigma\). These forms are square-integrable:
for any smooth positive volume \(dV\) near \(\overline\Omega\)
the densities \(dV_\sigma/dV\) and \(dV_\tau/dV\) are bounded
on the compact closure. Choose holomorphic coordinates near $x$, and for
$\eta\in A^2(\Omega,K_\Omega)$ write
$\eta=f_\eta(z)\,dz^1\wedge\cdots\wedge dz^n$.
Define the evaluation functional by
$E_z(\eta)=f_\eta(z)$.
Cauchy's estimates on a slightly larger coordinate
ball imply that $z\mapsto E_z$ is holomorphic
with values in the dual Hilbert space. In particular,
\[
E_x(\sigma)\ne0,\qquad E_x(\tau)=0,\qquad
(dE_x(v))(\tau)=df_x(v)E_x(\sigma)\ne0.
\]
Thus \(E_x\ne0\) and \(dE_x(v)\notin\C E_x\). By the definition
of the Bergman form and the Levi formula for a holomorphic Hilbert
norm,
\begin{equation}\label{eq:bergman-evaluation-levi}
-\ii\omega_B(v,\bar v)
=\frac{\|E_x\|^2\|dE_x(v)\|^2-
       |\langle dE_x(v),E_x\rangle|^2}{\|E_x\|^4}>0 .
\end{equation}
Change of variables in $A^2(\Omega,K_\Omega)$ shows that $g_B$ is
invariant under $\Aut(\Omega)$. The holomorphic elements form a
closed subgroup of $\Isom(\Omega,g_B)$ in the compact-open topology
since a locally uniform limit lying in the isometry group is holomorphic.
Thus $\Aut(\Omega)$ is a Lie group
\cite[Theorem~2.12]{AlexandrinoBettiol2015}.
We give the properness argument for $g_B$, which may be
incomplete. Compare \cite[Proposition~3.62]{AlexandrinoBettiol2015}.

Let $d_B$ be the distance of $g_B$ and let $B_B(q,r)$ be the
corresponding ball. Suppose $x_j\to x$ and $A_j x_j\to y$
with $A_j\in\Aut(\Omega)$. Then $A_jx\to y$ and $A_j^{-1}y\to x$.
Put
\[
 S=\{q\in\Omega:\{A_jq:j\ge1\}\text{ has compact closure in }\Omega\}.
\]
If $q\in S$ and $\overline{B_B(q,2r)}$ is compact, cover
$\{A_jq\}$ by finitely many balls $B_B(A_{j_\ell}q,r)$.
For $d_B(q,z)<r$ and $A_jq\in B_B(A_{j_\ell}q,r)$,
\[
 d_B(A_{j_\ell}^{-1}A_jz,q)
 \le d_B(z,q)+d_B(A_jq,A_{j_\ell}q)<2r .
\]
Thus $\{A_jz\}$ lies in finitely many translates of the compact
ball. This proves that $S$ is open. If $q_*\in\overline S$,
choose $r>0$ with $\overline{B_B(q_*,4r)}$ compact and choose
$q\in S$ with $d_B(q,q_*)<r$. The same estimate, using
$\overline{B_B(q,2r)}\subset\overline{B_B(q_*,4r)}$, shows that
$q_*\in S$. Since $x\in S$ and $\Omega$ is connected, $S=\Omega$.

Apply this argument to $A_j^{-1}$ starting at $y$. Since the maps
and their inverses are isometries, Arzel\`a--Ascoli on compact
subsets and a diagonal subsequence give locally uniform limits
$A$ and $B$ for $A_j$ and $A_j^{-1}$. The limits are holomorphic,
and $A\circ B=B\circ A=\id_\Omega$ follows by taking limits of
the compositions on compact neighborhoods. Thus $A_j$ has a
subsequence converging in $\Aut(\Omega)$, which proves properness.

If \(X,JX\) are complete real holomorphic fields, then
\([X,JX]=(\Lie_X J)X+J[X,X]=0\). Their commuting flows give
the map
\[
H_x(s+it)=\exp(sX)\exp(tJX)x,\qquad
\partial_t H_x=J\,\partial_s H_x .
\]
This is an entire curve, hence constant by the first part.
Differentiating at zero gives \(X(x)=0\) for every \(x\).
\end{proof}

\subsection{Lie-group structure and proper flows}

\begin{lemma}\label{lem:compact-adjoint-group}
Let $G$ be a connected real Lie group with Lie algebra
$\mathfrak g$. Suppose that for every $X\in\mathfrak g$,
the operator $\ad_X:\mathfrak g_\C\to\mathfrak g_\C$
is diagonalizable over $\C$ with spectrum in $\ii\R$.
Then $G\simeq\R^p\times C$ for some $p\ge0$ and a
compact connected Lie group $C$.
\end{lemma}

\begin{proof}
Write \(\mathfrak g\) for the Lie algebra of \(G\).
Let $\mathfrak n$ denote the nilradical, i.e., the largest nilpotent ideal, and let $\mathfrak r$ denote the
radical of $\mathfrak g$, namely its largest solvable ideal. Fix $X\in\mathfrak n$. Since $\mathfrak n$ is an ideal
of $\mathfrak g$, we have
$[X,\mathfrak g]\subset\mathfrak n$.
Define the lower central series of $\mathfrak n$ by
$\mathfrak n_1=\mathfrak n$ and
$\mathfrak n_{k+1}=[\mathfrak n,\mathfrak n_k]$.
Since $X\in\mathfrak n$, induction gives
\[
\operatorname{ad}_X^k(\mathfrak g)
\subset\mathfrak n_k,\qquad k\ge1.
\]
By definition of nilpotence, $\mathfrak n_m=0$
for some $m\ge1$. Hence $\operatorname{ad}_X^m=0$
on $\mathfrak g$, and therefore also on
$\mathfrak g_{\C}$. Since $\operatorname{ad}_X$ is diagonalizable over
$\C$, it follows that $\operatorname{ad}_X=0$.
Writing
\[
\mathfrak z(\mathfrak g)
=\{Y\in\mathfrak g:[Y,V]=0
\text{ for all }V\in\mathfrak g\}
\]
for the center of $\mathfrak g$, we conclude that
$\mathfrak n\subset\mathfrak z(\mathfrak g)$. By \cite[Corollary~1.41]{Knapp2002},
$[\mathfrak g,\mathfrak r]\subset\mathfrak n$.
Thus, for every $X\in\mathfrak r$,
\[
\operatorname{ad}_X(\mathfrak g)
\subset\mathfrak n\subset\mathfrak z(\mathfrak g),
\]
and hence $\operatorname{ad}_X^2=0$.
Since $\operatorname{ad}_X$ is diagonalizable
over $\C$, it follows that $\operatorname{ad}_X=0$.
Therefore $\mathfrak r\subset\mathfrak z(\mathfrak g)$.
Conversely, the center is an abelian ideal and
hence is contained in the radical by definition.
Thus $\mathfrak r=\mathfrak z(\mathfrak g)$.

By \cite[Theorem~B.2]{Knapp2002}, there exists a
semisimple real Lie subalgebra $\mathfrak s\subset\mathfrak g$
such that $\mathfrak g=\mathfrak s\oplus\mathfrak r$
as real vector spaces. Since
$\mathfrak r=\mathfrak z(\mathfrak g)$,
\[
[\mathfrak g,\mathfrak s]
=[\mathfrak s,\mathfrak s]\subset\mathfrak s.
\]
Thus $\mathfrak s$ is also an ideal of $\mathfrak g$,
and the decomposition is a direct sum of Lie algebras. For $X\in\mathfrak s$, the
complexification of $\operatorname{ad}_X|_{\mathfrak s}$
has spectrum in $i\R$, so its Killing form satisfies
\[
B_{\mathfrak s}(X,X)
=\tr\bigl((\operatorname{ad}_X|_{\mathfrak s})^2\bigr)
\le0.
\]
As $\mathfrak s$ is a semisimple Lie algebra, its
Killing form is nondegenerate
\cite[Theorem~2.41]{AlexandrinoBettiol2015},
and therefore negative definite. By Lie's third theorem
\cite[Theorem~1.14]{AlexandrinoBettiol2015},
there exists a connected, simply connected Lie group
$K$ with Lie algebra isomorphic to $\mathfrak s$.
Since $\mathfrak s$ is semisimple and its Killing
form is negative definite, $K$ is compact by
\cite[Theorem~2.35]{AlexandrinoBettiol2015}.

Since
$\mathfrak g=\mathfrak z(\mathfrak g)\oplus\mathfrak s$
is a direct sum of Lie algebras, the group
$\R^r\times K$, where
$r=\dim_{\R}\mathfrak z(\mathfrak g)$,
is connected and simply connected and has Lie
algebra isomorphic to $\mathfrak g$.
By \cite[Propositions~1.99, 1.100]{Knapp2002},
it is therefore the universal covering group
$\widetilde G$ of $G$. Consequently,
\cite[Proposition~1.101]{Knapp2002} gives
\[
G\simeq(\R^r\times K)/\Gamma,
\]
where $\Gamma$ is a discrete central subgroup.

Let $\Lambda$ be the projection of $\Gamma$
to $\R^r$. We claim that $\Lambda$ is discrete.
Otherwise, since $\Lambda$ is an additive subgroup,
there are nonzero $\lambda_j\in\Lambda$ with
$\lambda_j\to0$. Choose $k_j\in K$ such that
$(\lambda_j,k_j)\in\Gamma$.
By compactness of $K$, a subsequence of these
distinct elements converges in $\R^r\times K$,
contradicting the fact that $\Gamma$ is closed
and discrete. Put \(E=\operatorname{span}_{\R}\Lambda\) and choose a real basis
\(\lambda_1,\ldots,\lambda_q\) of \(E\) from \(\Lambda\).
The compact parallelepiped
\[
D=\left\{\sum_{j=1}^q t_j\lambda_j:0\le t_j\le1\right\}
\quad\text{satisfies}\quad
E=D+\sum_{j=1}^q\Z\lambda_j=D+\Lambda .
\]
Thus \(\Lambda\) is a lattice in \(E\).

Choose a real subspace $F\subset\R^r$ such that
$\R^r=F\oplus E$, and set $p=\dim_{\R}F$.
Since $\Gamma\subset E\times K$, the map
\[
F\times((E\times K)/\Gamma)
\longrightarrow(\R^r\times K)/\Gamma,
\qquad
(f,[(e,k)])\longmapsto[(f+e,k)]
\]
is a Lie-group isomorphism. Hence
\[
G\simeq\R^p\times C,
\qquad C=(E\times K)/\Gamma.
\]
The group $C$ is connected, since $E\times K$
is connected. For \((e,k)\in E\times K\), write \(e=d+\lambda\)
with \(d\in D\) and \(\lambda\in\Lambda\). If
\((\lambda,k_\lambda)\in\Gamma\), then
\([(e,k)]=[(d,kk_\lambda^{-1})]\).
Hence \(D\times K\) maps onto \(C\), so \(C\) is compact. The $\R^p$-factor is central,
since $F$ lies in the central factor $\R^r$
of $\widetilde G$.
\end{proof}

The preceding decomposition is also useful for
finding closed one-parameter subgroups of
noncompact Lie groups.

\begin{lemma}\label{stein:lem:realflow}\label{ft:lem:proper-flow}
A noncompact connected real Lie group \(G\) contains a closed
one-parameter subgroup isomorphic to \((\R,+)\).
In any real Lie group, a prescribed one-parameter subgroup is a
closed embedded real line if and only if its closure is noncompact.
When the group acts properly, these embedded lines act properly.
\end{lemma}
\begin{proof}
For a prescribed homomorphism \(t\mapsto g_t\), its closure \(H\)
is a closed connected abelian Lie group, hence isomorphic to
\(\R^a\times\mathbb T^b\) by
\cite[Corollary~1.103(b)]{Knapp2002}. Its projection to \(\R^a\)
is \(tv\) for some \(v\in\R^a\). Since the line \(\R v\) is closed and dense in \(\R^a\),
\(\R v=\R^a\). If \(a=0\) then \(H\) is compact. Otherwise
\(a=1\) and \(v\ne0\), and a compact set \(L\subset H\) has
\[
\{t:g_t\in L\}\subset
 \bigl[-\sup_{u\in\pi(L)}|u|/|v|,\ \sup_{u\in\pi(L)}|u|/|v|\bigr],
\]
where \(\pi:H\to\R^a\) is the projection. Its inverse image
is closed, hence compact. The derivative of \(\pi(g_t)=tv\)
is nonzero and the map is injective, so this is a closed
embedded real line. Conversely, an embedded closed
line has noncompact closure. A closed subgroup inherits properness of the action. If a
prescribed flow has compact closure, the orbit of any point over
all real times lies in a compact set, which rules out properness
of its \(\R\)-action.

Suppose that $\overline{\exp(\R X)}$ is compact for every
$X\in\operatorname{Lie}(G)$. For each $X$, average an inner
product on $\operatorname{Lie}(G)$ over this compact group
in the adjoint representation. Then $\exp(t\ad_X)$ is
orthogonal and $\ad_X$ is skew-adjoint for the averaged
inner product. Lemma~\ref{lem:compact-adjoint-group} gives
$G\simeq\R^p\times C$. Since $G$ is noncompact, $p>0$, and
for $0\ne v\in\R^p$ the subgroup
\[
 \{(tv,e_C):t\in\R\}
\]
has noncompact closure. This contradicts the assumption.
Thus some one-parameter subgroup has noncompact closure,
and the projection argument above makes it a closed embedded
real line. For bounded Euclidean domains, see also
\cite[Lemma~6.1]{MiebachOeljeklaus2009}.
\end{proof}

\subsection{The adjoint representation}

Let $\Omega$ be a connected complex manifold whose bounded
holomorphic functions separate tangent directions, and suppose
$G=\Aut_0(\Omega)$ is a Lie group acting properly. By
Lemma~\ref{stein:lem:bounded}, these conditions hold for
$\Omega\Subset Z$ when a neighborhood of $\overline\Omega$
admits a holomorphic immersion into Euclidean space.

The complex-linear extension of
$X\mapsto X^{1,0}$ defines a map from
$\mathfrak g_\C$ to the space of holomorphic
vector fields on $\Omega$.

We first show that this map is injective.
Suppose that
\[
X^{1,0}+\ii V^{1,0}=0,
\qquad X,V\in\mathfrak g.
\]
Then $X=-JV$. Since $V$ and $JV=-X$ are
complete real holomorphic fields, and
$[V,JV]=0$, the map
\[
s+\ii t\longmapsto
\exp(sV)\exp(tJV)(z)
\]
defines an entire curve in $\Omega$ for
every $z\in\Omega$. Every bounded
holomorphic function is constant along
this curve. Since these functions separate
tangent directions, $V$ vanishes identically,
and hence $X=0$.

For real holomorphic fields $X,V$, the mixed
brackets of their $(1,0)$- and $(0,1)$-parts
vanish in local holomorphic coordinates.
Consequently,
\[
[X^{1,0},V^{1,0}]=[X,V]^{1,0}.
\]
We therefore identify $\mathfrak g_\C$
with its image in the Lie algebra of
holomorphic vector fields on $\Omega$.

\begin{proposition}\label{prop:adjoint-branches}\label{stein:prop:adjointambient}
Let $X\in\mathfrak g$ generate a proper
real flow $\Phi_t$, and put $Z=X^{1,0}$.
\begin{enumerate}[label=\textup{(\roman*)}]
\item Suppose that a nonzero holomorphic vector field $W$
satisfies
\[
[Z,W]=\alpha W
\]
for some $\alpha\in\C$ with $\Re\alpha\ne0$.
Then, for every $K>0$, there exists a nonzero current
in $\mathcal T_{X,K}(\Omega)$ if $\Re\alpha>0$,
and in $\mathcal T_{-X,K}(\Omega)$ if $\Re\alpha<0$.

\item Suppose that $\ad_X$ has a nontrivial Jordan block
corresponding to a purely imaginary eigenvalue.
Then, for every $K>0$, there exist nonzero smooth
positive closed real $(1,1)$-forms in
$\mathcal T_{X,K}(\Omega)$ and
$\mathcal T_{-X,K}(\Omega)$.
\end{enumerate}
\end{proposition}

\begin{proof}
For \textup{(i)}, reversing the flow if necessary,
we may assume that $a:=\Re\alpha>0$.
With the convention \eqref{eq:lie-derivative-convention},
the pushforward of a holomorphic field $V$ satisfies
\begin{equation}\label{eq:adjoint-pushforward}
\frac{d}{dt}(\Phi_t)_*V=-(\Phi_t)_*[Z,V].
\end{equation}
Indeed \( (\Phi_t)_*V=\Phi_{-t}^*V\), and
\([\bar Z,V]=0\) in holomorphic coordinates. Applying
\eqref{eq:adjoint-pushforward} to \(W\) gives
\[
(\Phi_t)_*W=e^{-\alpha t}W,\qquad
D\Phi_t(z)W(z)=e^{-\alpha t}W(\Phi_t(z)).
\]
Write $\gamma_\Omega$ for the infinitesimal
Carath\'eodory--Reiffen metric,
\[
\gamma_\Omega(z;v)=
\sup_{f\in\mathcal O(\Omega,\D)}
\frac{|df_z(v)|}{1-|f(z)|^2},
\]
and set $u(z)=\gamma_\Omega(z;W(z))$;
see \cite[Propositions~10--11]{Kratz1997}. Biholomorphic
invariance and complex homogeneity give
\[
u(z)=\gamma_\Omega(\Phi_t(z);D\Phi_t(z)W(z))
    =e^{-at}u(\Phi_t(z)).
\]
Since bounded holomorphic functions separate
tangent directions, $\gamma_\Omega(z;v)>0$
for $v\ne0$. Thus $\Omega$ is
Carath\'eodory--Reiffen hyperbolic, and
\cite[Proposition~11]{Kratz1997} implies that
$\log u$ is plurisubharmonic on
$\Omega\setminus\{W=0\}$. To bound $\gamma_\Omega$ on a compact set $L$, postcompose each test map
$f:\Omega\to\D$ with the disc automorphism $A_{f(z)}$ sending
$f(z)$ to zero. Then
\[
\frac{|df_z(v)|}{1-|f(z)|^2}
=|d(A_{f(z)}\circ f)_z(v)|\le C_L|v|_h,
\]
by Cauchy's estimate on slightly larger coordinate balls,
with $h$ any fixed Hermitian metric on $\Omega$.
If \(\Omega\subset\C^n\), one can take
\(\delta=\operatorname{dist}(L,\partial\Omega)\) and restrict
to the disc in the direction \(v\) of radius \(\delta\) to get
\(\gamma_\Omega(z;v)\le |v|/\delta\) for \(z\in L\).
Thus \(\log u\) is locally bounded above
and tends to $-\infty$ at the zeros of $W$.
Since $\{W=0\}$ is a closed analytic,
hence pluripolar, subset of $\Omega$,
Theorem~5.24 of
Demailly~\cite[Chapter~I, \S5.F]{DemaillyCADG}
shows that $\log u$ extends by $-\infty$
across this set as a plurisubharmonic function. As the exponential is convex and increasing,
\[
U=\exp\left(\frac{K}{2a}\log u\right)
=u^{K/(2a)}
\]
is a nonzero, nonnegative, locally bounded
plurisubharmonic function satisfying
$U\circ\Phi_t=e^{Kt/2}U$.
Lemma~\ref{lem:positive-character}
therefore yields the required nonzero
positive closed current $\ddbar(U^2)$.

For \textup{(ii)}, let $i\lambda$,
$\lambda\in\R$, be an eigenvalue corresponding
to a nontrivial Jordan block. Choose
$W\in\mathfrak g_{\C}$ and $m\ge1$ such that
\[
V_j=(\operatorname{ad}_X-i\lambda I)^jW,
\qquad 0\le j\le m+1,
\]
satisfy $V_m\ne0$ and $V_{m+1}=0$.
Since \(V_m\) is a nonzero holomorphic vector field,
tangent separation yields a bounded holomorphic function \(h_0\)
such that \(V_m h_0\not\equiv0\).

Fix left Haar measure $dg$ on $G$. For $\eta\in C_c^\infty(G)$
with $\eta\ge0$ and $\int_G\eta\,dg=1$, put
\[
 h_\eta(z)=\int_G\eta(g)h_0(gz)\,dg .
\]
The function is bounded and holomorphic, by integration over
the fixed compact support of $\eta$. For $V\in\mathfrak g$,
write $D_V f(g)=\left.\tfrac{d}{ds}\right|_{s=0}f(g\exp(sV))$.
Thus $V_z[h_0(gz)]=D_V[h_0(\,\cdot\,z)](g)$. If $D_V^*$ is the
formal adjoint with respect to $dg$, integration by parts gives
\[
Vh_\eta(z)=
\int_G(D_V^*\eta)(g)h_0(gz)\,dg.
\]
Here
\[
D_V^*\eta=-D_V\eta-(\operatorname{div}_{dg}D_V)\eta .
\]
Consequently
$\|Vh_\eta\|_\infty\le\|D_V^*\eta\|_{L^1(G,dg)}\|h_0\|_\infty$.
On holomorphic functions, $V$ acts as $V^{1,0}$, so these
bounds extend complex linearly to $\mathfrak g_\C$.
As $\eta$ ranges over a smooth approximate
identity at the identity of $G$,
$h_\eta\to h_0$ with all derivatives
locally uniformly on $\Omega$.
We may therefore choose $\eta$ so that
$V_mh_\eta\not\equiv0$. Fix such an
$\eta$ and write $h=h_\eta$.

By \eqref{eq:adjoint-pushforward}, on $\mathfrak g_\C$ we have
\begin{equation}\label{eq:adjoint-exponential}
(\Phi_t)_*V=e^{-t\operatorname{ad}_X}V.
\end{equation}
This agrees with the abstract adjoint formula of
\cite[Proposition~1.91]{Knapp2002}: the identification
of $\operatorname{Lie}(G)$ with infinitesimal fields
for the left action reverses the bracket.
Using \eqref{eq:adjoint-exponential} and the definition of $V_j$,
we obtain
\begin{equation}\label{eq:adjoint-jordan-expansion}
\begin{aligned}
W(h\circ\Phi_t)
&=\bigl((e^{-t\operatorname{ad}_X}W)h\bigr)
  \circ\Phi_t=e^{-\ii\lambda t}
\sum_{j=0}^{m}\frac{(-t)^j}{j!}
(V_jh)\circ\Phi_t.
\end{aligned}
\end{equation}
Since $h$ is bounded on $\Omega$, local
Cauchy estimates give
\[
\sup_{t\in\R}\sup_{z\in L}
|W(h\circ\Phi_t)(z)|<\infty
\]
for every $L\Subset\Omega$.
Taking absolute values in
\eqref{eq:adjoint-jordan-expansion}
and isolating the term $j=m$, we obtain,
for $t\ne0$,
\[
\begin{aligned}
\frac{|t|^m}{m!}|(V_mh)(\Phi_t(z))|
&\le |W(h\circ\Phi_t)(z)| \quad+\sum_{j=0}^{m-1}
\frac{|t|^j}{j!}
|(V_jh)(\Phi_t(z))|.
\end{aligned}
\]
The first term is uniformly bounded for
$z\in L$ and $t\in\R$ by the local
Cauchy estimates above. The functions
$V_jh$, $0\le j<m$, are bounded on
$\Omega$ by the group-averaging argument.
Consequently,
\[
\begin{aligned}
\sup_{z\in L}|(V_mh)(\Phi_t(z))|
&\le \frac{m!}{|t|^m}
\left(
C_L+\sum_{j=0}^{m-1}
\frac{|t|^j}{j!}\|V_jh\|_\infty
\right) \longrightarrow0
\qquad(|t|\to\infty).
\end{aligned}
\]
Thus the nonzero holomorphic function
$f=V_mh$ satisfies
\eqref{eq:two-ended-holomorphic}.
Lemma~\ref{lem:two-ended-holomorphic}
produces the required currents for
both $\Phi_t$ and $\Phi_{-t}$.
\end{proof}

It remains to consider the case in which neither
adjoint construction applies to any proper real
flow. The corresponding adjoint operators are
diagonalizable with purely imaginary spectrum.

\begin{proposition}\label{prop:compact-type-structure}\label{stein:prop:compactambient}
Suppose that $G$ is noncompact and that, for every $X\in\mathfrak g$ whose flow defines a proper holomorphic $\R$-action on $\Omega$,
$\operatorname{ad}_X:\mathfrak g_{\C}\to\mathfrak g_{\C}$
is diagonalizable over $\C$ and has spectrum
contained in $i\R$. Then
\eqref{eq:central-decomposition} holds, where the $\R^p$-factor
is central and $1\le p\le n$. Its action on $\Omega$ is free
and proper, and the $(1,0)$-parts of its infinitesimal
generators are complex-linearly independent at every point.
\end{proposition}

\begin{proof}
For $X\in\mathfrak g$ let $H$ be the closure of its flow.
If $H$ is noncompact, Lemma~\ref{stein:lem:realflow} gives
properness, so $\ad_X$ has the stated spectrum and
diagonalizability by hypothesis. If $H$ is compact, average an
inner product on $\mathfrak g$ over its pushforward
representation. By \eqref{eq:adjoint-exponential},
$e^{-t\ad_X}$ is orthogonal for the resulting inner product, so
\[
\langle \ad_XV,W\rangle_H+\langle V,\ad_XW\rangle_H=0 .
\]
Thus $\ad_X$ is skew-adjoint and its complexification has the
same spectral properties. The abstract adjoint changes the
sign, so Lemma~\ref{lem:compact-adjoint-group} applies and gives
\eqref{eq:central-decomposition}, with $p\ge1$ by
noncompactness.

As observed above, $G$ acts properly on $\Omega$.
Since the $\R^p$-factor is a closed subgroup of $G$,
its action on $\Omega$ is also proper.
Its point stabilizers are compact subgroups of
$\R^p$, hence trivial. Thus the $\R^p$-action on $\Omega$ is free.
Let $X_1,\ldots,X_p$ be the complete real
holomorphic vector fields corresponding to the
standard basis of the Lie algebra of the
$\R^p$-factor, and put $Z_j=X_j^{1,0}$.
To prove their pointwise complex-linear independence,
fix $z\in\Omega$ and suppose
\[
\sum_{j=1}^p(a_j+ib_j)Z_j(z)=0,
\qquad a_j,b_j\in\R.
\]

Set $U=\sum_{j=1}^p a_jX_j$ and
$V=\sum_{j=1}^p b_jX_j$.
The assumed relation gives $U(z)=-J_zV(z)$.
Since $[U,V]=0$, their flows preserve both
vector fields. These flows are holomorphic
and hence preserve $J$. Consequently,
$U=-JV$ at every point
$\exp(sV)\exp(-tU)z$, $s,t\in\R$. Thus
\[
f(s+it)=\exp(sV)\exp(-tU)z
\]
satisfies
\[
 \frac{\partial f}{\partial t}=-U(f)=JV(f)
 =J\frac{\partial f}{\partial s}.
\]
Thus $f$ is holomorphic. Each
bounded holomorphic function is constant on its image, and
tangent separation gives $df=0$. Hence $U(z)=V(z)=0$.
Freeness of the $\R^p$-action forces all $a_j,b_j$ to vanish.
Thus $Z_1(z),\ldots,Z_p(z)$ are complex-linearly independent,
and $p\le n$.
\end{proof}

\section{Complexification of Proper Real Actions}
\label{sec:corank-complexification}

Let $A=\R^p$ act properly and holomorphically on $\Omega$.
The central factor obtained in
Proposition~\ref{prop:compact-type-structure}
provides one such action. We identify the complexification
of $A$ with $\C^p=A+\ii A$, where $a+\ii b$
corresponds to $(a,b)\in A\times A$.
The subgroup $\Z^p\subset A$ is a
discrete cocompact lattice, and
\[
\C^p/\Z^p\longrightarrow(\C^*)^p,
\qquad
[\zeta]\longmapsto
(e^{2\pi\ii\zeta_1},\ldots,
e^{2\pi\ii\zeta_p})
\]
is a biholomorphism. Since $(\C^*)^p$
is Stein, $A$ is pseudo-compact in
the sense of
\cite[Section~3]{IannuzziSpiroTrapani2004}.

The infinitesimal generators
$X_1,\ldots,X_p$ of the $A$-action
commute and are real holomorphic.
Their $(1,0)$-parts define commuting
local holomorphic flows, giving a
local $\C^p$-action on $\Omega$
extending the given $A$-action;
see \cite[Section~3]{CTIT2000}.

A holomorphic map between manifolds with local and global
$\C^p$-actions is locally equivariant if it intertwines the
actions wherever the local action is defined.

A universal globalization of the local action on $\Omega$
consists of a possibly non-Hausdorff complex manifold $P$
with a global holomorphic $\C^p$-action and a locally
equivariant open holomorphic embedding $\Omega\hookrightarrow P$.
Every locally equivariant holomorphic map from $\Omega$ to
a possibly non-Hausdorff complex manifold with a global
holomorphic $\C^p$-action must extend uniquely to a
$\C^p$-equivariant holomorphic map from $P$;
see \cite[Definition~1.3]{CTIT2000}.

If $\Omega$ is a connected Stein
manifold with a proper holomorphic
$A$-action, then
\cite[Corollary~3.5 and
Theorem~1.7]{CTIT2000}
shows that the induced local
$\C^p$-action is univalent and
admits a universal globalization
$P$, possibly non-Hausdorff.
We identify $\Omega$ with its
open image in $P$. The global
action agrees with the local
action wherever the latter
is defined, and
\[
P=\C^p\cdot\Omega.
\]
Write $\Phi_\zeta$ for the global action on $P$, so
$\Phi_{\zeta+\eta}=\Phi_\zeta\circ\Phi_\eta$ for
$\zeta,\eta\in\C^p$. We will use two criteria to establish that the universal
globalization $P$ is Hausdorff. When $\Omega$ admits a
positive definite Bergman metric, this follows from
\cite[Theorem~4.4]{IannuzziSpiroTrapani2004}.
We now give a direct proof under the alternative hypothesis
of a complete K\"ahler--Einstein metric with negative
Einstein constant. The argument also provides moment-level
slices for the complexified action on $P$.

\subsection{The K\"ahler--Einstein moment map and local slices}
\label{stein:sec:moment-slices-Stein}\label{sec:stein-structural}

Let $\Omega$ be a connected Stein manifold of complex dimension $n$,
and let $A=\R^p$, $p\ge1$, act properly and holomorphically on $\Omega$.
Write $\Phi_t$ for the action and $X_v$ for the infinitesimal generator
of $s\mapsto\Phi_{sv}$, $v\in\R^p$. Suppose in addition that $\Omega$
carries a complete K\"ahler--Einstein form $\omega_E$ normalized by
\begin{equation}\label{eq:section5-KE-normalization}
\Ric_{\omega_E}=-\omega_E.
\end{equation}
Such a metric exists for a relatively compact Stein domain whose
closure has a neighborhood admitting a holomorphic immersion into
$\C^N$ \cite[Remarks~3.6(2)]{GangulyKRsol2026Paper21}.
By uniqueness of the complete K\"ahler--Einstein metric with
normalization $\Ric_{\omega_E}=-\omega_E$, every $\Phi_t$
preserves $\omega_E$
\cite{CY80}.

For a real holomorphic generator $X=X_v$, put
$Z_X=X^{1,0}$, so $X=Z_X+\overline{Z_X}$. In local holomorphic
coordinates write
\[
\omega_E=\ii g_{E,j\bar k}\,dz^j\wedge d\bar z^k,
\qquad \rho_E=\log\det(g_{E,j\bar k}).
\]
The Einstein equation gives
\begin{equation}\label{eq:intro-KE-potential}
\ddbar\rho_E=\omega_E.
\end{equation}
Define the function $\mu_X^E$ by
\begin{equation}\label{eq:KE-canonical-moment}
\mu_X^E
=\ii\,\frac{\Lie_{Z_X}\omega_E^n}{\omega_E^n}
=\ii\left(\operatorname{div}_{\C}Z_X+Z_X\rho_E\right),
\qquad
\operatorname{div}_{\C}Z_X=\sum_{j=1}^n\partial_j Z_X^j.
\end{equation}
To verify the coordinate formula, let $\nu_E=\omega_E^n$. Since
$\Lie_{Z_X}\nu_E=\bigl(\sum_j\partial_j Z_X^j+Z_X\rho_E\bigr)\nu_E$
and $\Lie_X\nu_E=0$, the coefficient in parentheses is purely
imaginary. Thus $\mu_X^E$ is real. The divergence of $Z_X$ is
holomorphic; consequently,
\[
\bar\partial\mu_X^E
=\ii\sum_{j,k} Z_X^j(\rho_E)_{j\bar k}\,d\bar z^k
=\ii\sum_{j,k} Z_X^j g_{E,j\bar k}\,d\bar z^k
=(\iota_X\omega_E)^{0,1}.
\]
Both $d\mu_X^E$ and $\iota_X\omega_E$ are real one-forms.
Hence, for $X_j=X_{e_j}$,
\begin{equation}\label{stein:eq:intrinsic-moment-identity}
 d\mu_X^E=\iota_X\omega_E,
 \qquad d\mu_{X_j}^E(JX_k)=g_E(X_j,X_k).
\end{equation}
Since $A$ is abelian and preserves $\nu_E$,
naturality of the Lie derivative gives
\begin{equation}\label{eq:central-moment-invariance}
\mu_{X_j}^E\circ\Phi_t=\mu_{X_j}^E\qquad(t\in A).
\end{equation}
The Einstein normalization also yields a Laplacian identity not
available for a general invariant K\"ahler metric. By K\"ahler
symmetry $\partial_i g_{E,j\bar k}=\partial_jg_{E,i\bar k}$,
\begin{equation}\label{eq:KE-moment-laplacian}
\begin{aligned}
\Delta_E\mu_X^E
&=\ii g_E^{i\bar k}\bigl((\partial_i Z_X^j)g_{E,j\bar k}
       +Z_X^j\partial_i g_{E,j\bar k}\bigr) =\ii\bigl(\operatorname{div}_{\C}Z_X+Z_X\rho_E\bigr)
=\mu_X^E.
\end{aligned}
\end{equation}
For $\Ric_{\omega_E}=-\kappa\omega_E$, $\kappa>0$, the
moment normalized by $\omega_E$ is
$\mu_X^{E,\kappa}=\ii\kappa^{-1}(\Lie_{Z_X}\omega_E^n)/\omega_E^n$;
the same calculation gives
$d\mu_X^{E,\kappa}=\iota_X\omega_E$ and
$\Delta_E\mu_X^{E,\kappa}=\kappa\mu_X^{E,\kappa}$.

\begin{proposition}
\label{stein:prop:ambient-moment-slices}\label{prop:moment-slices}
Let $\Omega$ be a connected Stein manifold with a complete
K\"ahler--Einstein metric of negative Einstein constant, and let
$A=\R^p$ act properly and holomorphically on $\Omega$.
Rescale the metric to $\omega_E$ with $\Ric_{\omega_E}=-\omega_E$,
and use the moments in \eqref{eq:KE-canonical-moment}.
Then the universal globalization $P=\C^p\cdot\Omega$ is a connected
Hausdorff complex manifold and the following statements hold.
\begin{enumerate}[label=\textup{(\roman*)},leftmargin=0pt,
                  labelindent=0pt,labelsep=0.45em,itemindent=*,
                  widest=ii,align=left]

\item For every $x\in\Omega$, put
\begin{equation}\label{eq:imaginary-return-set}
D_x=\{\eta\in\R^p:\Phi_{\ii\eta}(x)\in\Omega\}.
\end{equation}
Then $D_x$ is open and convex, and
\[
\{\zeta\in\C^p:\Phi_\zeta(x)\in\Omega\}=A+\ii D_x.
\]
\item For $c\in\mu^E(\Omega)$ put
$N_c=(\mu^E)^{-1}(c)$, where
$\mu^E=(\mu_{X_1}^E,\ldots,\mu_{X_p}^E)$.
Then $N_c$ is a closed smooth $A$-invariant submanifold of
$\Omega$. The map
\[
\mathcal F_c:\R^p\times N_c\longrightarrow P,
\qquad (\eta,x)\longmapsto\Phi_{\ii\eta}(x)
\]
is a diffeomorphism onto the open invariant subset
$U_c=\C^p\cdot N_c$, and the action there is given by
\[
(a+\ii b)\cdot(\eta,x)=(\eta+b,\Phi_a(x)).
\]
The sets $U_c$ cover $P$.
\end{enumerate}
In particular, the complexified action is free and locally proper
(that is, every point has a $\C^p$-invariant open neighborhood
on which the restricted action is proper), and $p\le n$.
\end{proposition}

\begin{proof}
The real action is free, because its stabilizers are compact
subgroups of $\R^p$. Orbit-convexity
\cite[Corollary~3.6]{CTIT2000} holds in the possibly
non-Hausdorff globalization: for $x\in\Omega$ and
$b\in\R^p$, the set of $s\in\R$ for which
$\Phi_{\ii sb}(x)\in\Omega$ is an interval. If
$\eta_0,\eta_1\in D_x$, apply this assertion at
$\Phi_{\ii\eta_0}(x)$ in the direction $\eta_1-\eta_0$.
The entire segment from $\eta_0$ to $\eta_1$ belongs to $D_x$.
Thus $D_x$ is convex. Its openness follows from $\Omega$
being open in $P$, while real invariance gives the equality
in \textup{(i)}. This proves (i) without using Hausdorffness.

We prove that $P$ is Hausdorff using completeness of $g_E$.
Fix $\zeta=a+\ii b\in\C^p$, and suppose that
\[
 x_j\longrightarrow x,\qquad
 y_j=\Phi_\zeta(x_j)\longrightarrow y
 \qquad\text{in }\Omega.
\]
Put $w_j=\Phi_{-a}(y_j)$. Then $w_j=\Phi_{\ii b}(x_j)$
and $w_j\to\Phi_{-a}(y)=:w$ in $\Omega$. By (i),
\[
\gamma_j(s)=\Phi_{\ii sb}(x_j)\in\Omega,
\qquad 0\le s\le1.
\]
Holomorphicity of the global action gives
$\dot\gamma_j=JX_b(\gamma_j)$. From
\eqref{stein:eq:intrinsic-moment-identity},
\begin{equation}\label{eq:KE-imaginary-segment-energy}
\int_0^1|\dot\gamma_j(s)|_{g_E}^2\,ds
=\langle\mu^E(w_j)-\mu^E(x_j),b\rangle.
\end{equation}
The right-hand side is bounded independently of $j$.
For $0\le s<t\le1$, Cauchy--Schwarz yields
\[
 d_{g_E}(\gamma_j(s),\gamma_j(t))
 \le |t-s|^{1/2}
 \left(\int_0^1|\dot\gamma_j|_{g_E}^2\,du\right)^{1/2}.
\]
The length estimate places all $\gamma_j([0,1])$ in a fixed
closed $g_E$-ball about $x$. This ball is compact by completeness
and Hopf--Rinow. Arzel\`a--Ascoli gives a uniformly convergent
subsequence with limit $\gamma:[0,1]\to\Omega$,
$\gamma(0)=x$, $\gamma(1)=w$. On a small time interval in a
coordinate chart, write $u_j$ for the coordinates of $\gamma_j$
and $V_b$ for the coordinate expression of $JX_b$. There
\[
 u_j(t)-u_j(s)=\int_s^t V_b(u_j(r))\,dr .
\]
The curves on that interval lie in a compact subset of the
chart for large $j$, so the integrands converge uniformly.
It follows that $\dot\gamma=JX_b(\gamma)$ on $[0,1]$.
For each $s\in[0,1]$ and sufficiently small $h$ with
$s+h\in[0,1]$, uniqueness for the equation
$\dot\gamma=JX_b(\gamma)$ within $\Omega$ gives
\[
\gamma(s+h)=\Phi_{\ii hb}(\gamma(s)),\qquad
\Phi_{-\ii(s+h)b}(\gamma(s+h))
 =\Phi_{-\ii sb}(\gamma(s)).
\]
Thus $s\mapsto\Phi_{-\ii sb}(\gamma(s))$ is locally
constant in $P$ and hence constant on $[0,1]$.
In particular $w=\Phi_{\ii b}(x)$ and $y=\Phi_\zeta(x)$.
We have proved that
\begin{equation}\label{eq:section5-fixed-time-closed}
\{(u,v)\in\Omega\times\Omega:v=\Phi_\zeta(u)\}
\quad\text{is closed for every fixed }\zeta\in\C^p.
\end{equation}
For $\alpha,\beta\in\C^p$, the diagonal of $P\times P$
in the open product of charts
$\Phi_\alpha(\Omega)\times\Phi_\beta(\Omega)$ is precisely
the relation $y=\Phi_{\alpha-\beta}(x)$ in
$\Omega\times\Omega$. The products of these translated charts
cover $P\times P$; by \eqref{eq:section5-fixed-time-closed}
the diagonal is closed. Thus $P$ is Hausdorff. Connectedness
follows from $P=\C^p\cdot\Omega$.

For (ii), if $u\in\R^p\setminus\{0\}$ then $X_u$ is nowhere
zero, by freeness. The Gram identity
\begin{equation}\label{eq:moment-gram}
\sum_{j,k=1}^p u_ju_k\,d\mu_{X_j}^E(JX_k)
=g_E(X_u,X_u)>0
\qquad(u\ne0)
\end{equation}
shows that $\mu^E$ is a submersion. The level $N_c$ is therefore
smooth, closed in $\Omega$, and $A$-invariant by
\eqref{eq:central-moment-invariance}. The proper real action
restricts to $N_c$.

At $(0,x)\in\R^p\times N_c$,
\[
 d\mathcal F_c(u,W)=\sum_{k=1}^p u_k JX_k(x)+W,
 \qquad W\in\ker d\mu_x^E.
\]
If $G_x=(g_E(X_j,X_k)(x))_{j,k}$, then any
$V\in T_x\Omega$ has the unique decomposition
\[
 u=G_x^{-1}d\mu_x^E(V),\qquad
 W=V-\sum_k u_k JX_k(x)\in\ker d\mu_x^E.
\]
Thus $d\mathcal F_c$ is an isomorphism at $(0,x)$. For
$\eta,\xi\in\R^p$ and $y\in N_c$, the group law gives
$\mathcal F_c(\eta+\xi,y)=\Phi_{\ii\eta}(\mathcal F_c(\xi,y))$.
Differentiating at $(\xi,y)=(0,x)$ yields
\[
d\mathcal F_c|_{(\eta,x)}(u,W)
=d(\Phi_{\ii\eta})_x\bigl(d\mathcal F_c|_{(0,x)}(u,W)\bigr).
\]
Hence $\mathcal F_c$ is a local diffeomorphism.

Suppose $\mathcal F_c(\eta,x)=\mathcal F_c(\eta',y)$.
With $\xi=\eta-\eta'$ the group law gives
$y=\Phi_{\ii\xi}(x)$. By (i), the imaginary segment from
$x$ to $y$ lies in $\Omega$. If $\xi\ne0$, then
\begin{equation}\label{eq:KE-moment-increase}
\langle\mu^E(y)-\mu^E(x),\xi\rangle
=\int_0^1|X_\xi|_{g_E}^2(\Phi_{\ii s\xi}(x))\,ds>0,
\end{equation}
contrary to $x,y\in N_c$. Hence $\xi=0$, and $x=y$.
Thus $\mathcal F_c$ is an injective local diffeomorphism onto
$U_c$. Since $N_c$ is invariant under $A$, the group law gives
the stated action formula. The real action is proper on $N_c$,
and the imaginary parameters translate the first factor;
therefore the $\C^p$-action is free and proper on $U_c$.
These sets cover $P$ because every point is a complex translate
of some $x\in N_{\mu^E(x)}$. The complexified action is
therefore free and locally proper.

Finally, real invariance gives
$\omega_E(X_j,X_k)=d\mu_{X_j}^E(X_k)=0$.
The real span of the $X_j$ is $g_E$-orthogonal to the real span
of the $JX_j$, and both have dimension $p$ by
\eqref{eq:moment-gram}. Thus $2p\le2n$.
\end{proof}

For a proper real flow with generator $X$, the moment $\mu_X^E$
descends to $\Omega/\R$. Its derivative on the projection of $JX$
is $|X|_{g_E}^2>0$, so the quotient cannot be compact.
For the vector action, put
$m_j(\eta)=\mu_{X_j}^E(\Phi_{\ii\eta}(x))$ on $D_x$.
Then
\[
\frac{\partial m_j}{\partial\eta_k}
=g_E(X_j,X_k)(\Phi_{\ii\eta}(x))
=\frac{\partial m_k}{\partial\eta_j}.
\]
Put $\mathcal U_x(\eta)=\int_0^1\sum_k\eta_km_k(t\eta)\,dt$.
The symmetry above gives
\[
 \partial_{\eta_j}\mathcal U_x
 =\int_0^1\left(m_j(t\eta)+t\sum_k\eta_k
          (\partial_{\eta_k}m_j)(t\eta)\right)\,dt
 =m_j(\eta).
\]
The Hessian of $\mathcal U_x$ is the positive definite Gram matrix
at $\Phi_{\ii\eta}(x)$, so $\mathcal U_x$ is strictly convex on $D_x$.

\begin{lemma}\label{lem:holomorphic-principal-quotient}
A free proper holomorphic action of a complex Lie group $H$
on a complex manifold $X$ has a Hausdorff complex-manifold
quotient and is a holomorphic principal bundle over it.
\end{lemma}

\begin{proof}
Fix $x\in X$ and consider the holomorphic orbit map
\[
\alpha_x:H\longrightarrow X,\qquad h\longmapsto h\cdot x.
\]
Its differential at the identity $e\in H$ is injective.
Indeed, if $v\in T_eH$ belongs to its kernel, the
infinitesimal vector field generated by $v$ vanishes at
$x$. Uniqueness for its integral curves implies that
$\exp(tv)\cdot x=x$ for every $t\in\R$. By freeness,
$\exp(tv)=e$ for all $t$, and hence $v=0$.

Choose a complex submanifold $S_0$ through $x$ such that
\[
T_xX=\operatorname{im}(d\alpha_x)_e\oplus T_xS_0.
\]
Consider the holomorphic map
\[
F:H\times S_0\longrightarrow X,\qquad
F(h,s)=h\cdot s.
\]
Its differential at $(e,x)$ is an isomorphism. By the
holomorphic inverse function theorem, there exist an
identity neighborhood $N\subset H$ and a neighborhood
$S_1$ of $x$ in $S_0$ such that
\[
F:N\times S_1\longrightarrow F(N\times S_1)
\]
is a biholomorphism onto an open subset of $X$.

We claim that $S_1$ can be replaced by a smaller
neighborhood $S$ of $x$ such that
\begin{equation}\label{eq:principal-slice-overlap}
\{h\in H:hS\cap S\ne\varnothing\}\subset N.
\end{equation}
Otherwise, by choosing neighborhoods of $x$ in $S_1$
that shrink to $x$, we obtain sequences
$s_j,s'_j\to x$ and $h_j\in H\setminus N$ satisfying
\[
h_j\cdot s_j=s'_j.
\]
Since the action map
\[
H\times X\longrightarrow X\times X,\qquad
(h,z)\longmapsto(z,h\cdot z)
\]
is proper and $(s_j,s'_j)\to(x,x)$, a subsequence
of $h_j$ converges to some $h\in H\setminus N$.
Passing to the limit gives $h\cdot x=x$, so freeness
implies $h=e$, contradicting $e\in N$.
This proves \eqref{eq:principal-slice-overlap}.

We show that
\[
F:H\times S\longrightarrow HS,\qquad
HS=\{h\cdot s:h\in H,\ s\in S\},
\]
is a biholomorphism. Suppose
$F(h,s)=F(h',s')$. Then
\[
(h'^{-1}h)\cdot s=s',
\]
so \eqref{eq:principal-slice-overlap} gives
$h'^{-1}h\in N$. Since $F$ is injective on
$N\times S$, comparison of
\[
F(h'^{-1}h,s)=F(e,s')
\]
shows that $h=h'$ and $s=s'$.
Moreover, $F$ is a local biholomorphism at every
point of $H\times S$: this holds on $N\times S$
and extends to every $h\in H$ by left translation.
Thus $HS$ is open in $X$, and
$F:H\times S\to HS$ is a biholomorphism.

Let $Q=X/H$, and let $\pi:X\to Q$ be the quotient
map. The orbit relation
\[
R=\{(z,h\cdot z):z\in X,\ h\in H\}
\]
is closed in $X\times X$, since it is the image
of the proper action map. If $x,y\in X$ belong
to distinct orbits, there are open neighborhoods
$U$ of $x$ and $V$ of $y$ such that
$(U\times V)\cap R=\varnothing$.
The sets $\pi(U)$ and $\pi(V)$ are disjoint.
They are open because, for every open $W\subset X$,
\[
\pi^{-1}(\pi(W))=\bigcup_{h\in H}hW
\]
is open. Hence $Q$ is Hausdorff. Since $\pi$ is
an open surjection and $X$ is second countable,
$Q$ is also second countable.

For each slice $S$ constructed above, put
$W=\pi(HS)$. The set $W$ is open in $Q$, and
the biholomorphism $F:H\times S\to HS$ identifies
$\pi|_{HS}$ with the projection
\[
H\times S\longrightarrow S,\qquad (h,s)\longmapsto s.
\]
Consequently, $\pi|_S:S\to W$ is a homeomorphism.
These maps provide complex coordinates on $Q$.

To verify compatibility, let $S_i,S_j$ be two
such slices, with corresponding biholomorphisms
$F_i:H\times S_i\to HS_i$ and
$F_j:H\times S_j\to HS_j$. On
$S_i\cap HS_j$, the map
\[
F_j^{-1}|_{S_i\cap HS_j}:
S_i\cap HS_j\longrightarrow H\times S_j
\]
is holomorphic. Its $S_j$-component is precisely
the change of quotient coordinates. Thus the
charts define a complex-manifold structure on
$Q$, for which $\pi$ is a holomorphic submersion.

Finally, the inverse maps
\[
\sigma_i=(\pi|_{S_i})^{-1}:W_i\longrightarrow X
\]
are holomorphic local sections. The maps
\[
H\times W_i\longrightarrow\pi^{-1}(W_i),
\qquad (h,q)\longmapsto h\cdot\sigma_i(q)
\]
are biholomorphic and $H$-equivariant.
On $W_i\cap W_j$, the unique functions
$g_{ij}:W_i\cap W_j\to H$ satisfying
\[
\sigma_i(q)=g_{ij}(q)\cdot\sigma_j(q)
\]
are holomorphic, since they are the
$H$-components of $F_j^{-1}\circ\sigma_i$.
These are the holomorphic principal-bundle
transition functions. Therefore
$\pi:X\to Q$ is a holomorphic principal
$H$-bundle.
\end{proof}

\subsection{Fourier transforms of square-integrable holomorphic
\texorpdfstring{$n$}{n}-forms}

\begin{proposition}
\label{stein:prop:Fourier}\label{lem:spectral-Fourier}
Under the hypotheses of Lemma~\ref{stein:lem:intrinsic-mass}, for every
\(\theta\in\R^p\), the integral over the action parameter
\begin{equation}\label{stein:eq:Fourier}
S_{\sigma,\theta}(z)=\int_{\R^p}e^{-\ii\langle\theta,t\rangle}
  \bigl(\Phi_t^*\sigma\bigr)^{\otimes2}(z)\,dt
\end{equation}
defines \(S_{\sigma,\theta}\in H^0(M,K_M^{\otimes2})\). The
integral converges absolutely and locally uniformly in
holomorphic trivializations, and
\(\theta\mapsto S_{\sigma,\theta}\) is continuous in the
compact-open topology. Moreover,
\begin{equation}\label{stein:eq:real-character}
\Phi_u^*S_{\sigma,\theta}
=e^{\ii\langle\theta,u\rangle}S_{\sigma,\theta},
\qquad u\in\R^p.
\end{equation}
If \(\sigma(x)\ne0\), then \(S_{\sigma,\theta}(x)\ne0\) for at least one
\(\theta\in\R^p\).
\end{proposition}
\begin{proof}
On a relatively compact coordinate chart, write
\(\Phi_t^*\sigma=a_t(z)\,dz^1\wedge\cdots\wedge dz^n\). The norm of the
integrand in \eqref{stein:eq:Fourier} is \(|a_t(z)|^2\) in the induced frame of
\(K_M^{\otimes2}\). By \eqref{stein:eq:derivative-mass} with \(m=0\),
\(\int_{\R^p}\sup_Q|a_t(z)|^2dt<\infty\) for every compact \(Q\)
in this chart. Thus \eqref{stein:eq:Fourier} converges locally uniformly
and its coefficients are holomorphic. The same bound and dominated
convergence give compact-open continuity in \(\theta\).
By \eqref{eq:real-action-group-law} and the change of variables $v=t+u$,
\begin{align*}
\Phi_u^*S_{\sigma,\theta}
&=\int_{\R^p}e^{-\ii\langle\theta,t\rangle}
             (\Phi_{t+u}^*\sigma)^{\otimes2}\,dt\\
&=e^{\ii\langle\theta,u\rangle}
  \int_{\R^p}e^{-\ii\langle\theta,v\rangle}
             (\Phi_v^*\sigma)^{\otimes2}\,dv
=e^{\ii\langle\theta,u\rangle}S_{\sigma,\theta}.
\end{align*}
This is \eqref{stein:eq:real-character}.
Fix \(x\in M\) with \(\sigma(x)\ne0\), and choose a holomorphic frame
\(e\) of \(K_M\) near \(x\). Write
\((\Phi_t^*\sigma)(x)=a_x(t)e(x)\). Then
\(t\mapsto a_x(t)^2\) belongs to \(L^1(\R^p)\) by
Lemma~\ref{stein:lem:intrinsic-mass}, is continuous, and satisfies
\(a_x(0)^2\ne0\). If every \(S_{\sigma,\theta}(x)\) vanished, the
Fourier transform of this \(L^1\) function would vanish at every
real frequency. Fourier uniqueness \cite[Chapter~I, Corollary~1.22]{SteinWeiss1971} would give \(a_x(t)^2=0\)
almost everywhere, hence everywhere by continuity, a contradiction.
\end{proof}

For a bounded domain \(\Omega\subset\C^n\), the coordinate frame
\(\sigma_0=dz^1\wedge\cdots\wedge dz^n\) is square-integrable.
From \eqref{eq:shared-jacobian-L2},
\begin{equation}\label{eq:spectral-finite-mass}
\int_{\R^p}\sup_{z\in Q}|J_t(z)|^2\,dt<\infty\qquad(Q\Subset\Omega).
\end{equation}
The section in \eqref{stein:eq:Fourier} is then
\begin{equation}\label{eq:scalar-fourier-jacobian}
S_{\sigma_0,\theta}(z)=
\left(\int_{\R^p}e^{-\ii\langle\theta,t\rangle}J_t(z)^2\,dt\right)
\sigma_0^{\otimes2}.
\end{equation}
The square in \eqref{eq:scalar-fourier-jacobian} is the holomorphic
square. Since \(J_0(z)=1\), Proposition~\ref{stein:prop:Fourier}
shows that for every $z\in\Omega$, at least one of the
sections $S_{\sigma_0,\theta}$ is nonzero at $z$.

\begin{corollary}\label{stein:cor:finite}
Let $\Omega\Subset Z$ be a domain, and let $V\subset Z$ be
an open neighborhood of $\overline\Omega$. Suppose that
$V$ is Stein or admits a holomorphic immersion into $\C^N$.
Then there exist finitely many holomorphic $n$-forms
$\sigma_1,\ldots,\sigma_m\in H^0(V,K_V)$ with the following
property: for every proper holomorphic $\R^p$-action
$\Phi$ on $\Omega$, the Fourier integrals
\begin{equation}\label{stein:eq:finite-family}
S_{j,\theta}
=\int_{\R^p}e^{-\ii\langle\theta,t\rangle}
\bigl(\Phi_t^*(\sigma_j|_\Omega)\bigr)^{\otimes2}\,dt,
\qquad \theta\in\R^p,
\end{equation}
define holomorphic sections of $K_\Omega^{\otimes2}$,
depending continuously on $\theta$ in the compact-open
topology and satisfying \eqref{stein:eq:real-character}.
Moreover, the family $\{S_{j,\theta}:1\le j\le m,\,
\theta\in\R^p\}$ has no common zero on $\Omega$.
\end{corollary}

\begin{proof}
We first choose finitely many holomorphic $n$-forms on
$V$ having no common zero on $\overline\Omega$.
If $F:V\to\C^N$ is a holomorphic immersion, its
coordinate differentials span the holomorphic cotangent
space at every point. Hence the forms
\[
dF_{j_1}\wedge\cdots\wedge dF_{j_n},
\qquad 1\le j_1<\cdots<j_n\le N,
\]
have no common zero on $V$.
If $V$ is Stein, Cartan's Theorem A
\cite[Chapter~IV]{GrauertRemmert2004}
applied to $K_V$ gives, at each point of
$\overline\Omega$, a global holomorphic section
that does not vanish there. Compactness of
$\overline\Omega$ then gives the required finite
family in this case as well.

Each $\sigma_j|_\Omega$ is square-integrable,
since $\sigma_j$ is holomorphic on a neighborhood
of the compact set $\overline\Omega$.
Proposition~\ref{stein:prop:Fourier} therefore
applies to each $\sigma_j|_\Omega$ and gives
the asserted convergence, holomorphicity,
continuity, and transformation law.

Finally, fix $x\in\Omega$. By construction,
$\sigma_j(x)\ne0$ for some $j$. The nonvanishing
assertion of Proposition~\ref{stein:prop:Fourier}
then gives $\theta\in\R^p$ such that
$S_{j,\theta}(x)\ne0$. Thus the sections
$S_{j,\theta}$ have no common zero on $\Omega$.
\end{proof}

We next extend the real transformation law
\eqref{stein:eq:real-character} to the complexified action. The connectedness assumption in the following lemma ensures
that the extensions agree on overlapping translates of $M$.
\begin{lemma}\label{stein:lem:extension}\label{lem:spectral-extend}
Let \(P\) be a Hausdorff complex manifold with a holomorphic
\(\C^p\)-action \(\Phi_\zeta\). Suppose \(M\subset P\) is open
and invariant under the real subgroup \(\R^p\), which acts
properly on \(M\), and \(P=\C^p\cdot M\). Assume that for every
\(x\in M\) the set
\[
T_x=\{\zeta\in\C^p:\Phi_\zeta(x)\in M\}
\]
is connected. For each square-integrable holomorphic $n$-form
$\sigma$ on $M$, the sections $S_{\sigma,\theta}$ of
Proposition~\ref{stein:prop:Fourier} extend uniquely to
holomorphic sections of \(K_P^{\otimes2}\). These extensions
depend continuously on $\theta$ and satisfy
\begin{equation}\label{stein:eq:complex-character}
\Phi_\zeta^*S_{\sigma,\theta}
=e^{\ii\langle\theta,\zeta\rangle}S_{\sigma,\theta},
\qquad \zeta\in\C^p,
\end{equation}
where the pairing is the complex-bilinear extension of the usual pairing.
If the sections $S_{\sigma,\theta}$, with $\sigma$ ranging over
a collection of square-integrable holomorphic $n$-forms and
$\theta\in\R^p$, have no common zero on $M$, their extensions
have no common zero on $P$.
\end{lemma}
\begin{proof}
For $\zeta\in\C^p$ and $x\in P$, write
$(\Phi_\zeta)^*_x:(K_P^{\otimes2})_{\Phi_\zeta(x)}
\to (K_P^{\otimes2})_x$ for pullback by $D_x\Phi_\zeta$.
Since $\Phi_{\zeta+\eta}=\Phi_\eta\circ\Phi_\zeta$, we get
\begin{equation}\label{eq:canonical-fiber-composition}
(\Phi_{\zeta+\eta})^*_x
=(\Phi_\zeta)^*_x\circ(\Phi_\eta)^*_{\Phi_\zeta(x)} .
\end{equation}
Fix $x\in M$ and $\theta\in\R^p$. The expression
\[
F_{\theta,x}(\zeta)=(\Phi_\zeta)^*_x
\bigl(S_{\sigma,\theta}(\Phi_\zeta(x))\bigr)
-e^{\ii\langle\theta,\zeta\rangle}S_{\sigma,\theta}(x),
\qquad \zeta\in T_x,
\]
is holomorphic with values in the fixed fiber
$(K_P^{\otimes2})_x$. By \eqref{stein:eq:real-character}, it
vanishes on $\R^p\subset T_x$. Choose a polydisc centered
at zero in the open set $T_x$. The one-variable identity
theorem, applied successively in each coordinate, gives
$F_{\theta,x}=0$ on this polydisc. Since $T_x$ is connected,
$F_{\theta,x}$ vanishes on all of $T_x$. In particular,
\begin{equation}\label{eq:spectral-local-complex-equivariance}
(\Phi_\zeta)^*_x S_{\sigma,\theta}(\Phi_\zeta x)
=e^{\ii\langle\theta,\zeta\rangle} S_{\sigma,\theta}(x)
\qquad(\zeta\in T_x).
\end{equation}
For \(\eta\in\C^p\) set \(U_\eta=\Phi_\eta(M)\). These open sets
cover \(P\). On \(U_\eta\) put
\[
S_{\sigma,\theta}^{(\eta)}
=e^{\ii\langle\theta,\eta\rangle}
(\Phi_{-\eta}|_{U_\eta})^*S_{\sigma,\theta}.
\]
If \(q\in U_\eta\cap U_\xi\), put \(x=\Phi_{-\eta}(q)\)
and \(y=\Phi_{-\xi}(q)=\Phi_{\eta-\xi}(x)\).
Then $\eta-\xi\in T_x$. At \(x\) we have
\[
(\Phi_\eta)^*_x S_{\sigma,\theta}^{(\eta)}(q)
=e^{\ii\langle\theta,\eta\rangle} S_{\sigma,\theta}(x).
\]
For the other local section, \eqref{eq:canonical-fiber-composition}
and \eqref{eq:spectral-local-complex-equivariance} give
\begin{align*}
(\Phi_\eta)^*_x S_{\sigma,\theta}^{(\xi)}(q)
&=e^{\ii\langle\theta,\xi\rangle}
  (\Phi_{\eta-\xi})^*_x S_{\sigma,\theta}(y)\\
&=e^{\ii\langle\theta,\eta\rangle}S_{\sigma,\theta}(x).
\end{align*}
Since $(\Phi_\eta)^*_x$ is an isomorphism, the local sections
agree on overlaps and define a holomorphic section on $P$.
For $\eta=0$ the local definition is the original section on $M$. For \(\gamma\in\C^p\), \(\Phi_\gamma(U_\eta)=U_{\eta+\gamma}\).
On \(U_\eta\) the pullback of the glued section is
\[
e^{\ii\langle\theta,\eta+\gamma\rangle}
  (\Phi_{-(\eta+\gamma)}\circ\Phi_\gamma)^*S_{\sigma,\theta}
  =e^{\ii\langle\theta,\gamma\rangle}S_{\sigma,\theta}^{(\eta)} .
\]
This proves \eqref{stein:eq:complex-character}. For a compact
$Q\subset P$, choose finitely many compact sets
$Q_i\subset U_{\eta_i}$ whose union contains $Q$. On $Q_i$
the section is a fixed pullback of $S_{\sigma,\theta}$ from
the compact set $\Phi_{-\eta_i}(Q_i)$, multiplied by
$e^{\ii\langle\theta,\eta_i\rangle}$. This proves compact-open
continuity in $\theta$. Uniqueness follows from
\eqref{stein:eq:complex-character} and $P=\C^p\cdot M$. For $q\in P$, write $q=\Phi_\eta(x)$ with $x\in M$ and choose
a section with $S_{\sigma,\theta}(x)\ne0$. By
\eqref{stein:eq:complex-character},
$(\Phi_\eta)^*_x S_{\sigma,\theta}(q)
=e^{\ii\langle\theta,\eta\rangle}S_{\sigma,\theta}(x)\ne0$.
Thus the extended sections have no common zero.
\end{proof}

\subsection{Properness of the complexified action}

We apply the following lemma to the extended Fourier sections
$S_{\sigma,\theta}$. Ratios of sections with nearby frequencies
transform by characters whose exponents are the differences
of those frequencies. This gives freeness and properness on
invariant open neighborhoods. When $\dim_\C P=p+1$, dividing
out the minimum vanishing order along a transverse disc
extends these ratios across their common zeros and allows
the properness argument to be applied globally.

\begin{lemma}
\label{stein:lem:spectral-several}\label{lem:spectral-abstract}
Let \(P\) be a connected Hausdorff complex manifold with a
holomorphic \(\C^p\)-action. Let \(L\to P\) be a holomorphic
line bundle with a holomorphic \(\C^p\)-action lifting the
action on $P$ and acting linearly on each fiber.
Let $J$ be any index set, and let
\(\{s_{j,\theta}\}_{j\in J,\theta\in\R^p}\)
be a collection of holomorphic sections of \(L\), with no common zero,
such that for each fixed \(j\), the map \(\theta\mapsto s_{j,\theta}\)
is continuous in the compact-open topology and
\begin{equation}\label{stein:eq:lift}
\Phi_\zeta^*s_{j,\theta}
=e^{\ii\langle\theta,\zeta\rangle}s_{j,\theta}
\quad(\zeta\in\C^p).
\end{equation}
Then the \(\C^p\)-action on \(P\) is free and locally proper.
If $\dim_\C P=p+1$, the action is globally proper.

\end{lemma}
\begin{proof}
We first record an elementary consequence of convergence in the
frequency. Suppose that on a cube $I$ centered at $\theta_0$ the
functions
\[
E_k(\theta)=e^{\ii\langle\theta-\theta_0,\zeta_k\rangle}
\]
converge uniformly to a continuous $H$, with $H(\theta_0)=1$.
After shrinking to a smaller cube about $\theta_0$, we have
$|H(\theta)-1|<1/4$ and, for all large $k$,
$|E_k(\theta)-1|<1/2$. The principal logarithm is then defined,
and
\begin{equation}\label{eq:spectral-logarithmic-parameter}
\log E_k(\theta)=\ii\langle\theta-\theta_0,\zeta_k\rangle.
\end{equation}
The difference of the two sides is continuous with values in
$2\pi\ii\Z$ and vanishes at $\theta_0$. For a fixed
$\varepsilon>0$ with $\theta_0+\varepsilon e_\ell$ in the smaller
cube for each $\ell$, uniform convergence gives
\[
 \zeta_{k,\ell}
 =\frac{\log E_k(\theta_0+\varepsilon e_\ell)}{\ii\varepsilon}
 \longrightarrow
 \frac{\log H(\theta_0+\varepsilon e_\ell)}{\ii\varepsilon}.
\]

Fix $q\in P$. Choose an index $j$ and a frequency $\theta_0$
with $s_{j,\theta_0}(q)\ne0$. Compact-open continuity gives a
closed cube $I$ centered at $\theta_0$ such that
$s_{j,\theta}(q)\ne0$ for all $\theta\in I$. Set
\[
V=\{x\in P:s_{j,\theta}(x)\ne0\text{ for every }\theta\in I\}.
\]
This is an open neighborhood of $q$: in a local trivialization,
the continuous coefficient has a positive minimum over $I$,
and the same holds on a smaller neighborhood of $q$. The
character law \eqref{stein:eq:lift} makes $V$ invariant under
$\C^p$. For $\theta\in I$ the ratios
$r_\theta=s_{j,\theta}/s_{j,\theta_0}$ are holomorphic and
nowhere zero on $V$, and satisfy
\begin{equation}\label{eq:spectral-local-ratios}
 r_\theta(\Phi_\zeta x)
 =e^{\ii\langle\theta-\theta_0,\zeta\rangle}r_\theta(x).
\end{equation}
If $\Phi_\zeta(q)=q$, variation of $\theta$ on $I$ gives
$e^{\ii\langle\theta-\theta_0,\zeta\rangle}=1$ for every
$\theta\in I$. Differentiation in each real frequency yields
$\zeta=0$. Thus the action is free.

To prove properness on $V$, suppose
\[
 x_k\to x,\qquad y_k=\Phi_{\zeta_k}(x_k)\to y,
 \qquad x,y,x_k,y_k\in V.
\]
By \eqref{eq:spectral-local-ratios},
\[
 e^{\ii\langle\theta-\theta_0,\zeta_k\rangle}
 =\frac{r_\theta(y_k)}{r_\theta(x_k)}
 \longrightarrow\frac{r_\theta(y)}{r_\theta(x)}.
\]
The convergence is uniform in $\theta\in I$: the ratios are
jointly continuous on a compact neighborhood of each limiting
point times $I$, and their denominators are bounded away from
zero. The limit equals $1$ at $\theta_0$. By
\eqref{eq:spectral-logarithmic-parameter}, $\zeta_k$
converges. The sequential criterion proves properness of the
action on $V$. Such invariant sets cover $P$.

Assume now that $\dim_\C P=p+1$. If the action were not
proper, there would be sequences $x_k\to x$,
$y_k=\Phi_{\zeta_k}(x_k)\to y$ with $|\zeta_k|\to\infty$.
Choose \(j\) and \(\theta_1\) with \(s_{j,\theta_1}(x)\ne0\).
By continuity choose a closed cube \(I\subset\R^p\) centered at
\(\theta_1\), with nonempty interior, on which
\(s_{j,\theta}(x)\ne0\). Fix this \(j\) and write
\(s_\theta=s_{j,\theta}\) for \(\theta\in I\).
Each of these sections is nonzero as a global section of \(L\).

Choose an invariant proper-action neighborhood of the orbit through
\(y\). Lemma~\ref{lem:holomorphic-principal-quotient} gives, after shrinking the one-dimensional quotient
about the image of \(y\), a trivialization
\(\Delta_z\times\C_w^p\) in which \(y=(0,0)\) and
\(\Phi_\zeta(z,w)=(z,w+\zeta)\). Choose a holomorphic frame
\(e_0(z)\) of \(L\) over the zero section. For the covariant
lift \(\widetilde\Phi\) of the action to \(L\), put
\(e(z,w)=\widetilde\Phi_w(e_0(z))\).
By definition the lift acts on sections by
\[
(\Phi_\zeta^*s_\theta)(q)
  =\widetilde\Phi_{-\zeta}\bigl(s_\theta(\Phi_\zeta(q))\bigr),
\qquad
\widetilde\Phi_{-w}\bigl(s_\theta(z,w)\bigr)
  =e^{\ii\langle\theta,w\rangle}s_\theta(z,0),
\]
where the second equality is \eqref{stein:eq:lift} at
\(q=(z,0)\), \(\zeta=w\).
Write \(s_\theta(z,0)=h_\theta(z)e_0(z)\) and apply
\(\widetilde\Phi_w\) to get
\begin{equation}\label{eq:spectral-slice-coefficients}
s_\theta(z,w)=e^{\ii\langle\theta,w\rangle}h_\theta(z)e(z,w),
\qquad h_\theta\in\mathcal O(\Delta).
\end{equation}
None of the \(h_\theta\) is identically zero: otherwise \(s_\theta\)
vanishes on a nonempty open subset of connected \(P\), hence identically,
contrary to its value at \(x\).

Since $h_\theta\not\equiv0$ for every $\theta\in I$,
the vanishing orders $\operatorname{ord}_0h_\theta$
are finite nonnegative integers. Put
\[
m=\min_{\theta\in I}\operatorname{ord}_0h_\theta,
\]
and choose $\theta_*\in I$ with
$\operatorname{ord}_0h_{\theta_*}=m$.

Fix $0<r<1$. By Cauchy's formula, the coefficient
of $z^m$ in the Taylor expansion of $h_\theta$ is
\[
c_m(\theta)
=\frac{1}{2\pi\ii}
\int_{|\xi|=r}
\frac{h_\theta(\xi)}{\xi^{m+1}}\,d\xi.
\]
The compact-open continuity of $h_\theta$ in
$\theta$ implies that $c_m$ is continuous on $I$.
Since $c_m(\theta_*)\ne0$ and
$\operatorname{int}I$ is dense in $I$, we may
choose $\theta_0\in\operatorname{int}I$ such that
\[
c_m(\theta_0)\ne0.
\]

For every $\theta\in I$, the function $h_\theta$
vanishes to order at least $m$ at the origin.
Consequently, we may write
\[
h_\theta(z)=z^m a_\theta(z),
\qquad a_\theta\in\mathcal O(\Delta).
\]
Our choice of $\theta_0$ gives
$a_{\theta_0}(0)\ne0$. After shrinking the
coordinate disc, $a_{\theta_0}$ is nowhere
zero there. Thus, for every $\theta\in I$,
\[
\frac{h_\theta(z)}{h_{\theta_0}(z)}
=\frac{a_\theta(z)}{a_{\theta_0}(z)}
\]
extends holomorphically across $z=0$.

In this common neighborhood the quotients are
\[
r_\theta(z,w)=e^{\ii\langle\theta-\theta_0,w\rangle}
\frac{h_\theta(z)/z^m}{h_{\theta_0}(z)/z^m}.
\]
Globally, the quotients
\[
r_\theta=s_\theta/s_{\theta_0}
\]
are meromorphic functions on \(P\) that extend holomorphically near
\(y\) and are holomorphic near \(x\). For a fixed \(0<r<1\), the divisions satisfy
\[
\frac{h_\theta(z)}{z^m}
 =\frac{1}{2\pi\ii}\int_{|\xi|=r}
    \frac{h_\theta(\xi)}{\xi^m(\xi-z)}\,d\xi ,
 \qquad |z|<r .
\]
The right side depends continuously on \(\theta\), uniformly on
smaller closed discs. Hence the extended ratios depend continuously
on \(\theta\) in both compact-open topologies. Since
\(r_{\theta_0}\equiv1\), choose a
smaller closed cube \(J_0\) centered at \(\theta_0\), contained in
\(\operatorname{int}I\), for which both \(r_\theta(x)\) and
\(r_\theta(y)\) are nonzero for all \(\theta\in J_0\).

From \eqref{stein:eq:lift}, wherever the quotients are defined,
\begin{equation}\label{eq:spectral-ratio-equivariance}
r_\theta\circ\Phi_\zeta
=e^{\ii\langle\theta-\theta_0,\zeta\rangle}r_\theta.
\end{equation}
For large \(k\), \(s_{\theta_0}(x_k)\ne0\), and its
\eqref{stein:eq:lift} gives \(s_{\theta_0}(y_k)\ne0\) as well. Applying
\eqref{eq:spectral-ratio-equivariance} at \(x_k\) and using the
holomorphic extensions at the limiting points gives
\begin{equation}\label{stein:eq:freq-uniform}
e^{\ii\langle\theta-\theta_0,\zeta_k\rangle}
=\frac{r_\theta(y_k)}{r_\theta(x_k)}
\longrightarrow H(\theta):=\frac{r_\theta(y)}{r_\theta(x)}
\end{equation}
uniformly on \(J_0\). Indeed, the extended ratios are jointly continuous on the
product of $J_0$ with compact neighborhoods of $x$ and $y$.
Also $\inf_{\theta\in J_0}|r_\theta(x)|>0$, which bounds the
denominators at $x_k$ away from zero for large $k$. Thus the
convergence in \eqref{stein:eq:freq-uniform} is uniform,
with $H$ continuous and nonzero on $J_0$ and $H(\theta_0)=1$.
The logarithm argument preceding
\eqref{eq:spectral-logarithmic-parameter}, now applied to
\eqref{stein:eq:freq-uniform} on a smaller frequency cube,
shows that $\zeta_k$ converges. This contradicts
$|\zeta_k|\to\infty$.
\end{proof}

For proper real actions on hyperbolic Stein surfaces, an earlier
properness theorem uses tautness or the Bergman metric together
with a cohomological hypothesis
\cite[Theorem~1.2]{MiebachOeljeklaus2009}. The Fourier
construction gives the following result for $p=n-1$, when the
local quotient by $\C^p$ has complex dimension one, without
that cohomological hypothesis.

\begin{theorem}
\label{stein:thm:corank}
Let $\Omega$ be a connected Stein manifold of dimension $n\ge2$
with a proper holomorphic $\R^{n-1}$-action. Assume either
\begin{enumerate}[label=\textup{(\roman*)}]
\item the Bergman kernel form $\mathcal B_\Omega$ is nowhere zero
and the associated Bergman metric is positive definite; or
\item $\Omega$ admits a complete K\"ahler--Einstein metric of
negative Einstein constant, and its Bergman kernel form
$\mathcal B_\Omega$ is nowhere zero.
\end{enumerate}
Then the induced local $\C^{n-1}$-action has a Hausdorff universal
globalization $P$, and its action on $P$ is free and proper.
Consequently
\(P/\C^{n-1}\) is a Hausdorff Riemann surface and
\(P\to P/\C^{n-1}\) is a holomorphic principal
\(\C^{n-1}\)-bundle.
\end{theorem}
\begin{proof}
The universal globalization is Hausdorff by
\cite[Theorem~4.4]{IannuzziSpiroTrapani2004} under \textup{(i)}
and by Proposition~\ref{stein:prop:ambient-moment-slices} under
\textup{(ii)}. In both cases, the sets $T_x$ of
Lemma~\ref{stein:lem:extension} are connected by
\cite[Corollary~3.6]{CTIT2000}.

In either case $\mathcal B_\Omega$ is nowhere zero. For each
$x\in\Omega$, choose $\sigma\in A^2(\Omega,K_\Omega)$ nonzero at
$x$. Proposition~\ref{stein:prop:Fourier} then gives $\theta$
with $S_{\sigma,\theta}(x)\ne0$. As $\sigma$ ranges over $A^2(\Omega,K_\Omega)$ and
$\theta$ ranges over $\R^{n-1}$, the sections
$S_{\sigma,\theta}$ have no common zero on $\Omega$.
Lemma~\ref{stein:lem:extension} extends them to $P$,
where they satisfy \eqref{stein:eq:complex-character}
and still have no common zero.
The action on $K_P^{\otimes2}$ is given by the lift
\[
\widetilde\Phi_\zeta\big|_{(K_P^{\otimes2})_x}
       =\big((\Phi_\zeta)^*_x\big)^{-1}:
       (K_P^{\otimes2})_x\longrightarrow
       (K_P^{\otimes2})_{\Phi_\zeta x}.
\]
The composition law \eqref{eq:canonical-fiber-composition} gives
the group law for this holomorphic lift, and its pullback is
exactly the one in \eqref{stein:eq:complex-character}. Apply
Lemma~\ref{stein:lem:spectral-several} with
\(J=A^2(\Omega,K_\Omega)\), \(p=n-1\) and
\(L=K_P^{\otimes2}\) to obtain freeness and global properness. 
By Lemma~\ref{lem:holomorphic-principal-quotient}, the
free proper action gives
a holomorphic principal $\C^{n-1}$-bundle
\[
\pi:P\longrightarrow S=P/\C^{n-1}.
\]
The quotient $S$ is connected and Hausdorff, and
\[
\dim_\C S=n-(n-1)=1.
\]
Thus $S$ is a connected Riemann surface.
\end{proof}

\begin{corollary}
\label{stein:cor:unconditional-corank}
\label{lem:corank-one-full-proper}
Let $\Omega$ be a connected Stein manifold of dimension $n$,
and let $A=\R^p$, $1\le p\le n$, act properly and
holomorphically on $\Omega$. Suppose that either
\begin{enumerate}[label=\textup{(\roman*)}]
\item $\Omega$ admits a complete K\"ahler--Einstein metric
of negative Einstein constant; or
\item the Bergman kernel form $\mathcal B_\Omega$ is
nowhere zero and the Bergman metric is positive definite.
\end{enumerate}
Then the universal globalization $P$ of the induced local
$\C^p$-action is Hausdorff. Moreover, the following hold.
\begin{enumerate}[label=\textup{(\alph*)}]
\item If $p=n-1\ge1$ and $\mathcal B_\Omega$ is nowhere
zero, the $\C^{n-1}$-action on $P$ is free and proper.
The quotient $S=P/\C^{n-1}$ is a connected Hausdorff
Riemann surface, and
\[
\pi:P\longrightarrow S
\]
is a holomorphic principal $\C^{n-1}$-bundle.

\item If $p=n$, then $P$ is equivariantly biholomorphic
to $\C^n$ with its translation action. In particular,
the $\C^n$-action on $P$ is free and proper, and its
quotient is a point.
\end{enumerate}
\end{corollary}

\begin{proof}
Under \textup{(i)}, the universal globalization is
Hausdorff by Proposition~\ref{stein:prop:ambient-moment-slices}.
Under \textup{(ii)}, Hausdorffness follows from
\cite[Theorem~4.4]{IannuzziSpiroTrapani2004}.

Suppose first that $p=n-1\ge1$ and
$\mathcal B_\Omega$ is nowhere zero.
The hypotheses of Theorem~\ref{stein:thm:corank}
are then satisfied, so the $\C^{n-1}$-action is
free and proper. The assertions concerning the
quotient and the principal bundle follow from
the same theorem.

It remains to prove \textup{(b)}. Suppose that $p=n$.
Under \textup{(i)}, freeness of the complexified action
follows from
Proposition~\ref{stein:prop:ambient-moment-slices}.
Under \textup{(ii)}, we apply the Fourier-section
argument. Since $\mathcal B_\Omega$ is nowhere zero,
for every $x\in\Omega$ there exists
$\sigma\in A^2(\Omega,K_\Omega)$ with $\sigma(x)\ne0$.
By Proposition~\ref{stein:prop:Fourier}, some Fourier
section $S_{\sigma,\theta}$ is nonzero at $x$.
The sets $T_x$ of Lemma~\ref{stein:lem:extension}
are connected by \cite[Corollary~3.6]{CTIT2000}.
That lemma therefore extends the Fourier sections
to $P$, where they have no common zero and satisfy
\eqref{stein:eq:complex-character}.
The first assertion of
Lemma~\ref{stein:lem:spectral-several}
gives freeness of the $\C^n$-action.

Fix $z\in P$ and consider the holomorphic orbit map
\[
f_z:\C^n\longrightarrow P,\qquad
f_z(\zeta)=\Phi_\zeta(z).
\]
Since the action is free, the differential
$(df_z)_0$ is injective. As $\dim_\C P=n$,
it is an isomorphism. The holomorphic inverse
function theorem shows that the orbit through
$z$ is open in $P$. The same holds for every
orbit. Since $P$ is connected, there can be
only one orbit, and hence $f_z$ is surjective.

Freeness also makes $f_z$ injective. Thus $f_z$
is a bijective local biholomorphism, and consequently
a biholomorphism. Finally, the group law gives
\[
\Phi_\eta(f_z(\zeta))=f_z(\zeta+\eta),
\qquad \zeta,\eta\in\C^n.
\]
Hence $f_z$ identifies the action on $P$ with
translations on $\C^n$, proving \textup{(b)}.
\end{proof}

\subsection{Rank-one forms from holomorphic characters}

Once the complexified action is free and proper, its quotient
is a point or a Riemann surface in the cases considered below.
A holomorphic character of $\C^p$ gives a line bundle on
that quotient, from which we construct the required forms.

\begin{theorem}\label{thm:central-corank-one}\label{stein:prop:corank-current-ambient}
Let $\Omega$ be a connected Stein manifold of complex dimension $n$
with a proper holomorphic action of $A=\R^p$,
$p\ge\max\{1,n-1\}$. Assume that $\Omega$ has either a complete
K\"ahler--Einstein metric of negative Einstein constant or
a nowhere-vanishing Bergman kernel form $\mathcal B_\Omega$ and
a positive definite Bergman metric. In the K\"ahler--Einstein
alternative with $p=n-1$, assume additionally that
$\mathcal B_\Omega$ is nowhere zero. For every nonzero $v\in\R^p$ and $K>0$,
there exist smooth positive closed real $(1,1)$-forms of rank one
at every point, belonging to \(\mathcal T_{X_v,K}(\Omega)\) and
\(\mathcal T_{-X_v,K}(\Omega)\), respectively.
\end{theorem}

\begin{proof}

In the K\"ahler--Einstein case,
Proposition~\ref{prop:moment-slices} gives $p\le n$.
In the Bergman case, $P$ is Hausdorff by
\cite[Theorem~4.4]{IannuzziSpiroTrapani2004}.
As in the proof of Theorem~\ref{stein:thm:corank},
Proposition~\ref{stein:prop:Fourier} and
Lemma~\ref{stein:lem:extension} give extended sections
$S_{\sigma,\theta}$ on $P$ with no common zero.
The first assertion of Lemma~\ref{stein:lem:spectral-several}
gives freeness of the complexified action without any
restriction on $p$. For each $z\in P$, the differential
at zero of the orbit map
\[
f_z:\C^p\longrightarrow P,\qquad
\zeta\longmapsto\Phi_\zeta(z),
\]
is therefore injective. Since $\dim_\C P=n$, we obtain
$p\le n$. Together with $p\ge\max\{1,n-1\}$, this gives
$p\in\{n-1,n\}$ for $n\ge2$ and $p=1=n$ for $n=1$.
By Corollary~\ref{lem:corank-one-full-proper},
the quotient map
\[
\pi:P\longrightarrow S=P/\C^p
\]
is a holomorphic principal $\C^p$-bundle.
If $p=n$, then $S$ is a point, whereas if
$p=n-1$, then $S$ is a connected Riemann
surface. Fix $v\in\R^p\setminus\{0\}$
and choose a complex-linear functional
$\alpha:\C^p\to\C$ such that
$\alpha(v)=K/2$.

Define the holomorphic character
$\chi:\C^p\to\C^*$ by
$\chi(\zeta)=e^{\alpha(\zeta)}$.
The associated holomorphic line bundle
$L_\chi\to S$ is the quotient
\[
L_\chi=(P\times\C)/{\sim},
\]
where
\[
(\Phi_\zeta(q),w)\sim
(q,\chi(\zeta)w),
\qquad q\in P,\quad
\zeta\in\C^p,\quad w\in\C.
\]
We write $[q,w]$ for the equivalence
class of $(q,w)$. The bundle projection
is $[q,w]\mapsto\pi(q)$.

Choose an open covering $\{V_i\}$ of $S$
and holomorphic sections
$s_i:V_i\to P$ of $\pi$.
For $x\in V_i\cap V_j$, we have
$\pi(s_i(x))=\pi(s_j(x))=x$.
Since $S=P/\C^p$ and the action is free,
there is a unique $g_{ij}(x)\in\C^p$
such that
\[
s_j(x)=\Phi_{g_{ij}(x)}(s_i(x)).
\]
The holomorphic local trivializations of
the principal bundle show that $g_{ij}$
is holomorphic. On triple intersections,
the group law and freeness give
$g_{ij}+g_{jk}=g_{ik}$.
The local holomorphic frames
$e_i(x)=[s_i(x),1]$ of $L_\chi$
therefore satisfy
\[
e_j=e^{\alpha(g_{ij})}e_i.
\]
Thus the transition functions of $L_\chi$
are $e^{f_{ij}}$, where
$f_{ij}=\alpha(g_{ij})$ is a holomorphic
additive cocycle. Since
\[
f_{ij}+f_{jk}-f_{ik}=0,
\]
the integer-valued cocycle obtained by
taking logarithms in the exponential
sequence is
\[
\frac{f_{ij}+f_{jk}-f_{ik}}{2\pi\ii}=0.
\]
Consequently, $c_1(L_\chi)=0$.

If $S$ is a point, the line is trivial. If $S$ is noncompact,
Behnke--Stein \cite{BehnkeStein1949} and Cartan's Theorem B
\cite[Chapter~IV]{GrauertRemmert2004} give
$H^1(S,\mathcal O)=0$. Thus, after refinement, write
$f_{ij}=b_j-b_i$ with $b_i\in\mathcal O(V_i)$. The local
frames $e^{-b_i}e_i$ agree on overlaps and trivialize $L_\chi$.
In either case the trivialization gives a flat Hermitian metric.

If $S$ is a compact Riemann surface, choose a smooth Hermitian
metric $h_0$ on $L_\chi$. In a holomorphic frame $e$, its real curvature is
$\Theta_0=-\ddbar\log |e|_{h_0}^2$. Write $\Theta_0=f\eta$
for a K\"ahler form $\eta$ on $S$. Since $c_1(L_\chi)=0$,
\[
 \int_S f\,\eta=\int_S\Theta_0=0.
\]
The operator $\Delta_\eta$ of \eqref{eq:trace-laplacian-conventions}
is self-adjoint and elliptic. Its kernel on the connected compact
surface consists of constants, so $\Delta_\eta u=-f$ has a
smooth real solution. Since $\dim_\C S=1$,
$\ddbar u=(\Delta_\eta u)\eta$, and hence
\begin{equation}\label{eq:character-line-flatness}
\Theta(e^{-u}h_0)=\Theta_0+\ddbar u=(f+\Delta_\eta u)\eta=0 .
\end{equation}
Thus in all cases $L_\chi$ has a flat Hermitian metric $h$.

The pullback bundle $\pi^*L_\chi$
has a nowhere-vanishing holomorphic
section $\sigma$, defined by
\[
\sigma(q)=(q,[q,1])\in\pi^*L_\chi.
\]
Define
\[
U(q)=|\sigma(q)|_{\pi^*h}^2
=|[q,1]|_h^2>0.
\]
Since $\pi(\Phi_\zeta(q))=\pi(q)$,
the vectors $[\Phi_\zeta(q),1]$
and $[q,1]$ belong to the same
fiber of $L_\chi$. By the defining
equivalence relation,
\[
[\Phi_\zeta(q),1]
=[q,\chi(\zeta)]
=\chi(\zeta)[q,1].
\]
Taking squared norms gives
\begin{equation}\label{eq:central-corank-U-equivariance}
U(\Phi_\zeta(q))
=e^{2\Re\alpha(\zeta)}U(q),
\qquad q\in P,\quad\zeta\in\C^p.
\end{equation}
In particular, since $\alpha(v)=K/2$,
\begin{equation}\label{eq:central-corank-U-scaling}
U\circ\Phi_{tv}=e^{Kt}U,
\qquad t\in\R.
\end{equation}

Set $\rho=\log U$. Since $\sigma$ is
a nowhere-vanishing holomorphic frame
of $\pi^*L_\chi$, the Chern curvature
of the pullback metric $\pi^*h$ is
$-\ddbar\rho$. As $h$ has zero curvature,
$\ddbar\rho=0$.
Differentiating $U=e^\rho$ and
$U^{-1}=e^{-\rho}$ gives
\begin{equation}\label{eq:central-corank-current-forms}
\begin{aligned}
T_+:=\ddbar U
&=U\,\ii\,\partial\rho
\wedge\bar\partial\rho\ge0,\\
T_-:=\ddbar U^{-1}
&=U^{-1}\,\ii\,\partial\rho
\wedge\bar\partial\rho\ge0.
\end{aligned}
\end{equation}
Thus $T_+$ and $T_-$ are smooth
positive closed real $(1,1)$-forms
on $P$. Differentiating
\eqref{eq:central-corank-U-scaling}
at $t=0$ gives $X_v\rho=K$.
Since $\rho$ is real and $K>0$, this gives $\partial\rho\ne0$
everywhere. By \eqref{eq:central-corank-current-forms},
both $T_+|_\Omega$ and $T_-|_\Omega$ have rank one everywhere. Finally, since $\Phi_{tv}$ is holomorphic,
\eqref{eq:central-corank-U-scaling} gives
\begin{equation}\label{eq:central-corank-current-scaling}
\begin{aligned}
\Phi_{tv}^*T_+
&=\ddbar(U\circ\Phi_{tv})
=e^{Kt}T_+,\\
\Phi_{-tv}^*T_-
&=\ddbar(U^{-1}\circ\Phi_{-tv})
=e^{Kt}T_-,
\qquad t\in\R.
\end{aligned}
\end{equation}
Restricting \eqref{eq:central-corank-current-scaling}
to $\Omega$ proves
\eqref{eq:current-homogeneity} for
the flows generated by $X_v$ and
$-X_v$, respectively.
\end{proof}

\subsection{Averaging over compact subgroups}

In $G\simeq\R^p\times C$, the flow of $X_v+W$ is the product of
the commuting flows of $X_v$ and $W\in\operatorname{Lie}(C)$.
We average over $C$ to keep the character of the first flow.

\begin{proposition}\label{prop:compact-component}
Let $\Phi_t$ commute with a compact group $H$ of automorphisms,
and let $T$ be a nonzero positive closed $(1,1)$-current with
$\Phi_t^*T=e^{Kt}T$. For every one-parameter subgroup $c_t$ of $H$,
the flow $c_t\circ\Phi_t$ has a nonzero positive closed current
of character $e^{Kt}$. This current is smooth when $T$ is smooth.
\end{proposition}

\begin{proof}
For normalized Haar measure \(dh\), define
\[
\overline T(\beta)=\int_H T((h^{-1})^*\beta)\,dh .
\]
For a test form supported in a compact set \(Q\), all transported
forms are supported in the compact set \(H\cdot Q\), with uniform
bounds on their derivatives. Thus the integral defines a current
and commutes with \(d\). It is positive, and if \(T(\beta)>0\)
for a strongly positive test form \(\beta\), its integrand remains
positive on a neighborhood of the identity. Hence
\(\overline T\ne0\). Right invariance of Haar measure and
commutation with \(\Phi_t\) give
\[
c^*\overline T=\int_H (hc)^*T\,dh=\overline T,\qquad
\Phi_t^*\overline T=\int_H h^*(\Phi_t^*T)\,dh=e^{Kt}\overline T .
\]
These give \((c_t\circ\Phi_t)^*\overline T=e^{Kt}\overline T\).
If \(T\) is smooth, the same compact bounds allow derivatives
under the integral and \(\overline T\) is smooth.
\end{proof}

In the decomposition $G\simeq\R^p\times C$, a flow is
proper precisely when its generator has nonzero $\R^p$-component,
by Lemma~\ref{stein:lem:realflow}. If both cones for $X_v$ and
$-X_v$ contain nonzero currents, Proposition~\ref{prop:compact-component}
gives the same assertion for $X_v+W$ and $-X_v-W$ for any
$W\in\operatorname{Lie}(C)$.

\begin{proof}[Proof of Theorem~\ref{stein:thm:ambient-unified}]
Under the identification of $\mathfrak g_\C$ with holomorphic
vector fields established in Section~\ref{sec:central-reduction},
the operator $\ad_Y$ acts by $W\mapsto[Y^{1,0},W]$.
Proposition~\ref{prop:adjoint-branches} therefore gives both
the eigenvalue and Jordan-block conclusions in \textup{(i)}.
For \textup{(ii)}, Proposition~\ref{prop:compact-type-structure}
gives \eqref{eq:central-decomposition} and pointwise independence.
When $\Omega$ is Stein, $\mathcal B_\Omega$ is nowhere zero,
the Bergman metric is positive definite, and $p\ge n-1$, consider a field $Y=X_v+W$ with $v\ne0$ and
$W\in\operatorname{Lie}(C)$. Apply
Theorem~\ref{thm:central-corank-one} to $X_v$ and $-X_v$.
Averaging the two forms by Proposition~\ref{prop:compact-component}
gives the required forms for $Y$ and $-Y$.

For the final assertion, suppose additionally that
$\Omega\Subset Z$ is Stein and that a neighborhood of
$\overline\Omega$ admits a holomorphic immersion into
$\C^N$. Whenever \textup{(i)} or \textup{(ii)}
provides a nonzero current in $\mathcal T_{Y,K}(\Omega)$,
\cite[Theorem~1.2]{GangulyKRsol2026Paper21}
gives a complete nontrivial expander with soliton
field $Y$ and expansion constant $K$.
\end{proof}

\begin{proof}[Proof of Corollary~\ref{cor:unified-existence}]
By Lemma~\ref{stein:lem:bounded}, bounded holomorphic
functions on $\Omega$ separate tangent directions,
and $G=\Aut_0(\Omega)$ acts properly.
Since $G$ is noncompact,
Lemma~\ref{stein:lem:realflow} supplies a proper
real holomorphic flow.

Fix $K>0$. If the adjoint of some proper generator
has an eigenvalue with nonzero real part or a
nontrivial Jordan block corresponding to a purely
imaginary eigenvalue, Theorem~\ref{thm:unified-current}
\textup{(i)} gives a nonzero homogeneous current
for at least one orientation of that generator.

Otherwise, for every proper generator, the adjoint
operator is diagonalizable over $\C$ with purely
imaginary spectrum. By Theorem~\ref{thm:unified-current}
\textup{(ii)},
\[
G\simeq\R^p\times C,\qquad 1\le p\le n,
\]
where $C$ is compact and connected.
The hypothesis excluding $1\le p\le n-2$
forces $p\ge n-1$.

Lemma~\ref{stein:lem:bounded} shows that
$\mathcal B_\Omega$ is nowhere zero and that
the Bergman metric is positive definite.
Theorem~\ref{thm:unified-current}\textup{(ii)}
therefore gives a nonzero homogeneous current
for a nonzero generator of the central real factor.

In either case,
\cite[Theorem~1.2]{GangulyKRsol2026Paper21}
produces a complete nontrivial expander with
complete soliton field and expansion constant $K$.
\end{proof}

\begin{proof}[Proof of Corollary~\ref{thm:stein-surface-intro}]
By Lemma~\ref{stein:lem:bounded}, $\Aut(\Omega)$ acts properly.
Suppose that $\Omega$ is Stein and that the flow of $Y$ has
noncompact closure. By Lemma~\ref{stein:lem:realflow},
this flow defines a proper $\R$-action.
The immersion hypothesis also gives a complete
K\"ahler--Einstein metric $\omega_E$, normalized by
$\Ric_{\omega_E}=-\omega_E$
\cite[Remarks~3.6(2)]{GangulyKRsol2026Paper21}.
Moreover, the proof of Corollary~\ref{stein:cor:finite}
provides a finite family of square-integrable holomorphic
$n$-forms with no common zero on $\Omega$.
Since $p=1=n-1$, Theorem~\ref{thm:central-corank-one}
gives rank-one forms in $\mathcal T_{Y,K}(\Omega)$ and
$\mathcal T_{-Y,K}(\Omega)$ for every $K>0$.
\cite[Theorem~1.2]{GangulyKRsol2026Paper21} then gives
the two solitons. If $\Aut_0(\Omega)$ is noncompact,
Lemma~\ref{stein:lem:realflow} supplies a proper real flow.

Conversely, a complete nontrivial expander makes $\Omega$
Stein by \cite[Remarks~3.6(2)]{GangulyKRsol2026Paper21}. Its
flow has noncompact closure by
\cite[Lemma~3.1]{GangulyKRsol2026Paper21}, and hence the
identity component is noncompact.
\end{proof}

\section{Boundary Limits of Automorphisms}
\label{sec:boundary-dynamics}

We study locally uniform limits in $Z$ of holomorphic
self-maps of $\Omega\Subset Z$. The criteria for constant
boundary limits established below will be used in
Section~\ref{sec:geometric-applications}.

\subsection{Limits of holomorphic self-maps}
For a holomorphic flow $\Phi$ on $\Omega\Subset Z$, let
$\mathscr N_\pm(\Phi)$ consist of the maps $g:\Omega\to Z$
for which some sequence $t_j\to\pm\infty$ satisfies
$\Phi_{t_j}\to g$ locally uniformly in $Z$.

\begin{lemma}
\label{lem:ambient-normal}
Let $\Omega\Subset Z$ be a domain, and let
$d$ be the distance induced by a smooth Hermitian metric on $Z$.
Suppose that for every $b>0$ there is $\delta_b>0$ such that
\begin{equation}\label{eq:stein-bounded-separation}
 k_\Omega(z,w)\ge\delta_b
 \qquad\bigl(z,w\in\Omega,\ d(z,w)\ge b\bigr).
\end{equation}
Every sequence of holomorphic self-maps of $\Omega$ has a
subsequence converging locally uniformly in $Z$ to a holomorphic
map $F_\infty:\Omega\to Z$ with
$F_\infty(\Omega)\subset\overline\Omega$.
Condition~\eqref{eq:stein-bounded-separation} holds if a neighborhood
of $\overline\Omega$ admits a holomorphic immersion $\iota$ into
$\C^N$, or if the Kobayashi hyperbolic domain \(\Omega\) satisfies
relative \((1,\kappa_0)\)-visibility for one \(\kappa_0>0\).
\end{lemma}

\begin{proof}
Let \(B_{2r}\Subset\Omega\) be a coordinate ball, with
coordinate map \(z\), and let \(x,y\in B_r\). Comparison with
the ball metric and contraction under \(F:\Omega\to\Omega\) give
\begin{equation}\label{eq:normal-chart-contraction}
k_\Omega(Fx,Fy)\le k_\Omega(x,y)
 \le k_{B_{2r}}(x,y)\le C_r|z(x)-z(y)|.
\end{equation}
For each $b>0$, \eqref{eq:stein-bounded-separation} gives
$d(Fx,Fy)<b$ whenever $C_r|z(x)-z(y)|<\delta_b$. This
estimate, applied on finitely many coordinate balls covering
each compact subset of $\Omega$, gives equicontinuity there.

Let $(F_j)$ be a sequence of holomorphic self-maps of
$\Omega$. Since their images lie in the compact set
$\overline\Omega$, the equicontinuity established above
and Arzel\`a--Ascoli on a compact exhaustion give a
subsequence $(F_{j_\ell})$ converging locally uniformly
in $Z$ to a continuous map
$F_\infty:\Omega\to\overline\Omega$.

Fix $a\in\Omega$ and choose a holomorphic coordinate
chart $(W,w)$ of $Z$ containing $F_\infty(a)$.
By continuity, there is an open neighborhood
$U\Subset\Omega$ of $a$ such that
$F_\infty(\overline U)$ is a compact subset of $W$.
Uniform convergence on $\overline U$ then gives
$F_{j_\ell}(\overline U)\subset W$ for all sufficiently
large $\ell$. Consequently, $w\circ F_\infty$ is the
locally uniform limit on $U$ of the holomorphic maps
$w\circ F_{j_\ell}$, and hence is holomorphic.
Since $a$ was arbitrary, $F_\infty$ is holomorphic
on $\Omega$.

Suppose first that an immersion \(\iota\) exists on a
neighborhood of \(\overline\Omega\). Choose a Euclidean ball
\(B\) containing \(\iota(\overline\Omega)\). There are
\(c_1,c_2>0\) with \(\kappa_B(\zeta;\xi)\ge c_1|\xi|\) on
this compact image and \(|d\iota_x u|\ge c_2|u|_Z\) on
\(\overline\Omega\). For every holomorphic disc \(f:\D\to\Omega\),
\[
\kappa_{\D}(\zeta;v)
 \ge \kappa_B(\iota(f(\zeta));d(\iota\circ f)_\zeta v)
 \ge c_1c_2 |df_\zeta v|_Z .
\]
Integrating on a Poincar\'e geodesic from \(\alpha\) to \(\beta\)
gives \(k_{\D}(\alpha,\beta)\ge c_1c_2d(f(\alpha),f(\beta))\).
Summing over a chain of discs between \(x,y\) and taking the
infimum gives
\begin{equation}\label{eq:immersion-distance-floor}
k_\Omega(x,y)\ge c_1c_2d(x,y).
\end{equation}
In particular \(\Omega\) is hyperbolic and satisfies
\eqref{eq:stein-bounded-separation}. For a bounded Euclidean
domain one can take the identity immersion and \(d\) Euclidean.

Now assume relative \((1,\kappa_0)\)-visibility on the
Kobayashi hyperbolic \(\Omega\). If
\eqref{eq:stein-bounded-separation} fails, there are
\(z_j\to z\), \(w_j\to w\) with \(z\ne w\) in the compact
closure and \(k_\Omega(z_j,w_j)\to0\). If \(z\in\Omega\),
then \(k_\Omega(z,w_j)\le k_\Omega(z,z_j)+k_\Omega(z_j,w_j)\to0\).
Since $\Omega$ is Kobayashi hyperbolic, the topology
induced by $k_\Omega$ coincides with the original
topology of $\Omega$ by
\cite{Barth1972Kobayashi}. So $w_j\to z$, contradicting
$w\ne z$. The same argument rules out $w\in\Omega$.
Consequently, both $z$ and $w$ lie in
$\partial_Z\Omega$.

For a \((1,\kappa)\)-almost-geodesic \(\gamma:[0,T]\to\Omega\),
set \(x=\gamma(0)\), \(y=\gamma(T)\) and \(y_0=\gamma(u)\).
From \eqref{eq:almost-geodesic-convention},
\[
k_\Omega(x,y_0)\le u+\kappa,\quad
k_\Omega(y_0,y)\le T-u+\kappa,\quad
T\le k_\Omega(x,y)+\kappa.
\]
Adding gives
\begin{equation}\label{eq:almost-geodesic-splitting}
k_\Omega(x,y_0)+k_\Omega(y_0,y)\le k_\Omega(x,y)+3\kappa.
\end{equation}
Choose $(1,\kappa_j)$-almost-geodesics joining $z_j$
to $w_j$, where $0<\kappa_j\le\kappa_0$ and
$\kappa_j\to0$
\cite[Proposition~2.8]{Masanta2025Visibility}.
Each of these curves is also a
$(1,\kappa_0)$-almost-geodesic. Since $z_j\to z$
and $w_j\to w$, with $z,w$ distinct points of
$\partial_Z\Omega$, relative visibility provides
a compact set $Q\Subset\Omega$ meeting every
chosen curve for all sufficiently large $j$.
Choose an intersection point $y_j\in Q$ on
each such curve. Passing to a subsequence,
we may assume that $y_j\to y_*\in Q$. Then
\eqref{eq:almost-geodesic-splitting} gives
\[
k_\Omega(z_j,y_j)+k_\Omega(y_j,w_j)
 \le k_\Omega(z_j,w_j)+3\kappa_j\longrightarrow0 .
\]
Since $y_j\to y_*\in\Omega$, continuity of the
Kobayashi distance gives $k_\Omega(y_j,y_*)\to0$.
The triangle inequality therefore yields
\[
\begin{aligned}
k_\Omega(z_j,y_*)
&\le k_\Omega(z_j,y_j)+k_\Omega(y_j,y_*)
\longrightarrow0,\\
k_\Omega(w_j,y_*)
&\le k_\Omega(w_j,y_j)+k_\Omega(y_j,y_*)
\longrightarrow0.
\end{aligned}
\]
By \cite{Barth1972Kobayashi}, convergence with
respect to $k_\Omega$ implies convergence in
the original topology of $\Omega$. Thus
$z_j\to y_*$ and $w_j\to y_*$ in $\Omega$.
Since the inclusion $\Omega\hookrightarrow Z$
is continuous, both sequences converge to
$y_*$ in $Z$, contradicting
$z_j\to z$ and $w_j\to w$ in $Z$, where
$z,w\in\partial_Z\Omega$.
\end{proof}

\begin{lemma}
\label{lem:stein-normal-limit}\label{lem:backward-boundary-limit}
\label{stein:lem:intrinsic-jac-decay}
Let $\Omega\Subset Z$ be a domain such that every sequence
of holomorphic self-maps of $\Omega$ admits a subsequence
converging locally uniformly in $Z$ to a holomorphic map
into $\overline\Omega$, as in Lemma~\ref{lem:ambient-normal}.
Suppose that $F_j\in\Aut(\Omega)$ satisfies
\[
F_j(E)\cap L=\varnothing
\]
for all sufficiently large $j$, for every pair of compact
sets $E,L\Subset\Omega$.

If a subsequence of $F_j$ converges locally uniformly in
$Z$ to $g:\Omega\to\overline\Omega$, then
\[
g(\Omega)\subset\partial_Z\Omega,
\qquad
\operatorname{rank}_{\C}dg<n
\quad\text{on }\Omega.
\]
Moreover, $g$ is constant if either
\begin{enumerate}[label=\textup{(\roman*)}]
\item $\partial_Z\Omega$ contains no nonconstant analytic disc;
\item $\Omega$ is Kobayashi hyperbolic and satisfies relative
$(1,\kappa_0)$-visibility for some $\kappa_0>0$.
\end{enumerate}
\end{lemma}

\begin{proof}
Passing to a subsequence, we may assume that $F_j\to g$
locally uniformly in $Z$. Suppose that $g(a)\in\Omega$
for some $a\in\Omega$. Choose a compact neighborhood
$L\Subset\Omega$ of $g(a)$. Since $F_j(a)\to g(a)$,
we have $F_j(a)\in L$ for all sufficiently large $j$.
This contradicts the hypothesis with $E=\{a\}$.
Therefore
\[
g(\Omega)\subset\partial_Z\Omega.
\]
If $dg_a$ had rank $n$ at some $a\in\Omega$, the
holomorphic inverse function theorem would imply that
$g(\Omega)$ contains a nonempty open subset of $Z$.
This is impossible because $\partial_Z\Omega$ has
empty interior. Hence $\operatorname{rank}_{\C}dg<n$
throughout $\Omega$.

Suppose that the boundary contains no nonconstant
analytic disc. If $g$ were nonconstant, there would
exist $a\in\Omega$ and $v\in T_a^{1,0}\Omega$ such that
$dg_a(v)\ne0$. Choose holomorphic coordinates $z$
near $a$. For sufficiently small $\varepsilon>0$, the map
\[
\zeta\longmapsto
g\bigl(z^{-1}(z(a)+\varepsilon\zeta\,dz_a(v))\bigr),
\qquad \zeta\in\D,
\]
would be a nonconstant holomorphic disc contained in
$\partial_Z\Omega$, a contradiction. Thus $g$ is constant.

Now suppose that $\Omega$ satisfies condition~\textup{(ii)}.
If $g$ were nonconstant, we could choose $a,b\in\Omega$
with $g(a)\ne g(b)$. By
\cite[Proposition~2.8]{Masanta2025Visibility},
there exists a $(1,\kappa_0)$-almost-geodesic
$\gamma:[0,T]\to\Omega$ joining $a$ to $b$.

Set $\gamma_j=F_j\circ\gamma$. Since $F_j$ is an
automorphism, invariance of the Kobayashi distance
and infinitesimal metric gives
\[
\begin{aligned}
k_\Omega(\gamma_j(s),\gamma_j(t))
 &=k_\Omega(\gamma(s),\gamma(t)),\\
\kappa_\Omega(\gamma_j(t);\dot\gamma_j(t))
 &=\kappa_\Omega(\gamma(t);\dot\gamma(t))
 \le1
\end{aligned}
\]
for $s,t\in[0,T]$, with the second identity holding
for almost every $t$. Thus each $\gamma_j$ is a
$(1,\kappa_0)$-almost-geodesic.

Since
\[
\gamma_j(0)=F_j(a)\longrightarrow g(a),
\qquad
\gamma_j(T)=F_j(b)\longrightarrow g(b),
\]
and $g(a),g(b)$ are distinct boundary points, relative
$(1,\kappa_0)$-visibility gives a compact set
$Q\Subset\Omega$ such that
\[
\gamma_j([0,T])\cap Q\ne\varnothing
\]
for all sufficiently large $j$. Equivalently,
\[
F_j(\gamma([0,T]))\cap Q\ne\varnothing,
\]
contrary to the hypothesis with
$E=\gamma([0,T])$ and $L=Q$.
Therefore $g$ is constant.
\end{proof}

\begin{remark}\label{rem:normal-limit-psh-barrier}
Let $\Omega\Subset Z$ be a domain such that every
sequence of holomorphic self-maps of $\Omega$
admits a subsequence converging locally uniformly
in $Z$ to a holomorphic map
$\Omega\to\overline\Omega$, and let $\Phi_t$
be a proper real holomorphic flow on $\Omega$,
with generator $Y$.

Suppose that a neighborhood $U$ of $\overline\Omega$
admits a bounded nonnegative function
$\beta\in\PSH(U)$ such that
$\beta|_\Omega\not\equiv0$ and
\[
\beta|_{g(\Omega)}=0
\qquad\text{for every }g\in\mathscr N_-(\Phi).
\]
Then
\[
\beta\circ\Phi_t\longrightarrow0
\quad\text{locally uniformly as }t\to-\infty,
\]
and $\mathcal T_{Y,K}(\Omega)$ contains a nonzero
current for every $K>0$.

Indeed, if the convergence failed, there would exist
$\epsilon>0$, a compact set $L\Subset\Omega$, and
sequences $t_j\to-\infty$ and $z_j\in L$ such that
\[
\beta(\Phi_{t_j}(z_j))\ge\epsilon.
\]
After passing to a subsequence, we have
$z_j\to z\in L$ and $\Phi_{t_j}\to g$ locally
uniformly in $Z$, for some $g\in\mathscr N_-(\Phi)$.
In particular, $\Phi_{t_j}(z_j)\to g(z)$.
Upper semicontinuity of $\beta$ gives
\[
\epsilon
\le\limsup_{j\to\infty}\beta(\Phi_{t_j}(z_j))
\le\beta(g(z))=0,
\]
a contradiction.

Choose $q\in\Omega$ with $\beta(q)>0$ and
$0<\delta<\beta(q)$. Then
\[
b=(\beta|_\Omega-\delta)_+
\]
is nonzero and satisfies
\eqref{eq:backward-cutoff-vanishing}.
Lemma~\ref{lem:backward-cutoff-current} therefore
gives a nonzero current in
$\mathcal T_{Y,K}(\Omega)$ for every $K>0$.
\end{remark}

\subsection{Iterates on weakly visible domains}

\begin{lemma}\label{lem:weak-visibility-iterates}\label{lem:stein-visibility-iterates}
Let $\Omega\Subset Z$ be a Kobayashi hyperbolic domain
satisfying relative $(1,\kappa_0)$-visibility for some
$\kappa_0>0$, and let $F\in\Aut(\Omega)$.
Suppose that for every pair of compact sets $E,L\Subset\Omega$,
\begin{equation}\label{eq:weakvis-compact-divergence}
F^m(E)\cap L=\varnothing
\qquad\text{for all sufficiently large }m.
\end{equation}
Then there exists $q\in\partial_Z\Omega$ such that $F^m\to q$ locally uniformly on $\Omega$ as $m\to\infty$.
\end{lemma}

\begin{proof}
By Lemmas~\ref{lem:ambient-normal} and
\ref{lem:stein-normal-limit}, every sequence of positive
iterates admits a subsequence converging locally uniformly in
$Z$ to a constant map into $\partial_Z\Omega$.
In particular, suppose that $r_j\to\infty$ and
$F^{r_j}(a)\to\xi$ for some $a\in\Omega$. Every subsequence
of $(F^{r_j})$ has a locally uniformly convergent further
subsequence, whose limit is constant and must equal $\xi$
by evaluation at $a$. Hence $F^{r_j}\to\xi$ locally
uniformly. It remains to show that different sequences of
positive iterates cannot approach distinct boundary points.
For a related iteration result on taut weakly visible domains,
see \cite[Theorem~6.2]{Masanta2025Visibility}.

Fix a smooth Hermitian metric $h$ on $Z$,
and let $d=d_h$ be the induced distance.
By Lemma~\ref{lem:ambient-normal},
condition~\eqref{eq:stein-bounded-separation}
holds for this distance. Fix $o\in\Omega$
and set
\[
d_m=k_\Omega(o,F^m(o)),\qquad m\ge0.
\]
We first establish an estimate for points approaching distinct
boundary points. For distinct $\xi,\eta\in\partial_Z\Omega$, visibility gives
neighborhoods $U_\xi,U_\eta$ and a compact set $Q$ meeting the
$(1,\kappa_0)$-almost-geodesics between these neighborhoods. For
any pair $x\in U_\xi\cap\Omega$, $y\in U_\eta\cap\Omega$,
\cite[Proposition~2.8]{Masanta2025Visibility} supplies such a
curve. Take $y_0\in Q$ on it. By \eqref{eq:almost-geodesic-splitting}
and the triangle inequality at $o$, with
$C=2\max_Q k_\Omega(o,\cdot)+3\kappa_0$, we obtain
\begin{equation}\label{eq:weakvis-separated-boundary}
k_\Omega(x,y)\ge k_\Omega(o,x)+k_\Omega(o,y)-C,
\qquad x\in U_\xi\cap\Omega,\quad y\in U_\eta\cap\Omega.
\end{equation}

Suppose first that $(d_m)$ is unbounded.
Choose increasing integers $n_i$ such that
\[
d_{n_i}=\max_{0\le r\le n_i}d_r\to\infty.
\]
In particular, $d_r\le d_{n_i}$ whenever
$0\le r\le n_i$. After passing to a subsequence,
$F^{n_i}\to\xi\in\partial_Z\Omega$ locally uniformly.
Suppose that there exist $r_j\to\infty$ with
$d_{r_j}\to\infty$ and $F^{r_j}\to\eta\ne\xi$
locally uniformly. For all sufficiently large
$j$, we have $F^{r_j}(o)\in U_\eta$. Choose
$i=i(j)$ sufficiently large that $n_i>r_j$
and $F^{n_i}(o)\in U_\xi$. By invariance of
the Kobayashi distance and the choice of $n_i$,
\[
k_\Omega(F^{n_i}(o),F^{r_j}(o))
=d_{n_i-r_j}\le d_{n_i},
\]
whereas \eqref{eq:weakvis-separated-boundary}
gives
\[
k_\Omega(F^{n_i}(o),F^{r_j}(o))
\ge d_{n_i}+d_{r_j}-C.
\]
Thus $d_{r_j}\le C$, contradicting
$d_{r_j}\to\infty$. By Lemma~\ref{lem:ambient-normal},
for every sequence $r_j\to\infty$ satisfying
$d_{r_j}\to\infty$, the maps $F^{r_j}$
therefore converge locally uniformly to $\xi$.

Unboundedness of $(d_m)$ does not imply that $d_m\to\infty$.
We must therefore rule out a different boundary limit along
iterates whose distances from $o$ remain bounded. Suppose
that $F^{m_j}\to\eta\ne\xi$ locally uniformly.
Choose a neighborhood $U$ of $\eta$ whose
closure does not contain $\xi$. Then
\[
M:=\sup\{d_m:m\ge0,\ F^m(o)\in U\}<\infty.
\]
Indeed, otherwise there would exist
$r_j\to\infty$ such that $F^{r_j}(o)\in U$
and $d_{r_j}\to\infty$. The preceding
conclusion would imply $F^{r_j}(o)\to\xi$,
contradicting $\xi\notin\overline U$.
For each fixed $i$, we have
\[
F^{m_j-n_i}(o)
=F^{m_j}(F^{-n_i}(o))\to\eta
\qquad(j\to\infty),
\]
where $m_j\ge n_i$ for all sufficiently
large $j$. Consequently, both
$F^{m_j-n_i}(o)$ and $F^{m_j}(o)$ belong
to $U$ for sufficiently large $j$.
By the triangle inequality and invariance
of the Kobayashi distance,
\[
\begin{aligned}
d_{n_i}
&=k_\Omega(F^{m_j}(o),F^{m_j-n_i}(o)) \le d_{m_j}+d_{m_j-n_i} \le 2M.
\end{aligned}
\]
Since $i$ was arbitrary, this contradicts
$d_{n_i}\to\infty$.

We now consider the case
$R:=\sup_{m\ge0}d_m<\infty$.
We adapt the argument in
\cite[proof of Theorem~1.8, Case~2]{BharaliMaitra2021}
to the present setting. Suppose, for a contradiction, that
there are sequences $n_i,m_j\to\infty$ for which
$F^{n_i}\to\xi$ and $F^{m_j}\to\eta$ locally uniformly on
$\Omega$, with distinct points
$\xi,\eta\in\partial_Z\Omega$.

Relative visibility supplies neighborhoods of $\xi$ and $\eta$
and a compact set $Q\Subset\Omega$ meeting every
$(1,\kappa_0)$-almost-geodesic joining points in these
neighborhoods. For $a,b\in Q$ and $\delta>0$, define
\begin{equation}\label{eq:weakvis-H-definition}
\begin{aligned}
H_\delta(a,b)
&=\inf_{\substack{r\in\Z_{\ge1}\\
d(F^r(a),\xi)<\delta}}
k_\Omega(F^r(a),b),\\
H(a,b)
&=\lim_{\delta\downarrow0}H_\delta(a,b).
\end{aligned}
\end{equation}
Since $F^{n_i}(a)\to\xi$, for every
$\delta>0$ there exists an integer $r\ge1$
such that $d(F^r(a),\xi)<\delta$.
Moreover, by the triangle inequality and
invariance of the Kobayashi distance,
\[
k_\Omega(F^r(a),b)
\le k_\Omega(a,o)+d_r+k_\Omega(o,b)
\le R+2\max_{z\in Q}k_\Omega(o,z).
\]
Thus the set of integers over which the infimum defining
$H_\delta(a,b)$ is taken is nonempty, and the preceding
estimate bounds $H_\delta(a,b)$ uniformly for $a,b\in Q$
and $\delta>0$. As $\delta\downarrow0$, this infimum is
nondecreasing. Consequently, $H(a,b)$ exists and is finite.

Set $\rho=\inf_{y\in Q}d(\eta,y)>0$.
Since $F^{m_j}\to\eta$ uniformly on $Q$,
for all sufficiently large $j$ we have
\[
d(F^{m_j}(b),\eta)<\rho/2
\qquad(b\in Q).
\]
Consequently, for every $y,b\in Q$,
\[
d(y,F^{m_j}(b))
\ge d(y,\eta)-d(F^{m_j}(b),\eta)
\ge\rho/2.
\]
Thus, with $\varepsilon:=\delta_{\rho/2}>0$ from
\eqref{eq:stein-bounded-separation},
\begin{equation}\label{eq:weakvis-compact-separation}
k_\Omega(y,F^{m_j}(b))\ge\varepsilon
\qquad(y,b\in Q)
\end{equation}
for all sufficiently large $j$. In the bounded Euclidean case,
\eqref{eq:immersion-distance-floor} gives the choice
\(\varepsilon=c_1c_2\rho/2\).
Choose $a,b\in Q$ satisfying
\begin{equation}\label{eq:weakvis-near-minimum}
H(a,b)<\inf_{x,y\in Q}H(x,y)+\varepsilon/2.
\end{equation}

For each integer $s\ge1$, the definition of
$H_{1/s}(a,b)$ involves the nonempty set of
positive integers $r$ satisfying
\[
d(F^r(a),\xi)<1/s.
\]
We may therefore choose $p_s\in\Z_{\ge1}$ such that
\[
\begin{aligned}
d(F^{p_s}(a),\xi)&<1/s,\\
H_{1/s}(a,b)
&\le k_\Omega(F^{p_s}(a),b)
 <H_{1/s}(a,b)+1/s
 \le H(a,b)+1/s.
\end{aligned}
\]
Since $H_{1/s}(a,b)\to H(a,b)$ as $s\to\infty$,
it follows that
\[
k_\Omega(F^{p_s}(a),b)\longrightarrow H(a,b).
\]
Moreover, $F^{p_s}(a)\to\xi\in\partial_Z\Omega$,
so $p_s\to\infty$. By the first paragraph of the
proof, $F^{p_s}\to\xi$ locally uniformly on
$\Omega$. Hence, for each fixed $j$,
\[
F^{p_s+m_j}(a)
=F^{p_s}(F^{m_j}(a))\to\xi
\qquad(s\to\infty).
\]
Choose $s(j)\ge j$ sufficiently large that
\[
d(F^{p_{s(j)}+m_j}(a),\xi)<1/j,
\]
and put $r_j=p_{s(j)}$.
Then $F^{r_j+m_j}(a)\to\xi$.
Since $s(j)\to\infty$, we also have
\begin{equation}\label{eq:weakvis-minimizing-iterates}
k_\Omega(F^{r_j}(a),b)\to H(a,b).
\end{equation}

Set $A_j=F^{r_j+m_j}(a)$ and $B_j=F^{m_j}(b)$.
Since $A_j\to\xi$ and $B_j\to\eta$, the endpoints lie in
the neighborhoods used to choose $Q$ for all sufficiently
large $j$. Choose $(1,\kappa_j)$-almost-geodesics joining
$A_j$ to $B_j$, where $\kappa_j\downarrow0$ and
$\kappa_j\le\kappa_0$, by
\cite[Proposition~2.8]{Masanta2025Visibility}.
Each is also a $(1,\kappa_0)$-almost-geodesic and hence
meets $Q$. Choose intersection points $y_j\in Q$ and pass
to a subsequence such that $y_j\to y\in Q$.
By \eqref{eq:almost-geodesic-splitting},
\begin{equation}\label{eq:weakvis-minimizing-splitting}
k_\Omega(A_j,B_j)\ge k_\Omega(A_j,y_j)+k_\Omega(y_j,B_j)-3\kappa_j.
\end{equation}

Since $A_j=F^{r_j+m_j}(a)\to\xi$,
equation~\eqref{eq:weakvis-H-definition}
gives, for every fixed $\delta>0$,
\[
k_\Omega(A_j,y)\ge H_\delta(a,y)
\]
for all sufficiently large $j$.

As $y_j\to y\in\Omega$, continuity of the
Kobayashi distance
\cite{Barth1972Kobayashi}
gives $k_\Omega(y_j,y)\to0$.
The triangle inequality then yields
\[
k_\Omega(A_j,y_j)
\ge k_\Omega(A_j,y)-k_\Omega(y_j,y).
\]
Consequently,
\[
\liminf_{j\to\infty}k_\Omega(A_j,y_j)
\ge H_\delta(a,y).
\]
Letting $\delta\downarrow0$, we obtain
\[
\liminf_{j\to\infty}k_\Omega(A_j,y_j)
\ge H(a,y).
\]
Since $k_\Omega(A_j,B_j)=k_\Omega(F^{r_j}(a),b)$, equations
\eqref{eq:weakvis-minimizing-iterates} and
\eqref{eq:weakvis-minimizing-splitting}, together with
\eqref{eq:weakvis-compact-separation}, yield
\[
H(a,b)\ge H(a,y)+\varepsilon
\ge\inf_{x,z\in Q}H(x,z)+\varepsilon,
\]
contradicting \eqref{eq:weakvis-near-minimum}.
Therefore all locally uniform subsequential limits of
$(F^m)$ in $Z$ equal the same constant map with value
$q\in\partial_Z\Omega$. By Lemma~\ref{lem:ambient-normal},
every sequence of iterates has a locally uniformly convergent
subsequence. The uniqueness of its limit now implies that
the entire sequence $F^m$ converges locally uniformly to
$q$ as $m\to\infty$.
\end{proof}

\section{Geometric Criteria for Expanding Solitons}
\label{sec:geometric-applications}

\subsection{Boundary limits and displacement functions}
\label{sec:stein-geometric}\label{pcvx:subsec:finite-type}

We first construct homogeneous currents from the displacement
$f\circ\Phi_\tau-f$, without assuming an ambient
holomorphic immersion.

\begin{proposition}\label{prop:geometric-beyond-immersion}
Let $\Omega\Subset Z$ satisfy either of the following:
\begin{enumerate}[label=\textup{(\roman*)}]
\item The relative boundary $\partial_Z\Omega$
contains no nonconstant analytic disc, and
\eqref{eq:stein-bounded-separation} holds.
\item The domain $\Omega$ is Kobayashi hyperbolic
and has relative $(1,\kappa_0)$-visibility
for some $\kappa_0>0$.
\end{enumerate}
Let $Y$ generate a proper real holomorphic
flow $\Phi_t$ on $\Omega$. Suppose that
$f\in\mathcal O(\Omega)\cap C(\overline\Omega)$
is nonconstant on some real orbit.

Then, for every $K>0$, there exist nonzero
smooth positive closed real $(1,1)$-forms in
$\mathcal T_{Y,K}(\Omega)$ and
$\mathcal T_{-Y,K}(\Omega)$.
\end{proposition}

\begin{proof}
By hypothesis, there exist $z_0\in\Omega$
and $\tau\in\R$ such that
$f(\Phi_\tau(z_0))\ne f(z_0)$. Hence
\begin{equation}\label{eq:scalar-displacement}
h_\tau=f\circ\Phi_\tau-f
\end{equation}
is a nonzero holomorphic function on $\Omega$.

Let $L\Subset\Omega$ and let $|s_j|\to\infty$.
Properness of the real action implies that
the sequence $\Phi_{s_j}$ is compactly
divergent. By Lemmas~\ref{lem:ambient-normal}
and \ref{lem:stein-normal-limit}, after
passing to a subsequence, there exists
$q\in\partial_Z\Omega$ such that
$\Phi_{s_j}\to q$ locally uniformly in $Z$.

Put $L_\tau=L\cup\Phi_\tau(L)$. The set
$L_\tau$ is compact in $\Omega$, and
continuity of $f$ on $\overline\Omega$
gives
\[
\sup_{w\in L_\tau}
|f(\Phi_{s_j}(w))-f(q)|\longrightarrow0.
\]
Using the group law
\eqref{eq:real-action-group-law}, we obtain
\begin{equation}\label{eq:scalar-displacement-estimate}
\begin{aligned}
\sup_{z\in L}|h_\tau(\Phi_{s_j}(z))|
&=\sup_{z\in L}
|f(\Phi_{s_j}(\Phi_\tau(z)))
       -f(\Phi_{s_j}(z))|\\
&\le 2\sup_{w\in L_\tau}
|f(\Phi_{s_j}(w))-f(q)|
\longrightarrow0.
\end{aligned}
\end{equation}

Suppose that $h_\tau\circ\Phi_s$ does
not converge to zero locally uniformly
as $|s|\to\infty$. Then there exist
$L\Subset\Omega$, $\varepsilon>0$,
and a sequence $|s_j|\to\infty$ such that
\[
\sup_{z\in L}|h_\tau(\Phi_{s_j}(z))|
\ge\varepsilon
\qquad(j\ge1).
\]
Passing to a subsequence with a constant
boundary limit contradicts
\eqref{eq:scalar-displacement-estimate}.
Therefore
\begin{equation}\label{eq:scalar-two-ended-decay}
h_\tau\circ\Phi_s\longrightarrow0
\quad\text{locally uniformly as }s\to\pm\infty.
\end{equation}
Since $h_\tau\not\equiv0$,
Lemma~\ref{lem:two-ended-holomorphic}
gives the asserted nonzero smooth
positive closed forms.
\end{proof}

If $\Omega$ is biholomorphic to a relatively compact Stein domain
whose closure has a neighborhood admitting a holomorphic
immersion into Euclidean space, transport the flows and forms
to that realization. By \cite[Theorem~1.2]{GangulyKRsol2026Paper21},
complete nontrivial expanders exist there. Pullback gives
the expanders for $Y$ and $-Y$ on $\Omega$.

\begin{proof}[Proof of Theorem~\ref{thm:intro-manifold-geometric}]
Let $\iota=(f_1,\ldots,f_N)$ be a holomorphic
immersion defined on a neighborhood of
$\overline\Omega$. Since $\overline\Omega$ is
compact, each $f_j$ belongs to
$\mathcal O(\Omega)\cap C(\overline\Omega)$.

Under hypothesis \textup{(i)}, the immersion gives
\eqref{eq:stein-bounded-separation} by
Lemma~\ref{lem:ambient-normal}. Under \textup{(ii)},
the visibility assumption is precisely that of
Proposition~\ref{prop:geometric-beyond-immersion}\textup{(ii)}.

We show that at least one component $f_j$
is nonconstant along a real orbit.
Fix $x\in\Omega$ and define
\begin{equation}\label{eq:immersion-finite-fiber}
E_x=\overline\Omega\cap
\iota^{-1}(\iota(x)).
\end{equation}
The set $E_x$ is compact. By the holomorphic
constant-rank theorem, $\iota$ is injective
on a neighborhood of each point of $E_x$.
Consequently, $E_x$ is discrete and finite.

Suppose that
\[
\iota\circ\Phi_1=\iota.
\]
By the group law
\eqref{eq:real-action-group-law},
\[
\iota(\Phi_m(x))=\iota(x)
\qquad(m\in\Z).
\]
Thus
\[
\{\Phi_m(x):m\in\Z\}\subset E_x.
\]
A proper $\R$-action has compact
stabilizers. Since $\R$ has no nontrivial
compact subgroup, the action is free.
Therefore the points $\Phi_m(x)$,
$m\in\Z$, are pairwise distinct,
contradicting the finiteness of $E_x$.

Hence there exists $j_0\in\{1,\ldots,N\}$
such that
\[
f_{j_0}\circ\Phi_1-f_{j_0}\not\equiv0.
\]
In particular, $f_{j_0}$ is nonconstant
along some real orbit.
Proposition~\ref{prop:geometric-beyond-immersion}
now gives nonzero smooth positive closed
real $(1,1)$-forms in
$\mathcal T_{Y,K}(\Omega)$ and
$\mathcal T_{-Y,K}(\Omega)$
for every $K>0$.

For each $1\le j\le N$, the estimates
\eqref{eq:scalar-displacement-estimate} and
\eqref{eq:scalar-two-ended-decay} apply with
$f=f_j$ and $\tau=1$, even when
$f_j\circ\Phi_1-f_j\equiv0$. Hence
\begin{equation}\label{eq:unit-displacement-decay}
\iota\circ\Phi_{s+1}
-\iota\circ\Phi_s\longrightarrow0
\qquad(s\to\pm\infty)
\end{equation}
locally uniformly on $\Omega$.

Define
\begin{equation}\label{eq:immersion-displacement-seed}
\begin{aligned}
u(z)
&=|\iota(\Phi_1(z))-\iota(z)|^2,\\
u(\Phi_s(z))
&=|\iota(\Phi_{s+1}(z))
       -\iota(\Phi_s(z))|^2.
\end{aligned}
\end{equation}
Since $\iota(\overline\Omega)$ is compact,
$u$ is bounded. As the squared norm
of a holomorphic map into $\C^N$,
$u$ is nonnegative and plurisubharmonic.
The choice of $j_0$ shows that
$u\not\equiv0$, although $u$ may
vanish at some points because
$\iota$ need not be injective.

By \eqref{eq:unit-displacement-decay}
and \eqref{eq:immersion-displacement-seed},
\[
u\circ\Phi_s\longrightarrow0
\quad\text{locally uniformly as }s\to\pm\infty.
\]
Consequently, if $\Omega$ is Stein,
Theorem~\ref{thm:backward-psh-prescribed-intro}
applies to $u$ for both $\Phi_s$
and the reversed flow $\Phi_{-s}$.

In particular, when $\Omega\Subset\C^n$
and $\iota=\id_\Omega$, the function
in \eqref{eq:immersion-displacement-seed}
reduces to
\[
u(z)=|\Phi_1(z)-z|^2>0.
\]
The strict positivity follows from
freeness of the real action, and
\[
u(\Phi_s(z))
=|\Phi_{s+1}(z)-\Phi_s(z)|^2
\longrightarrow0
\qquad(s\to\pm\infty)
\]
locally uniformly on $\Omega$
by \eqref{eq:unit-displacement-decay}.

Suppose now that $\Omega$ is Stein.
The nonzero homogeneous forms obtained
above, together with
\cite[Theorem~1.2]{GangulyKRsol2026Paper21},
give complete nontrivial expanders
with soliton fields $Y$ and $-Y$
for every $K>0$.

If $\Aut_0(\Omega)$ is noncompact,
Lemmas~\ref{stein:lem:bounded} and
\ref{stein:lem:realflow}
supply a proper real holomorphic flow. The homogeneous-current construction therefore gives a complete nontrivial expander.

Conversely, the existence of a complete
nontrivial expander implies that
$\Omega$ is Stein by
\cite[Remarks~3.6(2)]{GangulyKRsol2026Paper21}.
The soliton flow has noncompact closure
in $\Aut(\Omega)$ by
\cite[Lemma~3.1]{GangulyKRsol2026Paper21}.
Since the soliton flow lies in
$\Aut_0(\Omega)$, the group
$\Aut_0(\Omega)$ is noncompact.
\end{proof}

\begin{proof}[Proof of Corollary~\ref{cor:intro-global-finite-type}]
Finite D'Angelo $1$-type excludes a nonconstant analytic
disc in the boundary. The immersion hypothesis gives the
proper Lie-group action of $\Aut(\Omega)$ by
Lemma~\ref{stein:lem:bounded}. Hence a complete field with
noncompact flow closure generates a proper action by
Lemma~\ref{stein:lem:realflow}. Apply
Theorem~\ref{thm:intro-manifold-geometric}\textup{(i)}.
The converse is \cite[Lemma~3.1]{GangulyKRsol2026Paper21}.
\end{proof}

For a Kobayashi hyperbolic domain, local weak visibility
and hyperbolic embedding imply relative weak visibility
\cite[Theorem~1.6]{Masanta2025Visibility}.

If $\Omega\Subset Z$, where $Z$ is Stein or Kobayashi hyperbolic, relative weak visibility also follows if the set of
points $p\in\partial_Z\Omega$ near which $\Omega$ is not weakly
visible in the sense of \cite[Definition~1.9]{Masanta2025Visibility} is totally disconnected
\cite[Theorem~1.10(ii)]{Masanta2025Visibility}.

\begin{remark}\label{rem:single-kappa-manifold}
Let $\Omega\Subset Z$ satisfy relative $(1,\kappa_0)$-visibility
for some $\kappa_0>0$, and let $\Phi_t$ be a proper real
holomorphic flow with generator $Y$. Fix a smooth Hermitian
metric on $Z$ with distance $d$. For $F=\Phi_1$, properness gives
\eqref{eq:weakvis-compact-divergence} for $F$ and $F^{-1}$.
Lemma~\ref{lem:weak-visibility-iterates} gives limits
$q_\pm\in\partial_Z\Omega$ for the integer maps $\Phi_m$ as
$m\to\pm\infty$. For $L\Subset\Omega$, the set
$L'=\bigcup_{0\le r\le1}\Phi_r(L)$ is compact.
For $t=m+r$, with $m\in\Z$ and $0\le r<1$, we have
\begin{equation}\label{eq:visibility-real-time}
\sup_{z\in L}d(\Phi_t(z),q_\pm)
\le \sup_{w\in L'}d(\Phi_m(w),q_\pm)\longrightarrow0
\qquad(t\to\pm\infty).
\end{equation}

If a neighborhood of $\overline\Omega$ admits a holomorphic
immersion into Euclidean space, let $f$ be a nonconstant
component of that immersion and put
\[
u_-(z)=|f(z)-f(q_-)|^2,\qquad
u_+(z)=|f(z)-f(q_+)|^2.
\]
Both functions are bounded, nonnegative, nonzero and
plurisubharmonic on $\Omega$. By
\eqref{eq:visibility-real-time},
\[
u_-\circ\Phi_s\longrightarrow0,\qquad
u_+\circ\Phi_{-s}\longrightarrow0
\quad\text{locally uniformly as }s\to-\infty.
\]
Thus, if $\Omega$ is Stein,
Theorem~\ref{thm:backward-psh-prescribed-intro}
gives complete nontrivial expanders with fields $Y$
and $-Y$, respectively, for every $K>0$.
\end{remark}

\begin{remark}
Let $\Omega\Subset Z$ be such that every sequence of holomorphic
self-maps of $\Omega$ has a subsequence converging locally
uniformly in $Z$ to a holomorphic map into $\overline\Omega$,
as in Lemma~\ref{lem:ambient-normal}. Let $\Phi_t$ be a proper
real holomorphic flow. Suppose that for some $\tau\ne0$,
$g\circ\Phi_\tau=g$ for every $g\in\mathscr N_-(\Phi)$.
If $f\in\mathcal O(\Omega)\cap C(\overline\Omega)$ and
$f\circ\Phi_\tau\ne f$, put $u=|f\circ\Phi_\tau-f|^2$.
Then $u$ is bounded, nonnegative, plurisubharmonic, and
nonzero.

If $u\circ\Phi_s$ did not tend to zero locally uniformly as
$s\to-\infty$, there would exist a compact $L\Subset\Omega$,
$\epsilon>0$, and $z_j\in L$, $s_j\to-\infty$ with
$u(\Phi_{s_j}z_j)\ge\epsilon$. Pass to a subsequence with
$z_j\to z\in L$ and $\Phi_{s_j}\to g\in\mathscr N_-(\Phi)$
locally uniformly. Using convergence on $L\cup\Phi_\tau(L)$,
\[
u(\Phi_{s_j}z_j)
=|f(\Phi_{s_j}(\Phi_\tau z_j))-f(\Phi_{s_j}z_j)|^2
\longrightarrow |f(g(\Phi_\tau z))-f(g(z))|^2=0 .
\]
Hence $u\circ\Phi_s\to0$ locally uniformly as $s\to-\infty$.
Choose $0<\delta<u(a)$ at a point $a$ where $u(a)>0$. Then
$(u-\delta)_+$ satisfies \eqref{eq:backward-cutoff-vanishing},
and Lemma~\ref{lem:backward-cutoff-current} gives a nonzero
current in $\mathcal T_{Y,K}(\Omega)$ for every $K>0$, where
$Y$ generates $\Phi_t$. The forward version is the same
assertion for $-Y$.

For a holomorphic immersion near $\overline\Omega$, the same
finite-fiber argument as for the time-one map shows that some
component $f$ is moved by $\Phi_\tau$: otherwise the orbit
$\{\Phi_{m\tau}(z):m\in\Z\}$ lies in one finite fiber. On a
bounded Euclidean domain one can instead use
$u_\tau(z)=|\Phi_\tau(z)-z|^2$. Freeness makes this function
strictly positive, and applying the displayed argument
componentwise gives $u_\tau\circ\Phi_s\to0$ locally uniformly
as $s\to-\infty$.
\end{remark}

Recall a subset \(E\) of a complex manifold \(U\) is complete pluripolar
if \(E=\{\psi=-\infty\}\) for some \(\psi\in\PSH(U)\).

\begin{corollary}\label{prop:stein-pluripolar}\label{thm:pluripolar-backward-intro}
\label{cor:pluripolar-backward}
Let $\Omega\Subset Z$ be a domain such that every
sequence of holomorphic self-maps of $\Omega$
admits a subsequence converging locally uniformly
in $Z$ to a holomorphic map
$\Omega\to\overline\Omega$, and let $\Phi_t$
be a proper real holomorphic flow on $\Omega$.
Suppose a neighborhood \(U\) of \(\overline\Omega\) contains a
complete pluripolar set \(E\) disjoint from \(\Omega\) and containing
every image \(g(\Omega)\) for \(g\in\mathscr N_-(\Phi)\).
Then, for every $K>0$, there exists a nonzero current
in $\mathcal T_{Y,K}(\Omega)$, where $Y$ is the generator
of $\Phi_t$. If, in addition, $\Omega$ is Stein and a
neighborhood of $\overline\Omega$ admits a holomorphic
immersion into Euclidean space, then $\Omega$ admits
a complete nontrivial expanding K\"ahler--Ricci soliton
with soliton field $Y$ and expansion constant $K$.
\end{corollary}
\begin{proof}
Choose $\psi\in\PSH(U)$ with $E=\{\psi=-\infty\}$ and an
open set $U_0$ with $\overline\Omega\subset U_0\Subset U$.
Choose $c$ with $\psi\le c$ on $U_0$. Then
$\beta=e^{\psi-c}\in\PSH(U_0)$ is bounded, vanishes on the
limit images, and is positive on $\Omega$.
Remark~\ref{rem:normal-limit-psh-barrier}
therefore gives a nonzero current in
$\mathcal T_{Y,K}(\Omega)$ for every $K>0$.
Under the additional Stein and immersion assumptions,
\cite[Theorem~1.2]{GangulyKRsol2026Paper21}
then gives a complete nontrivial expander with
soliton field $Y$ and expansion constant $K$.
\end{proof}

For a related discrete-quotient criterion on a bounded Stein
domain in $\C^n$, see \cite[Theorem~5]{Chen2008ProperActions}.
It uses a negative continuous plurisubharmonic function that
tends to $-\infty$ at the union of the orbit cluster sets for
points in a fundamental domain.

\subsection{Bounded convex and homogeneous domains}
\label{sec:geometric-criteria}

The supporting hyperplanes of a bounded convex domain $\Omega$
control the variation of its diagonal Bergman kernel $K_\Omega(z,z)$
in the Kobayashi distance.

\begin{lemma}\label{lem:convex-kernel}
Let $\Omega\Subset\C^n$ be a bounded convex domain, and let $k_\Omega$ denote its Kobayashi distance. Then
\begin{equation}\label{eq:convex-kernel}
\left|\log K_\Omega(z,z)-\log K_\Omega(w,w)\right|
\le 4n\,k_\Omega(z,w)
\qquad(z,w\in\Omega).
\end{equation}
Consequently, every $F\in\Aut(\Omega)$ satisfies
\begin{equation}\label{eq:convex-jac-distortion}
\left|\log|\det_{\C}DF(z)|-\log|\det_{\C}DF(w)|\right|
\le 4n\,k_\Omega(z,w).
\end{equation}
\end{lemma}

\begin{proof}
Let $\ell$ be a nonzero real-linear functional on $\C^n$, and set
\[
h(\ell)=\sup_{\zeta\in\Omega}\ell(\zeta),
\qquad
d_\ell(z)=h(\ell)-\ell(z).
\]
Since $\Omega$ is bounded and open and $\ell$ is nonzero, we have $0<d_\ell(z)<\infty$ for every $z\in\Omega$. Choose a complex-linear functional $L$ such that $\operatorname{Re}L=\ell$, and let $\mathbb H=\{\zeta\in\C:\operatorname{Re}\zeta>0\}$. The holomorphic map $f=h(\ell)-L$ sends $\Omega$ into $\mathbb H$, with $\operatorname{Re}f=d_\ell$. The distance-decreasing property of the Kobayashi metric gives
\[
k_{\mathbb H}(f(z),f(w))\le k_\Omega(z,w).
\]
For the Cayley transform $C(\zeta)=(\zeta-1)/(\zeta+1)$,
we have $C'(\zeta)=2/(\zeta+1)^2$ and
$1-|C(\zeta)|^2=4\Re\zeta/|\zeta+1|^2$. Therefore the
Poincar\'e metric in the disc normalization of
\cite[Chapter~1]{KrantzGeometricViewpoint2004} gives
\[
\kappa_{\mathbb H}(\zeta;v)=
\frac{|C'(\zeta)v|}{1-|C(\zeta)|^2}
=\frac{|v|}{2\Re\zeta}.
\]
Consequently, integration along piecewise smooth curves in $\mathbb H$ yields
\[
k_{\mathbb H}(\zeta_1,\zeta_2)\ge
\frac12\left|\log\operatorname{Re}\zeta_1-
\log\operatorname{Re}\zeta_2\right|.
\]
Therefore,
\[
\left|\log d_\ell(z)-\log d_\ell(w)\right|
\le 2k_\Omega(z,w).
\]
Writing $d=k_\Omega(z,w)$, we obtain
\begin{equation}\label{eq:convex-gap}
e^{-2d}d_\ell(z)\le d_\ell(w)\le e^{2d}d_\ell(z).
\end{equation}
Since $\Omega$ is open and convex, the supporting-hyperplane theorem gives
\[
\Omega=\bigcap_{\ell\ne0}\{u\in\C^n:\ell(u)<h(\ell)\},
\]
where the intersection is taken over all nonzero real-linear functionals $\ell$. For every $u\in\Omega$ and every such $\ell$,
\[
\begin{aligned}
\ell\bigl(w+e^{-2d}(u-z)\bigr)
&=\ell(w)+e^{-2d}\bigl(\ell(u)-\ell(z)\bigr) <\ell(w)+e^{-2d}d_\ell(z) \le\ell(w)+d_\ell(w)=h(\ell).
\end{aligned}
\]
Consequently, $e^{-2d}(\Omega-z)\subset\Omega-w$. Interchanging $z$ and $w$ yields
\begin{equation}\label{eq:convex-inclusions}
e^{-2d}(\Omega-z)\subset\Omega-w
\subset e^{2d}(\Omega-z).
\end{equation}
The preceding inclusions and monotonicity of the Bergman kernel under domain inclusion give
\[
K_{e^{2d}(\Omega-z)}(0,0)
\le K_{\Omega-w}(0,0)
\le K_{e^{-2d}(\Omega-z)}(0,0).
\]
Recall that for a biholomorphism $\Psi:U\to V$, the change of variables for square-integrable canonical forms gives
\begin{equation}\label{eq:convex-bergman-transformation}
K_V(\Psi(\zeta),\Psi(\zeta))
|\det_{\C}D\Psi(\zeta)|^2
=K_U(\zeta,\zeta).
\end{equation}
For $\Psi(\zeta)=r(\zeta-z)$, \eqref{eq:convex-bergman-transformation} gives $K_{r(\Omega-z)}(0,0)=r^{-2n}K_\Omega(z,z)$, while translation invariance gives $K_{\Omega-w}(0,0)=K_\Omega(w,w)$. Thus the preceding inequalities become
\[
e^{-4nd}K_\Omega(z,z)
\le K_\Omega(w,w)
\le e^{4nd}K_\Omega(z,z).
\]
Taking logarithms, we obtain
\[
\left|\log K_\Omega(z,z)-\log K_\Omega(w,w)\right|
\le4nd=4n\,k_\Omega(z,w),
\]
which proves \eqref{eq:convex-kernel}. For $F\in\Aut(\Omega)$, \eqref{eq:convex-bergman-transformation} gives
\[
K_\Omega(F(z),F(z))|\det_{\C}DF(z)|^2
=K_\Omega(z,z).
\]
Applying \eqref{eq:convex-kernel} to $(z,w)$ and $(F(z),F(w))$, and using the invariance of the Kobayashi distance under $F$, proves \eqref{eq:convex-jac-distortion}.
\end{proof}

\begin{proof}[Proof of Theorem~\ref{thm:geometric-prescribed}]
Let $F:M\to D$ be the given biholomorphism. The conjugate
flow $\Psi_t=F\circ\Phi_t\circ F^{-1}$ is proper and has
generator $F_*Y$. Fix $K>0$.

If $D$ is convex, let $j_t$ be the logarithmic Jacobian of
$\Psi_t$ from Lemma~\ref{lem:jacobian-log}. For $a\in D$,
\[
|\Re(j_t(z)-j_t(a))|\le4n\,k_D(z,a)
\qquad(z\in D,\ t\in\R)
\]
by \eqref{eq:convex-jac-distortion}. This gives
\eqref{eq:halfline-jacobian} for each of $\Psi_t$ and
$\Psi_{-t}$. Apply Proposition~\ref{prop:frame-halfline}
with the coordinate frame and the Bergman kernel form to get
nonzero smooth forms in $\mathcal T_{F_*Y,K}(D)$
and $\mathcal T_{-F_*Y,K}(D)$.

If $D$ is homogeneous, it is Stein by
\cite[Proposition~2.7 and Theorem~2.8]{MiebachHomogeneous2010}.
The subgroup $\langle\Psi_1\rangle$ lies in $\Aut_0(D)$
and is discrete, since it is the integer subgroup of a
proper real action. Its quotient is Stein by
\cite[Theorem~6.5]{MiebachHomogeneous2010}. Let $H$ be the time
function given by Proposition~\ref{prop:stein-cyclic-character}
for $\Psi_t$. The forms $\ddbar e^{\pm K\Re H}$ have the
required characters by \eqref{eq:stein-cyclic-time}.

In both cases $D$ is Stein. \cite[Theorem~1.2]{GangulyKRsol2026Paper21}
gives complete nontrivial expanders on $D$ with fields
$F_*Y$ and $-F_*Y$. Pullback by $F$ gives the forms and
complete solitons on $M$ with fields $Y$ and $-Y$.

Conjugation by $F$ and Lemma~\ref{stein:lem:bounded} make
$\Aut(M)$ a Lie group acting properly. If $\Aut_0(M)$
is noncompact, Lemma~\ref{stein:lem:realflow} gives a proper
real flow, and the construction just proved applies.
Conversely, a complete nontrivial expander has soliton flow
with noncompact closure by \cite[Lemma~3.1]{GangulyKRsol2026Paper21};
the flow lies in $\Aut_0(M)$.
\end{proof}

\section{Uniformization and Finite-Type Domains}
\label{sec:finite-type-extensions}

Let $\Omega\Subset Z$ have a $C^2$ strongly pseudoconvex
boundary orbit accumulation point, with $\dim_\C Z=n\ge2$.
By \cite[Theorem~1]{GaussierKimKrantz2002}, there is a
biholomorphism $F:\Omega\to\B^n$. If $Y$ is a complete real
holomorphic field on $\Omega$ whose flow has noncompact closure,
then the flow of $F_*Y$ is proper by
Lemmas~\ref{stein:lem:bounded} and~\ref{stein:lem:realflow}.
For every $K>0$, Theorem~\ref{thm:geometric-prescribed} gives
nonzero smooth forms $T_\pm\in\mathcal T_{\pm F_*Y,K}(\B^n)$
and complete nontrivial expanders with fields $\pm F_*Y$.
Pullback by $F$ gives
\[
F^*T_\pm\in\mathcal T_{\pm Y,K}(\Omega)
\]
and complete nontrivial expanders with fields $\pm Y$.

\subsection{Finite-type points on complex surfaces}

The following prescribed-field result strengthens the
imaginary-translation construction of
\cite[Corollary~5.2]{GangulyKRsol2026Paper21}; compare also
\cite[Corollary~5.3]{GangulyKRsol2026Paper21}.

\begin{theorem}\label{ft:thm:two-d-prescribed}
Let $\Omega\Subset Z$ be a domain in a complex surface, and suppose
that a neighborhood of $\overline\Omega$ admits a holomorphic
immersion into $\C^N$. Suppose that $\Omega$ has a boundary orbit
accumulation point $p\in\partial_Z\Omega$ at which the boundary is
$C^\infty$-smooth, pseudoconvex, and of finite D'Angelo $1$-type
in a neighborhood. Then $\Omega$ is biholomorphic to
\begin{equation}\label{ft:eq:model}
M_H=\{(z,w)\in\C^2:\Re w+H(z,\bar z)<0\},
\end{equation}
where $H$ is a nonzero real homogeneous subharmonic polynomial of
even degree $2m$ without harmonic terms. Moreover,
\begin{equation}\label{ft:eq:surface-aut-dimensions}
\dim_\R\Aut(\Omega)\in\{2,3,4,8\}.
\end{equation}
For every complete real holomorphic field $Y$ whose flow has
noncompact closure and every $K>0$, the cones
$\mathcal T_{Y,K}(\Omega)$ and $\mathcal T_{-Y,K}(\Omega)$
contain nonzero smooth positive closed real $(1,1)$-forms of
rank one at every point. Both $Y$ and $-Y$ occur as fields of
complete nontrivial expanders with expansion constant $K$.
Conversely, if $Y$ is the field of a complete nontrivial
expander, or if $\mathcal T_{Y,K}(\Omega)$ contains a nonzero
current for some $K>0$, then its flow has noncompact closure.
\end{theorem}

\begin{proof}
By \cite[Theorem~1.1]{SeshadriVerma2008}, $\Omega$ is
biholomorphic to
$\{(z,w):2\Re w+H(z,\bar z)<0\}$, with $H$ as in the
statement. Rescaling the $w$-coordinate gives
\eqref{ft:eq:model}. The model is Stein by
\cite[proof of Corollary~5.2]{GangulyKRsol2026Paper21}.

The automorphisms
\[
\delta_t(z,w)=(e^{t/(2m)}z,e^tw)
\]
satisfy $\delta_t(0,-1)\to(0,0)$ as $t\to-\infty$.
Put $\rho=\Re w+H(z,\bar z)$. Since $H$ is subharmonic,
$\partial M_H$ is smooth and pseudoconvex. For a nonconstant
holomorphic curve germ $\gamma=(z,w)$ at $(0,0)$ with
$z\not\equiv0$ and $w\not\equiv0$, write
$k=\operatorname{ord}_0z$ and $\ell=\operatorname{ord}_0w$.
The terms of degree $2mk$ in
$H(z(\zeta),\overline{z(\zeta)})$ are mixed monomials in
$\zeta,\bar\zeta$, whereas those of degree $\ell$ in
$\Re w(\zeta)$ are harmonic. Since $H$ has no harmonic
terms, they cannot cancel. Hence
\begin{equation}\label{ft:eq:surface-contact-type}
\frac{\operatorname{ord}_0(\rho\circ\gamma)}
{\operatorname{ord}_0\gamma}
=\frac{\min\{\ell,2mk\}}{\min\{k,\ell\}}\le2m.
\end{equation}
For $z\equiv0$ the quotient equals $1$; for $w\equiv0$
it equals $2m$. The curve $\gamma(\zeta)=(\zeta,0)$
therefore shows that $\Delta_1(\partial M_H,0)=2m$.
The type is locally bounded near the origin by
\cite[Theorem~6.5]{DAngeloKohn1999}, and
\cite[Corollary~1]{ThuDuc2016} gives
\eqref{ft:eq:surface-aut-dimensions}.

Since $\Omega$ is Stein, Corollary~\ref{thm:stein-surface-intro}
gives the rank-one forms and complete expanders for every $Y$
with noncompact flow closure. For a nontrivial soliton field,
noncompactness follows from
\cite[Lemma~3.1]{GangulyKRsol2026Paper21}.

Finally, suppose $0\ne T\in\mathcal T_{Y,K}(\Omega)$ and that
the flow $\Phi_t$ of $Y$ has compact closure. By the argument in \cite[proof of Lemma~3.1]{GangulyKRsol2026Paper21},
there are $s_j\to+\infty$ with $\Phi_{s_j}\to\id_\Omega$ in
$\Aut(\Omega)$. Since $\Phi_{-s_j}\to\id_\Omega$
smoothly on compact sets, and the pullbacks of any
fixed compactly supported test form have supports
in a common compact set, we have
$\Phi_{s_j}^*T\rightharpoonup T$. Choose
a compactly supported smooth test form $\beta$ with
$T(\beta)\ne0$. Then
\[
T(\beta)=\lim_{j\to\infty}(\Phi_{s_j}^*T)(\beta)
=\lim_{j\to\infty}e^{Ks_j}T(\beta),
\]
a contradiction.
\end{proof}

\begin{corollary}\label{ft:cor:topology}
Let $\Omega\Subset Z$ be a domain in a complex surface
with a boundary orbit accumulation point
$p\in\partial_Z\Omega$ at which the boundary is
$C^\infty$-smooth, pseudoconvex, and of finite
D'Angelo $1$-type in a neighborhood. Then
$\Omega$ is diffeomorphic to $\R^4$.
In particular, $\Omega$ is contractible.
\end{corollary}

\begin{proof}
By \cite[Theorem~1.1]{SeshadriVerma2008},
after rescaling the $w$-coordinate, $\Omega$ is
biholomorphic to $M_H$ in \eqref{ft:eq:model}.
The map
\begin{equation}\label{ft:eq:surface-diffeomorphism}
\C\times(-\infty,0)\times\R\longrightarrow M_H,
\qquad(z,r,t)\longmapsto(z,r-H(z,\bar z)+\ii t)
\end{equation}
is a diffeomorphism, with inverse
$(z,w)\mapsto(z,\Re w+H(z,\bar z),\Im w)$.
\end{proof}

\subsection{Egg domains}

\begin{proof}[Proof of Corollary~\ref{eggdomain_char}]
Suppose that $\Omega$ admits a complete nontrivial expander
with complete soliton field. By
\cite[Theorem~1.1]{GangulyKRsol2026Paper21}, $\Omega$ is
pseudoconvex, and by
\cite[Lemma~3.1]{GangulyKRsol2026Paper21}, $\Aut(\Omega)$
is noncompact. The local boundedness of D'Angelo type
\cite[Theorem~6.5]{DAngeloKohn1999} and compactness of
$\partial\Omega$ give $M<\infty$ with
\begin{equation}\label{ft:eq:egg-uniform-type}
\Delta_1(\partial\Omega,p)\le M
\qquad(p\in\partial\Omega).
\end{equation}
If an analytic set $V$ of positive dimension at
$p\in\partial\Omega$ satisfies $|\rho(z)|\le C|z-p|^q$ for a local defining function
$\rho$, choose a nonconstant holomorphic map
$\gamma:(\Delta,0)\to(V,p)$. Setting
$\nu=\operatorname{ord}_0(\gamma-p)>0$, we obtain
\begin{equation}\label{ft:eq:egg-contact-bound}
q\nu\le\operatorname{ord}_0(\rho\circ\gamma)
\le\Delta_1(\partial\Omega,p)\nu\le M\nu.
\end{equation}
Thus the contact order is uniformly bounded. The main theorem
of \cite{BedfordPinchukEgg1991} gives $\Omega\simeq E_m$ for
some $m\ge1$.

Conversely, let $F:\Omega\to E_m$ be a biholomorphism and
write $z'=(z_2,\ldots,z_n)$. Define
\begin{equation}\label{ft:eq:egg-cayley}
C(w,z_1,z')=
\left(\ii\frac{1+w}{1-w},\,
\frac{z_1}{(1-w)^{1/m}},\,
\frac{z'}{1-w}\right)=:(\zeta,q_1,q').
\end{equation}
The branch of $(1-w)^{1/m}$ is defined using the
holomorphic logarithm of $1-w$ on $|w|<1$. We have
\begin{equation}\label{ft:eq:egg-cayley-identity}
\Im\zeta=\frac{1-|w|^2}{|1-w|^2},\qquad
|q_1|^{2m}+|q'|^2
=\frac{|z_1|^{2m}+|z'|^2}{|1-w|^2}.
\end{equation}
Hence $C$ maps $E_m$ into
\[
\mathcal S_m=\{(\zeta,q_1,q'):
\Im\zeta>|q_1|^{2m}+|q'|^2\}.
\]
For $(\zeta,q_1,q')\in\mathcal S_m$, the inverse is
\[
w=\frac{\zeta-\ii}{\zeta+\ii},\qquad
z_1=(1-w)^{1/m}q_1,\qquad z'=(1-w)q'.
\]
The inequality $\Im\zeta>0$ gives $|w|<1$ and
\eqref{ft:eq:egg-cayley-identity} gives the defining
inequality for $E_m$. Thus $C:E_m\to\mathcal S_m$ is a
biholomorphism.

The translations
\begin{equation}\label{ft:eq:egg-translation-flow}
\Phi_t(\zeta,q_1,q')=(\zeta+t,q_1,q')
\end{equation}
define a proper holomorphic $\R$-action on $\mathcal S_m$:
if $x_j$ and $\Phi_{t_j}(x_j)$ remain in fixed compact sets,
then $\Re\zeta(x_j)$ and $\Re\zeta(x_j)+t_j$ are bounded,
so $(t_j)$ is bounded. For $K>0$, set
\begin{equation}\label{ft:eq:egg-explicit-forms}
T_{K,\pm}=\ddbar e^{\pm K\Re\zeta}
=\frac{K^2}{4}e^{\pm K\Re\zeta}\,
\ii\,d\zeta\wedge d\bar\zeta.
\end{equation}
These smooth positive closed forms have rank one, and
\begin{equation}\label{ft:eq:egg-character}
\Phi_t^*T_{K,+}=e^{Kt}T_{K,+},\qquad
\Phi_t^*T_{K,-}=e^{-Kt}T_{K,-}.
\end{equation}
With $G=C\circ F$ and $Y_0=G^*(\partial/\partial\Re\zeta)$,
\eqref{ft:eq:egg-character} gives
\[
G^*T_{K,+}\in\mathcal T_{Y_0,K}(\Omega),\qquad
G^*T_{K,-}\in\mathcal T_{-Y_0,K}(\Omega).
\]
By \cite[Theorem~1.2]{GangulyKRsol2026Paper21},
$Y_0$ and $-Y_0$ are soliton fields of complete nontrivial
expanders for every $K>0$.

If $Y$ is any complete real holomorphic field on $\Omega$
with noncompact flow closure, the flow of $F_*Y$ on $E_m$
is proper by Lemmas~\ref{stein:lem:bounded} and
\ref{stein:lem:realflow}. Since $E_m$ is convex,
Theorem~\ref{thm:geometric-prescribed} gives the homogeneous
forms and complete expanders for both $Y$ and $-Y$.
\end{proof}

In $\C^2$, Corollary~\ref{eggdomain_char} applies to bounded
domains with smooth finite-type boundary without a Levi-rank
assumption.

\subsection{Weighted polynomial models}

\begin{proposition}\label{ft:prop:weighted}
Let $\Omega\Subset Z$ be a domain of complex dimension $n\ge2$
such that a neighborhood of $\overline\Omega$ admits a
holomorphic immersion into $\C^N$. Suppose that
$F:\Omega\to M_P$ is a biholomorphism, where
\begin{equation}\label{ft:eq:weighted-model}
M_P=\{(z,w)\in\C^{n-1}\times\C:
\rho(z,w):=\Re w+P(z,\bar z)<0\},
\end{equation}
$P$ is a real-valued plurisubharmonic polynomial, and
\begin{equation}\label{ft:eq:weighted}
P(e^{\lambda_1t}z_1,\ldots,e^{\lambda_{n-1}t}z_{n-1})
=e^tP(z),\qquad t\in\R,
\end{equation}
for some $\lambda_1,\ldots,\lambda_{n-1}>0$.
Define
\begin{equation}\label{ft:eq:weighted-flows}
\begin{aligned}
\delta_t(z,w)&=(e^{\lambda_1t}z_1,\ldots,
 e^{\lambda_{n-1}t}z_{n-1},e^tw),\\
\tau_t(z,w)&=(z,w+\ii t).
\end{aligned}
\end{equation}
Let $D$ and $A$ denote the real infinitesimal
generators of $\delta_t$ and $\tau_t$, respectively.
For $a,b\in\R$ with $(a,b)\ne(0,0)$ and every $K>0$,
$\mathcal T_{F^*(aD+bA),K}(\Omega)$ contains a nonzero
positive closed real $(1,1)$-current, and $F^*(aD+bA)$ is
the field of a complete nontrivial expander with
expansion constant $K$. Moreover, $\Omega$ is diffeomorphic
to $\R^{2n}$.
\end{proposition}

\begin{proof}
By \eqref{ft:eq:weighted},
\begin{equation}\label{ft:eq:weighted-defining-character}
\rho\circ\delta_t=e^t\rho,\qquad
\rho\circ\tau_t=\rho.
\end{equation}
In the coordinates with $w'=w+\ii r$, the $(1,0)$-part of
$(\tau_r)_*D$ has $w'$ coefficient $w'-\ii r$, and that of
$A$ has coefficient $\ii$. Consequently,
\begin{equation}\label{ft:eq:weighted-conjugation}
(\tau_r)_*(aD+bA)=aD+(b-ar)A.
\end{equation}
If $a\ne0$, the flow of $aD+bA$ is conjugate by
$\tau_{b/a}$ to $\delta_{as}$; if $a=0$, the flow
is $\tau_{bs}$. For compact sets $C_1,C_2\Subset M_P$,
$x\in C_1$ and $\delta_{as}(x)\in C_2$ imply
\[
\frac{\min_{C_2}(-\rho)}{\max_{C_1}(-\rho)}
\le e^{as}\le
\frac{\max_{C_2}(-\rho)}{\min_{C_1}(-\rho)}.
\]
Hence $\delta_{as}$ is proper for $a\ne0$. The flow
$\tau_{bs}$ is proper for $b\ne0$ by translation of $\Im w$.

For $K>0$, put
\begin{equation}\label{ft:eq:weighted-potentials}
u_{a,b,K}(z,w)=
\begin{cases}
|z_1|^{K/(a\lambda_1)},&a>0,\\
(-\rho(z,w))^{K/a},&a<0,\\
e^{(K/b)\Im w},&a=0,\ b\ne0.
\end{cases}
\end{equation}
The character identities are
\begin{equation}\label{ft:eq:weighted-potential-character}
\begin{aligned}
u_{a,b,K}\circ\delta_{as}&=e^{Ks}u_{a,b,K} &&(a\ne0),\\
u_{0,b,K}\circ\tau_{bs}&=e^{Ks}u_{0,b,K} &&(b\ne0).
\end{aligned}
\end{equation}
For $a>0$, with $\gamma=K/(a\lambda_1)>0$,
$|z_1|^\gamma$ is plurisubharmonic and
\begin{equation}\label{ft:eq:weighted-positive-levi}
\ddbar|z_1|^\gamma=
\frac{\gamma^2}{4}|z_1|^{\gamma-2}\,
\ii\,dz_1\wedge d\bar z_1
\qquad(z_1\ne0).
\end{equation}
The coefficient is locally integrable across $z_1=0$,
since $\int_0^\epsilon r^{\gamma-1}\,dr<\infty$.
The positive closed current $\ddbar|z_1|^\gamma$ is
therefore nonzero.

For $a<0$, put $q=-K/a>0$. Since $\rho<0$ and
$\ddbar\rho=\ddbar P\ge0$,
\begin{equation}\label{ft:eq:negative-power-levi}
\begin{aligned}
\ddbar(-\rho)^{-q}
&=q(-\rho)^{-q-1}\ddbar\rho\\
&\quad+q(q+1)(-\rho)^{-q-2}
\ii\,\partial\rho\wedge\bar\partial\rho\ge0.
\end{aligned}
\end{equation}
The second term is nonzero everywhere because
$\partial\rho/\partial w=1/2$. For $a=0$ and $b\ne0$,
\begin{equation}\label{ft:eq:weighted-translation-levi}
\ddbar e^{(K/b)\Im w}
=\frac{K^2}{4b^2}e^{(K/b)\Im w}\,
\ii\,dw\wedge d\bar w.
\end{equation}
Thus $\ddbar u_{a,b,K}$ is positive, closed, and nonzero
in all three cases.

For $a\ne0$, set
$T=(\tau_{b/a}\circ F)^*(\ddbar u_{a,b,K})$;
for $a=0$, set $T=F^*(\ddbar u_{0,b,K})$.
By \eqref{ft:eq:weighted-conjugation} and
\eqref{ft:eq:weighted-potential-character},
\[
\Phi_s^*T=e^{Ks}T,
\]
where $\Phi_s$ is the flow of $F^*(aD+bA)$.
Hence $0\ne T\in\mathcal T_{F^*(aD+bA),K}(\Omega)$.
Since $P$ is plurisubharmonic, $M_P$ is a
pseudoconvex domain in $\C^n$, and is Stein
\cite{Hormander1990}. By \cite[Theorem~1.2]{GangulyKRsol2026Paper21},
$F^*(aD+bA)$ is the field of a complete nontrivial expander
with expansion constant $K$.

Finally, the map
\begin{equation}\label{ft:eq:weighted-diffeomorphism}
\C^{n-1}\times(-\infty,0)\times\R\longrightarrow M_P,
\quad(z,r,t)\longmapsto(z,r-P(z,\bar z)+\ii t)
\end{equation}
is a diffeomorphism.
\end{proof}

\begin{corollary}\label{ft:cor:zimmer}
Let $\Omega\Subset\C^n$, $n\ge2$, be bounded and
pseudoconvex with $C^\infty$ boundary of finite D'Angelo
type. If the boundary orbit accumulation set
$\mathcal L(\Omega)$ of $\Aut(\Omega)$ contains two distinct
points, then $\Omega$ is biholomorphic to $M_P$ as in
\eqref{ft:eq:weighted-model}--\eqref{ft:eq:weighted} and is
diffeomorphic to $\R^{2n}$. If
$F:\Omega\to M_P$ is such a biholomorphism, then
$F^*(aD+bA)$ is the field of a complete nontrivial expander
for every $(a,b)\in\R^2\setminus\{(0,0)\}$ and $K>0$, with
$D,A$ the generators of the flows in \eqref{ft:eq:weighted-flows}.
\end{corollary}

\begin{proof}
By \cite[Theorem~1.2(1)]{Zimmer2020}, $\Omega$ is
biholomorphic to a weighted homogeneous model
$\{(z,\zeta):\Im\zeta>P(z,\bar z)\}$,
with the weights of \cite[Definition~1.1]{Zimmer2020}.
The change of coordinates $w=\ii\zeta$ gives
\eqref{ft:eq:weighted-model} and \eqref{ft:eq:weighted}.

Pseudoconvexity of $\Omega$ implies pseudoconvexity
of $M_P$. Fix $z,v\in\C^{n-1}$ and choose
$(z,w)\in\partial M_P$. Since
$V=(v,-2\partial P_z(v))$ is complex tangent to
$\partial M_P$,
\begin{equation}\label{ft:eq:zimmer-levi}
0\le-\ii(\ddbar\rho)(V,\bar V)
=\sum_{j,k=1}^{n-1}P_{j\bar k}(z)v_j\bar v_k.
\end{equation}
Thus $P$ is plurisubharmonic. Proposition~\ref{ft:prop:weighted}
gives the diffeomorphism and the expanders with fields
$F^*(aD+bA)$ for every $(a,b)\ne(0,0)$ and $K>0$.
\end{proof}

\begin{corollary}\label{ft:cor:h-extendible}
Let $\Omega\subset\C^{n+1}$, $n\ge1$, be a pseudoconvex
domain with $C^\infty$ boundary, and let
$p\in\partial\Omega$ be $h$-extendible of finite Catlin
multitype $(1,m_1,\ldots,m_n)$. Suppose that
$x_j=F_j(q)\to p$, where $q\in\Omega$,
$F_j\in\Aut(\Omega)$, and for some $C>0$,
\begin{equation}\label{ft:eq:rz-cone}
|x_j-p|\le C\,\dist(x_j,\partial\Omega)
\qquad(j\gg1).
\end{equation}
Then $\Omega$ is biholomorphic to
\[
M_P=\{(z,w)\in\C^n\times\C:\Re w+P(z,\bar z)<0\},
\]
where $P$ is a real-valued plurisubharmonic polynomial
without pluriharmonic monomials, satisfying
\[
P(e^{t/m_1}z_1,\ldots,e^{t/m_n}z_n)=e^tP(z).
\]
In particular, $\Omega$ is diffeomorphic to
$\R^{2n+2}$. If $\Omega$ is bounded, then for any biholomorphism $F:\Omega\to M_P$, every
$F^*(aD+bA)$ with $(a,b)\ne(0,0)$ is the field
of a complete nontrivial expander for every $K>0$, with
$D,A$ the generators of the flows in \eqref{ft:eq:weighted-flows}.
\end{corollary}

\begin{proof}
Let $\chi$ be the distinguished coordinate map at $p$ from
\cite[Section~3.2]{ThuDieu2020}, with $\chi(p)=0$ and
$\chi(x_j)=(\alpha_j,\beta_j)$. Write
$d_j=\dist(x_j,\partial\Omega)$, and let $\delta_j$ be the
Euclidean distance from $\chi(x_j)$ to the image of the
boundary in this chart. Comparing with $p$ and with $0$ shows
that nearest boundary points in the two coordinate systems tend
to $p$ and $0$, respectively. The chart and its inverse are
Lipschitz on smaller neighborhoods of these points, so
$c d_j\le\delta_j\le C d_j$ for large $j$, with $c,C>0$.

Since $\chi(p)=0$, the Lipschitz bound and \eqref{ft:eq:rz-cone}
give $|\alpha_j|+|\beta_j|=O(d_j)=O(\delta_j)$.
Since $m_k\ge2$ and $\alpha_j\to0$,
\begin{equation}\label{ft:eq:weighted-nontangential}
|\Im\beta_j|+\sum_{k=1}^n|\alpha_{j,k}|^{m_k}
=O(\delta_j).
\end{equation}
With $\Lambda=(1/m_1,\ldots,1/m_n)$, this is the
$\Lambda$-nontangential condition of
\cite[Definition~3.4]{ThuDieu2020}. By
\cite[Theorem~1.1]{ThuDieu2020}, $\Omega$ is biholomorphic to
$M_P$ as stated.

The coordinates $(z,\Re w+P(z,\bar z),\Im w)$
give $M_P\simeq\R^{2n+2}$ diffeomorphically.
For bounded $\Omega$,
Proposition~\ref{ft:prop:weighted} gives the
prescribed complete nontrivial expanders.
\end{proof}

\section*{Acknowledgments}
The author thanks Max Hallgren and Yuan Yuan
for helpful discussions concerning the subject
of this paper. The author is also grateful to
Peter Ebenfelt and Ming Xiao for their
encouragement and support.

\section*{AI-use declaration}
The author used OpenAI's ChatGPT to assist with
locating relevant literature, translating passages
from non-English mathematical sources, and improving
the grammar, organization, and exposition of the
manuscript. The author takes full responsibility
for the contents of this paper.

\bibliographystyle{amsplain}
\bibliography{refKRsoliton}
\end{document}